\documentclass[review, 3p, times]{elsarticle}
\usepackage{amsmath}
\usepackage{amsfonts}
\usepackage{amssymb}
\usepackage{amsthm}
\usepackage{bm}
\usepackage{indentfirst}
\usepackage{graphicx}
\usepackage{subfigure}
\usepackage{array}
\usepackage{fullpage}
\usepackage{mathrsfs} 

\DeclareMathAlphabet{\mathsfsl}{OT1}{cmss}{m}{sl}

\newcommand{\PreserveBackslash}[1]{\let\temp=\\#1\let\\=\temp}
\newcolumntype{C}[1]{>{\PreserveBackslash\centering}p{#1}}
\newcolumntype{R}[1]{>{\PreserveBackslash\raggedleft}p{#1}}
\newcolumntype{L}[1]{>{\PreserveBackslash\raggedright}p{#1}}

\numberwithin{equation}{section}
\newtheorem{thm}{Theorem}[section]
\newtheorem{lem}[thm]{Lemma}
\newtheorem{cor}[thm]{Corollary}
\theoremstyle{definition}

\newtheorem{pro}[thm]{Proposition}

\def\cS{\mathcal{S}}
\def\cT{\mathcal{T}}
\def\cW{\mathcal{W}}
\def\sQ{\mathscr{Q}}
\def\sM{\mathscr{M}}
\def\sN{\mathscr{N}}
\def\sW{\mathscr{W}}
\def\R{\mathbb{R}}
\def\Z{\mathbb{Z}}

\newcommand{\del}{\delta}
\newcommand{\ga}{{\gamma}}
\newcommand{\Om}{\Omega}
\newcommand{\by}{\bm y}
\newcommand{\bx}{\bm x}
\newcommand{\bz}{\bm z}

\newcommand{\dist}{\text{dist}}

\usepackage{hyperref}
\usepackage{fancyvrb}
\usepackage{fancyhdr}
\usepackage{color}
\usepackage{tikz}
\usepackage{appendix}

\newcommand{\real}{\mathbb{R}}

\newcommand{\mcL}{\mathcal{L}}

\newcommand{\mcR}{\mathcal{R}}

\def\omg{{\Omega}}

\def \tx{\tilde{x}}
\def \ty{\tilde{y}}
\def \tz{\tilde{z}}
\def \tw{\tilde{w}}

\def \tl{\tilde{l}}

\def \wb{\bm{w}}

\def \xb{\bm{x}}

\def \hb{\bm{h}}

\def \zb{\bm{z}}

\def \yb{\bm{y}}

\def \gbd{{\ga^{\,\beta}_\del}}

\newcommand{\vertii}[1]{{\left\vert\left\vert #1
    \right\vert\right\vert}}

\newcommand{\verti}[1]{{\left\vert #1
    \right\vert}}

\usepackage{xcolor}
\definecolor{darkgreen}{RGB}{21, 233, 153}

\biboptions{sort&compress} 

\begin{document}

\theoremstyle{definition}
\newtheorem{remark}{Remark}

\begin{frontmatter}

\title{Nonlocal Trace Spaces and Extension Results for Nonlocal Calculus}

\address[qd]{Department of Applied Physics and Applied Mathematics and Data Science Institute, Columbia University, New York, NY, 10027}
\address[xt]{Department of Mathematics, The University of California, San Diego, CA, 92093}
\address[yy]{Department of Mathematics, Lehigh University, Bethlehem, PA 18015, USA}

\author[qd]{Qiang Du}\ead{qd2125@columbia.edu}
\author[xt]{Xiaochuan Tian}\ead{xctian@ucsd.edu}
\author[yy]{Cory Wright}\ead{cow218@lehigh.edu}
\author[yy]{Yue Yu\corref{cor1}}\ead{yuy214@lehigh.edu}

\begin{abstract}
For a given Lipschitz domain {$\Omega$}, it is a classical result that the trace space of $W^{1,p}(\Om)$ is $W^{1-1/p,p}(\partial\Om)$, namely  any $W^{1,p}(\Om)$ function has a well-defined $W^{1-1/p,p}(\partial\Om)$ trace on its codimension-1 boundary $\partial\omg$ and any $W^{1-1/p,p}(\partial\Om)$ function on $\partial\omg$ can be extended to a $W^{1,p}(\Om)$ function. Recently, \cite{dyda2019function} characterizes the trace space for nonlocal Dirichlet problems  involving integrodifferential operators with infinite interaction ranges, where the boundary datum is provided on the whole complement of the given domain $\real^d\backslash\omg$. In this work, we study function spaces for nonlocal Dirichlet problems with a finite range of nonlocal interactions, which naturally serves a bridging role between the classical local PDE problem and the nonlocal problem with infinite interaction ranges. For these nonlocal {Dirichlet} problems, the boundary conditions are normally imposed on a region with finite thickness volume which lies outside of the domain. We introduce a function space on the volumetric boundary region that serves as a trace space for these nonlocal problems and study the related extension results. Moreover, we discuss the consistency of the new nonlocal trace space with the classical $W^{1-1/p,p}(\partial\Om)$ space as the size of nonlocal interaction tends to zero. In making this connection, we conduct an investigation on the relations between nonlocal interactions on a larger domain and the induced interactions on its subdomain.  The various forms of trace, embedding and extension theorems may then be viewed as consequences in different scaling limits.
\end{abstract}
\begin{keyword}
Trace Theorem; Nonlocal Function Space; Inverse Trace Theorem; Trace Space; Finite Nonlocal Interactions; Extension Operator.\\
{\MSC[2010] 	46E35   \sep  47G10 \sep 35A23 \sep 35R11}
\end{keyword}

\end{frontmatter}

\tableofcontents

\section{Introduction}
 In \cite{gagliardo1957}, Gagliardo characterizes the trace space of the Sobolev space $W^{1,p}(\Om)$ ($p>1$) for a given {bounded} Lipschitz domain $\Om\subset \R^d$. The result consists of the following two parts.  First, the trace operator $T$ from $W^{1,p}(\Om)$ to $W^{1-1/p,p}(\partial\Om)$ is linear and continuous, and conversely, one can define a continuous linear extension operator $E$ from $W^{1-1/p,p}(\partial\Om)$ to $W^{1,p}(\Om)$.
Our goal in this work is to study the trace spaces of {some} nonlocal function spaces denoted by $\{\cS^{\,\beta}_\del(\hat\Om)\}_{\del>0}$ 
related to the Dirichlet energies of a class of nonlocal problems. {Here, $\del>0$ represents the horizon parameter that characterizes the ranges of nonlocal interactions, and
$\hat\Om=\Om\cup\Om_{\del}$
with $\Om_{\del}:= \left\{ \bm x\in \R^d\backslash\Om:  \text{dist}(\bm x, \partial\Om)  <\del\right\}$ being viewed as a nonlocal ``boundary'' set of the given  domain $\Om$.}
The function space $\cS^{\,\beta}_\del(\hat\Om)$ is defined as 
the completion of $C^1(\overline{\hat\Om})$ with respect to the norm 
\begin{align}
\label{eq:nnorm}
\| \cdot \|_{\cS^{\,\beta}_\del(\hat{\Om})} = (\| \cdot \|_{L^p(\hat{\Om})}^p + 
 | \cdot |^p_{\cS^{\,\beta}_\del(\hat{\Om})} )^{1/p},
\end{align} 
with the associated semi-norm $|\cdot |_{\cS^{\,\beta}_\del(\hat\Om)}$ given by
\begin{align}
\label{eq:nsnorm}
| u |^p_{\cS^{\,\beta}_\del(\hat\Om)}= \int_{\hat\Om}  \int_{\hat\Om} \gbd(| \bm y-\bm x|) |u(\bm y)-u(\bm x)|^p d\bm y d\bm x \,. 
\end{align}
Notice that the space $\cS^{\,\beta}_\del$ also depends on $p$, and we always assume $p>1$ in this paper.
The kernel function $\gbd$ in \eqref{eq:nsnorm},
for $\del>0$ and $\beta\in [0,d+p)$ is taken as 
 \begin{equation}\label{def:rescaledkernel}
\gbd(|\bm y - \bm x|) =\frac{C_{d,p,\,\beta}}{\delta^{d+p-\beta}}\frac{1}{|\yb-\xb|^{\,\beta}} 1_{\{|\yb-\xb|<\delta\}},
  \end{equation}
where for any $\del$, $1_{\{|\yb-\xb|<\delta\}}$ denotes the characteristic function on the set $\{|\yb-\xb|<\delta\}$, and $C_{d,p,\,\beta}$ normalizes the $p^{th}$ moment of $\gamma_\delta^{\,\beta}$.  In particular, we have 
\[
C_{d,p,\,\beta}=s^{-1}_{d-1}(d+p-\beta),
\]
where $s_{d-1}$ denotes the area of the $d-1$-sphere since
\[
\int_{\real^d}\frac{1}{\delta^{d+p-\beta}}\frac{1}{|\bz|^{\,\beta}} 1_{\{|\bz|<\delta\}}|\bz|^pd\bz=\frac{1}{\delta^{d+p-\beta}}\int_{B(\bm{0},\delta)} |\bz|^{p-\beta}d\bz=\frac{s_{d-1}}{\delta^{d+p-\beta}}\int_0^\delta r^{d+p-\beta-1}dr=\frac{s_{d-1}}{d+p-\beta}.
\]
We note that, for different $\del>0$, $\gbd$ can be obtained from a rescaling of a given $\del$-independent nonnegative kernel $\ga^{\,\beta}$ defined on $(0,1)$ by:
\begin{equation} \label{def:frackernel}
\gbd(|\bm y - \bm x|) =\frac{1}{\delta^{d+p}}\ga^{\,\beta}\left(\frac{|\bm y - \bm x|}{\delta} \right)\,, \quad\text{where }\ga^{\,\beta}(|\bm y - \bm x|) = \frac{C_{d,p,\,\beta}}{|\bm y - \bm x|^{\,\beta}} 1_{\{|\bm y - \bm x|<1\}}.
\end{equation}

{It is easy to see that the nonlocal function space $\cS^{\,\beta}_\del(\hat\Om)$ contains all $p$-integrable functions on $\hat\omg$ with finite norms with respect to
$\| \cdot \|_{\cS^{\,\beta}_\del(\hat{\Om})}$. Moreover, for any finite and given $\del>0$, the kernel $\gbd$ is integrable for
 $\beta\in [0,d)$ and the corresponding space $\cS^{\,\beta}_\del(\hat\Om)$ is equivalent to the $L^p(\hat\Om)$ space; while $\gbd$ is non-integrable  for $\beta\in (d, d+p)$ and  $\cS^{\,\beta}_\del(\hat\Om)$ is equivalent to the standard fractional Sobolev space $W^{(\beta-d)/p, p}(\hat\Om)$}. 
 Moreover, we have the convergence of the space $\cS^{\,\beta}_\del(\hat\Om)$ to 
{the local limit} $W^{1,p}(\Om)$ as $\del\to0$, see e.g., discussions in \cite{BBM01,ponce2004,Mengesha12,foss2016differentiability}.  

{Given any domain $\Om$,} let  $\cT^{\,\beta}_\del({\Om})$ denote the space of all $L^p({\Om})$-functions $u$ with the norm defined as
\begin{equation} 
\|u\|_{\cT^{\,\beta}_\del({\Om})}:= \left( \frac{1}{\del} \| u\|^p_{L^p({\Om})} + |u|^p_{\cT^{\,\beta}_\del({\Om})} \right)^{1/p}\,.
\end{equation}
Here the semi-norm is defined as
\begin{equation}\label{eqn:nonlocalhalf}
\left|u\right|_{\cT^{\,\beta}_\del({\Om})}:= \left( \del^{\,\beta-2}\int_{{\Om}}\int_{{\Om}} \frac{|u(\bm y) -u(\bm x)|^p}{(|\bm y-\bm x|\vee \del)^{d+p-2} (|\bm y -\bm x|\wedge\del)^{\,\beta}} \,d\yb d\xb \right)^{1/p},     \end{equation}
where $a\wedge b:=\min(a, b)$ and $a\vee b:=\max(a, b)$. Our main result is to show that $\cT^{\,\beta}_\del(\Om_{\del})$ is the trace space of $\cS^{\,\beta}_\del(\hat\Om)$, i.e., {to establish} the existence of trace operator $T$ and extension operator $E$ that define continuous linear maps in between $\cS^{\,\beta}_\del(\hat\Om)$ and $\cT^{\,\beta}_\del(\Om_{\del})$.

For the ease of presentations, in the following we denote:
\begin{equation}
    \label{stripenotation}
\mathcal{R}^L=(0,L)\times\real^{d-1},\, \mathcal{R}_a=(a,0)\times\real^{d-1}, \text{ for } a<0, \, \text{ and } \mathcal{R}_a^L=(a,L)\times\R^{d-1}, \text{ for } a < L.
\end{equation}
and we note that when $\omg=(0,\infty)\times\real^{d-1}=\mcR^{\infty}$, we have $\omg_{\delta}=\mathcal{R}_{-\del}$ and $\hat{\omg}=\mcR_{-\del}^{\infty}$.  {With these notations, 
we first show the main results on half spaces as follows.}

\begin{thm}[Trace theorem on half spaces]\label{mainthm_1}
Let $\del>0$ and $\beta\in [0, d+p)$, then there exists a constant $C$ independent of $\delta$ and $\beta$
such that for any $u\in \cS^{\,\beta}_\del(\mcR_{-\del}^{\infty})$,
\[
\frac{1}{\del}\|  u\|^p_{L^p(\mcR_{-\del})}\leq C |d+p-\beta|^{-1}\| u \|^{p-1}_{L^p(\mcR_{-\del}^{\infty})}\verti{u}_{\cS^{\,\beta}_\del(\mcR_{-\del}^{\infty})}\leq C |d+p-\beta|^{-1} \left(\| u \|^p_{L^p(\mcR_{-\del}^{\infty})}+\verti{u}^p_{\cS^{\,\beta}_\del(\mcR_{-\del}^{\infty})}\right),
\]
\[
\verti{u}^p_{\cT^{\,\beta}_\del(\mcR_{-\del})}\leq C |d+p-\beta|^{-1}\verti{u}^p_{\cS^{\,\beta}_\del(\mcR_{-\del}^{\infty})},
\]
and therefore
\[
\|  u\|_{\cT^{\,\beta}_\del(\mcR_{-\del})}\leq C |d+p-\beta|^{-1/p}\| u \|_{\cS^{\,\beta}_\del(\mcR_{-\del}^{\infty})}.
\]
\end{thm}


\begin{thm}[Inverse trace theorem on  half spaces]\label{mainthm_2}
Let {$\del\in (0,M)$ for some fixed number $M>0$}, {$\beta\in [0, d+p)$, and $\beta\neq d$,} then there exists an extension operator $E: \cT^{\,\beta}_\del(\mcR_{-\del})\to \cS^{\,\beta}_\del(\mcR_{-\del}^\infty)$ such that
\[
\| E u \|_{\cS^{\,\beta}_\del(\mcR_{-\del}^\infty)} \leq C |d-\beta|^{-1/p} \| u\|_{\cT^{\,\beta}_\del(\mcR_{-\del})}\,, {\quad \forall u \in 
\cT^{\,\beta}_\del(\mcR_{-\del})}
\]
where $C$ is a constant independent of $\del$, $\beta$ and $u\in \cT^{\,\beta}_\del(\mcR_{-\del})$.
\end{thm}



{Using partition of unity techniques, the above trace theorems in special domains can then be
extended to more general domains. which are stated in the theorems below.}
 
\begin{thm}[General trace and inverse trace theorems]\label{mainthm_1_general}
Assume that $\Om$ is a bounded {and} simply connected Lipschitz domain in $\R^d$ and $\Om_\del:=\{\xb\in\real^d\backslash\omg: \text{dist}(\xb,\omg)<\del\}$ is its nonlocal boundary set.
There exists a constant $\epsilon$ depending on the domain $\Om$, such that for any $\del\in(0,\epsilon)$ and $\beta$,
\[
\| u\|_{\cT^{\,\beta}_\del(\Om_\del)}\leq C_1|d+p-\beta|^{-1/p} \| u \|_{\cS^{\,\beta}_\del(\hat\Om)}, {\quad \forall u \in \cS^{\,\beta}_\del(\hat\Om),} 
\]
On the other hand, for any $\del\in (0,\epsilon)$, $\beta\in[0,d+p)$ {and $\beta\neq d$}, there exists an extension operator $E: \cT^{\,\beta}_\del(\Om_\del)\to \cS^{\,\beta}_\del(\hat\Om)$ such that
\[
\| E u \|_{\cS^{\,\beta}_\del(\hat\Om)} \leq C_2 |d-\beta|^{-1/p} \left( \| u \|_{\cT_\del^{\,\beta}(\Om_\del)} + |d+p-\beta|^{-1/p}\| u\|_{L^p(\Om_\del)}\right)\,, {\quad \forall u \in 
\cT^{\,\beta}_\del(\Om_\del).}
\]
Here $C_1$, $C_2$ are constants independent of $\del$, $\beta$.
\end{thm}

The paper is organized as follows. In Section \ref{sec:notation} we discuss the motivation of this work, together with {additional definitions  and notation} relevant to our main results. To provide some insights {on the various nonlocal spaces under consideration,  
we also investigate their scaling properties and consistency with the
classical trace spaces in the local limit, 
as $\del\rightarrow 0$.
 Moreover,
in making these connections, our study also represents an investigation on the relations between nonlocal interactions on a larger domain and the induced interactions on a subdomain of a smaller size or dimension. This further leads to a new way of viewing the
various forms of trace, embedding and extension theorems in different function spaces as consequences in different scaling limits, further illustrating the contributions of our study.}
Sections \ref{sec:traceRd}-\ref{sec:generaldomain} contain the proofs of the aforementioned trace theorems. In particular, to show these trace theorems, {while following} the footsteps of the proofs for the local trace theorems, {we take into account the effect of nonlocal interactions}.  In Section \ref{sec:traceRd} we provide the proof of Theorem \ref{mainthm_1} with a special case $\del=1$ first, which captures the intrinsic effect of nonlocal interactions defined on a larger domain for subdomains. We then extend the results to cases with general $\delta>0$ using a scaling argument.  For the inverse trace theorem, in Section \ref{sec:inversetrace} we present the proof of Theorem \ref{mainthm_2} by constructing an extension operator based on the Whitney decomposition. Then, in Section \ref{sec:generaldomain} we prove Theorem \ref{mainthm_1_general} for the general bounded simply connected Lipschitz domain $\omg$ {using partition of unity techniques}. 
Lastly, Section \ref{sec:conclusion} summarizes our findings and discusses future research directions.

\section{Motivation and Notation}\label{sec:notation}

{In this section, we first make a few comments on the motivation of our work. We then 
investigate the consistency and connection between our nonlocal trace space with the trace space in the classical calculus in Section \ref{sec:loclimit}. We also}
provide notations and several useful lemmas for the later proofs, including some scaling properties in Section \ref{sec:scale} and the Whitney decomposition of $\mcR^L$ in Section \ref{sec:whitney}.

\subsection{
Local, nonlocal and fractional modeling}\label{eqn:nonlocalcalc}
{A major motivation of our work comes from nonlocal modeling that are represented by integro-differential equations, in particular, equations involving nonlocal interactions with a finite interaction length. The latter have drawn much attention recently in modelling certain physical systems where the classical models are not most effective}. Comparing with the classical {local} partial differential equation (PDE) models, these equations have the ability to describe these physical phenomena in a setting with reduced regularity requirements allowing singularities and discontinuities to naturally occur \cite{du2013nonlocal,gunzburger2010nonlocal,foss2018existence}. On the other hand,  when comparing with the nonlocal integro-differential equations characterized by an infinite lengthscale, compactly supported nonlocal models are computationally more efficient and therefore a more feasible choice for scientific and engineering applications. These extra flexibility and efficiency allow this framework to be used in many different situations involving physical discontinuity such as dynamic fracture \cite{silling_2000,silling2007peridynamic,hu2012peridynamic,Ha2010peri,CHENG2015,ZHANG2018Rayleighcrack,Yu2018paper,Trask2018paper}, corrosion models \cite{CHEN2015pitcorrosion, Chen2015passivefilm,Jafarzadeh2017corrosion,Li2017corrosion}, and heat conduction \cite{Bobaru2010heatconduction}. The development in this subject has also produced other applications in image processing \cite{lou2010image} and population models \cite{Carrillo2005} among many other different fields which can be further seen in \cite{bobaru2015handbook}.  Particularly, nonlocal problems with boundary constraints have become of recent interest in works such as \cite{le2018surface,bobaru2016handbook,oterkus2010peridynamic,macek2007peridynamics,du2017peridynamic,madenci2014peridynamic,oterkus2014peridynamic,oterkus2015peridynamics,tao2017nonlocal,you2019neumann,you2020asymptotically,yu2021asymptotically}.  In nonlocal models, the boundary conditions are normally not imposed on a sharp interface.  Rather, they are imposed on a region with non-zero volume which lies outside of the domain, and treating the nonlocal boundary problem improperly can cause artificial phenomena such as a ``surface'' or ``skin'' effect \cite{bobaru2011adaptive,ha2011characteristics,prudhomme2020treatment,chen2020peridynamics}.  Differs from the local problems, in some nonlocal problems boundary effects play a major role. For example, in nonlocal minimal surface problems, the ``stickiness'' effect arises and the boundary datum may not be attained continuously \cite{dipierro2017boundary,borthagaray2019finite}. All the above examples indicate that studying the nonlocal boundary conditions and the {associated nonlocal trace spaces} are critical for the development of nonlocal models.

In this work, we aim to introduce {a function space that serves as a trace space} for nonlocal problems with constant finite interaction length (the so-called interaction radius or horizon $\delta$), and study related extension results. Extension and trace theorems are well-known in the study of classical local problems with boundary constraints. For the case of Sobolev spaces of integer order, these results are well-established long time ago (see, e.g., \cite{adams2003sobolev,slobodetskiui1958sl}). For Sobolev spaces with fractional order of differentiability, which can be seen as one type of nonlocal problems with infinite interaction length, the trace space and extension results are studied in \cite{koskela2017traces,dyda2019function,bogdan2020extension,rutkowskifunction}. The latter can be useful in studying nonlocal problems with non-homongeneous boundary data, such as those associated with the nonlocal Laplacian and nonlocal $p$-Laplacian, see for example \cite{Andreu2008,Andreu2009,Andreu2010,bogdan2020extension,DRV2017,Ros2016}. In \cite{du2021fractional,tian2017trace,foss2020traces}, trace theorems are developed for nonlocal problems with varying influence horizon $\delta(\xb)$, where $\delta(\xb)\rightarrow0$ as $\xb$ approaches the boundary, in a way that the trace spaces of classical Sobolev spaces are recovered. The trace results are also applied to the study of the coupling of nonlocal and local models \cite{TaTiDu19}. To our best knowledge, the definition of trace space and extension results for nonlocal problems with constant finite horizon have not been dealt with so far. These results would extend the knowledge on the trace space in nonlocal calculus and its connection with the trace space in classical calculus. Moreover, the trace theorem and the inverse trace theorem would also provide important mathematical tools for developing well-posed nonlocal models with volumetric boundary conditions, such as discussed in \cite{you2019neumann}. 

\subsection{Nonlocal Space $\cS^{\,\beta}_\del(\Om)$, Associated Nonlocal Problems and Their Local Limits}\label{eqn:nonlocalcalcf}


Before discussing their connections in the following sections, in this section we introduce the classical and nonlocal Laplacian operators and their corresponding nonlocal function spaces relevant to this paper. {The discussions in this subsection are restricted to the Hilbert space setting where $p=2$.}

{Given a scalar function $u(\xb):\omg\rightarrow\real$, the classical Laplacian operator is defined as $\Delta u:=\nabla\cdot\nabla u$ and
boundary value problems on the domain $\omg$ related to $\Delta$ are often associated with the Sobolev space $H^1(\omg)$ with its} norm defined by
\[\vertii{u}_{H^1(\omg)}:=\left(\vertii{u}^2_{L^2(\omg)}+\verti{u}^2_{H^1(\omg)}\right)^{1/2}.\]

On the other hand, when incorporating long-range interactions into the model such that where every point $\xb\in\omg$ is interacting with a finite neighborhood of points, 
a nonlocal Laplacian operator 
is then given by
$$\mcL[u](\xb):=C\int_{\hat{\omg}}\gamma(\xb,\yb)(u(\yb)-u(\xb))d\yb, \quad \xb\in\omg,$$
where {$\gamma(\xb,\yb)$ is a kernel function that will be prescribed shortly}, $\hat{\omg}=\omg\cup\omg_I$ and
$$\omg_I:=\{\yb\in\real^d\backslash\omg\text{ such that }\gamma(\xb,\yb)\neq0 \text{ for some }\xb\in\omg\}$$
is the interaction domain of $\omg$. The nonlocal Laplacian operator is associated with the following nonlocal norm
$$\vertii{u}_{\cS(\hat{\omg})}:=\left(\vertii{u}^2_{L^2(\hat{\omg})}+\verti{u}^2_{\cS(\hat{\omg})}\right)^{1/2}\text{ where }\verti{u}^2_{\cS(\hat{\omg})}:= \frac{C}{2}\int_{\hat{\omg}}  \int_{\hat{\omg}} \gamma(\xb,\yb) (u(\bm y)-u(\bm x))^2 d\bm y d\bm x.$$
In this paper we further assume that such neighborhood is a Euclidean ball surrounding $\xb$, i.e., $B(\xb,\delta):=\{\yb\in\real^d:|\yb-\xb|<\delta\}$. Here $\delta$ is the interaction radius or horizon. This fact has implications on the boundary conditions that are prescribed on a collar of thickness $\delta$ outside the domain $\omg$, that we have the interaction domain $\omg_I=\omg_\del {:= \{\yb\in\R^d\backslash\Om: \dist(\yb,\partial\Om)<\del }\}$. In particular, we can take a popular class of kernels $\gamma(\xb,\yb)=\ga_{\delta}^{\,\beta}(\verti{\yb-\xb})$ as in \eqref{def:rescaledkernel}. We note that when the constant $C=C^{\text{diff}}=2d$, we have the following property  
\begin{equation}\label{eqn:Cdiff}
C^{\text{diff}}\int_{B(\xb,\delta)}\gamma^{\,\beta}_\delta(|\yb-\xb|)|\yb-\xb|^2 d\yb=2d \quad {\forall \xb \in \real^d},
\end{equation}
since the $p^{th}$ moment of kernel $\gamma_\del^{\,\beta}$ in \eqref{def:rescaledkernel} is normalized to $1$. Then
it is well-known (see, e.g., 
\cite{du2013nonlocal}) that the nonlocal diffusion operator converge to its local for all  counterpart pointwise: for any $u\in C^\infty(\real^d)$ and $\xb\in\real^d$,
$$\mcL[u](\xb)=C^{\text{diff}}\int_{B(\xb,\delta)}(u(\yb)-u(\xb))\gamma^{\,\beta}_\del(|\yb-\xb|)d\yb\overset{\delta\rightarrow 0}{\longrightarrow} \Delta u(\xb).$$
Moreover, when $u\in H^1(\omg)$, its nonlocal norm converges to the $H^1$ norm:
$$\| u \|_{\cS^{\,\beta}_\del(\omg)}\overset{\delta\rightarrow 0}{\longrightarrow} \| u \|_{H^1(\omg)}.$$
Naturally, we can extend the above conclusion to more general cases of nonlocal and local $p$-Laplacians corresponding to $p>1$.

\subsection{{Nonlocal Space $\cT^{\,\beta}_\del(\Om_\del)$ and} Connection to Classical Local Trace Spaces}\label{sec:loclimit}

First of all, we may view \eqref{eqn:nonlocalhalf} as a nonlocal counterpart of the classical trace semi-norm 
$$|u|_{W^{1-1/p, p}(\partial\omg)}:=\left( \int_{\partial\omg}\int_{\partial\omg}\dfrac{|u(\yb)-u(\xb)|^p}{|\yb-\xb|^{d+p-2}}d\yb d\xb \right)^{1/p}$$
and seek a nonlocal analog of the classical trace theorem. The relation between the classical $W^{1-1/p, p}(\partial\Omega)$ trace space and the new nonlocal trace space $\cT^{\,\beta}_\del(\Om_{\del})$ can be seen from the
limiting process as $\del\to0$ in the following proposition. In the rest of the paper,
we use $f\lesssim g$ if $f\leq C g$ for a generic constant $C>0$  independent of $\del$ and $\beta$.
We also write $f\approx g$ if $f\lesssim g$ and $g\lesssim f$. 

\begin{pro}\label{thm:consist}
Let $\partial\mcR=\{0\}\times\R^{d-1}$ and $\mcR_{-\del}$ 
{be defined as in \eqref{stripenotation}}
for $\del\in (0,1)$, then
$$| u |^p_{\cT^{\,\beta}_\del(\mcR_{-\del})} \xrightarrow{\del\to0} | u |^p_{W^{1-1/p, p}(\partial\mcR)},$$
for any {$u\in C^1_c\left([-1,0]\times B^{d-1}(\bm 0,M)\right)$  for some $M>0$. Here $B^{d-1}(\overline{\xb},r)$ denotes the ball centered at $\overline{\xb}$ with radius $r$ in $\real^{d-1}$.}
\end{pro}
\begin{proof}
In this proof, {we denote any point $\bx\in \R^d$}  by ${\bx=}(\tilde{x},\overline{\bx})\in \R\times \R^{d-1}$. Similarly $\by\in\R^d$ is also denoted by $(\tilde{y},\overline{\by})\in\R\times\R^{d-1}$. We first have the estimate
\begin{align*}
&\left|\del^{\,\beta-2}\int_{\mcR_{-\del}}\int_{\mcR_{-\del}} \frac{|u(\bm y) -u(\bm x)|^p}{(|{\bm y}-{\bm x}|\vee\del)^{d+p-2} (|\bm y -\bm x|\wedge\del)^{\,\beta}} \,d\yb d\xb-\int_{\partial\mcR}\int_{\partial\mcR}\dfrac{|u({\yb})-u({\xb})|^p}{|{\yb}-{\xb}|^{d+p-2}}d{\yb} d{\xb}\right|\\
=&\left|\del^{\,-2}\int_{\mcR_{-\del}} \int_{{\mcR_{-\del}}\backslash B(\xb,\delta)}\frac{|u(\yb) -u(\xb)|^p}{|{\bm y}-{\bm x}|^{d+p-2}} \,d\yb d\xb+\del^{\,\beta-d-p}\int_{\mcR_{-\del}} \int_{{\mcR_{-\del}}\cap B(\xb,\delta)} \frac{|u(\yb) -u(\xb)|^p}{|\bm y -\bm x|^{\,\beta}} \,d\yb d\xb\right.\\
&\left.\qquad -\int_{\real^{d-1}}\int_{\real^{d-1}}\dfrac{|u(0,\overline{\yb})-u(0,\overline{\xb})|^p}{|\overline{\yb}-\overline{\xb}|^{d+p-2}}d\overline{\yb} d\overline{\xb}\right|\\
\leq&\left|\del^{\,\beta-d-p}\int_{\mcR_{-\del}} \int_{{\mcR_{-\del}}\cap B(\xb,\delta)}\frac{|u(\yb) -u(\xb)|^p}{|{\bm y}-{\bm x}|^{\,\beta}} \,d\yb d\xb-\del^{\,-d-p}\int_{\mcR_{-\del}} \int_{{\mcR_{-\del}}\cap B(\xb,\delta)}|u(\yb) -u(\xb)|^p \,d\yb d\xb\right|~~~~\leftarrow\;A_1\\
&+\left|\del^{\,-2} \int_{\mcR_{-\del}} \int_{\mcR_{-\del}} \frac{|u(\yb) -u(\xb)|^p}{(|{\bm y}-{\bm x}|\vee\del)^{d+p-2}} \,d\yb d\xb-\del^{\,-2}\int_{\real^{d-1}}\int_{-\del}^0 \int_{\real^{d-1}}\int_{-\del}^0 \frac{|u(0,\overline{\yb}) -u(0,\overline{\xb})|^p}{(|{\bm y}-{\bm x}|\vee\del)^{d+p-2}} \,d\ty d\overline{\yb} d\tx d\overline{\xb}\right|~~~~~~~~~\leftarrow\;A_2\\
&+\left|\del^{\,-2}\int_{\real^{d-1}}\int_{-\del}^0 \int_{\real^{d-1}}\int_{-\del}^0 \frac{|u(0,\overline{\yb}) -u(0,\overline{\xb})|^p}{(|{\bm y}-{\bm x}|\vee\del)^{d+p-2}} \,d\ty d\overline{\yb} d\tx d\overline{\xb}-\int_{\real^{d-1}}\int_{\real^{d-1}}\dfrac{|u(0,\overline{\yb})-u(0,\overline{\xb})|^p}{|\overline{\yb}-\overline{\xb}|^{d+p-2}}d\overline{\yb} d\overline{\xb}\right|.~~~~~~~~~~\leftarrow\;A_3
\end{align*}
To estimate the $A_1$ part, we first note that {$u\in C^1([-1,0]\times B^{d-1}(\bm 0,M))$} 
implies
\begin{equation}\label{eq:prop2.1}
 \verti{{u}(\yb) -{u}(\xb)}\leq \tilde{C} \verti{\yb-\xb}   
\end{equation}
for a constant $\tilde{C}$ independent of $\delta$, $\xb$, and $\yb$. 
Notice also that $\text{supp}({u})\subseteq [-1,0]\times B^{d-1}(\bm 0, M)$ for some $M>0$. Therefore
\begin{align*}
A_1=&\del^{\,-d-p}\int_{\mcR_{-\del}} \int_{\mcR_{-\del}\cap B(\xb,\del)}{|{u}(\yb) -{u}(\xb)|^p} \left|\dfrac{\del^{\,\beta}}{|\bm y -\bm x|^{\,\beta}}-1\right|\,d\yb d\xb\\
\leq&C \del^{\,-d-p+1} {\int_0^\del}r^{d-1+p} \left|\dfrac{\del^{\,\beta}}{r^{\,\beta}}-1\right|\,dr\leq C\del^{\,-d-p+1}\left(\dfrac{\del^{d+p}}{d+p-\beta}+ \dfrac{\delta^{d+p}}{d+p}\right)\leq C \del\xrightarrow{\del\to0} 0.
\end{align*}
For the $A_2$ part, we first want to show for any $\xb=(\tilde{x}, \overline{\xb})\in \mcR_{-\del}$ and  $\yb=(\tilde{y}, \overline{\yb}) \in \mcR_{-\del}$, we have 
\begin{equation}
\label{eq:prop2.1A2}
\begin{split}
&\left| |u(\ty,\overline{\yb}) -u(\tx,\overline{\xb})|^p-|u(0,\overline{\yb}) -u(0,\overline{\xb})|^p \right| \\
&
\left\{ 
\begin{aligned}
&\lesssim\, \max( \del \verti{\yb-\xb}^{p-1}, \del^p),\quad \text{ when }\xb,\yb\in (-\del,0)\times B^{d-1}(\bm 0,2M),\\
&\lesssim\, \del ,\quad \text{ when }\xb\in (-\del,0)\times B^{d-1}(\bm 0,M),\, \yb\in \omg_\del\backslash(-\del,0)\times B^{d-1}(\bm 0,2M) \text{ or vice verse},\\
&=0,\,\quad \text{ else}.
\end{aligned}\right.  
\end{split}
\end{equation}
To show \eqref{eq:prop2.1A2}, we can first assume $ |u(\ty,\overline{\yb}) -u(\tx,\overline{\xb})| \geq  |u(0,\overline{\yb}) -u(0,\overline{\xb})|$ and $u(0,\overline{\yb}) -u(0,\overline{\xb}) \geq 0$ without loss of generality. Then by rewriting $ u(\ty,\overline{\yb}) -u(\tx,\overline{\xb})$ as $ \left(u(\ty,\overline{\yb}) -u(0,\overline{\yb})  -( u(\tx,\overline{\xb}) -  u(0,\overline{\xb}))\right) + u(0,\overline{\yb}) -u(0,\overline{\xb}) $ and the fact that 
\[
\left| u(\ty,\overline{\yb}) -u(0,\overline{\yb})  -( u(\tx,\overline{\xb}) -  u(0,\overline{\xb}))\right|  \leq \tilde C( |\ty| + |\tx|) \leq 2 \del \tilde C ,
\]
we can estimate $u(\ty,\overline{\yb}) -u(\tx,\overline{\xb})$ by two different cases where $u(0,\overline{\yb}) -u(0,\overline{\xb}) > 4 \del \tilde C$ or $0\leq u(0,\overline{\yb}) -u(0,\overline{\xb})\leq 4 \del \tilde C$. 
If $u(0,\overline{\yb}) -u(0,\overline{\xb}) > 4 \del \tilde C$, then we must have $u(\ty,\overline{\yb}) -u(\tx,\overline{\xb}) >0 $  and therefore 
\[
\begin{split}
&|u(\ty,\overline{\yb}) -u(\tx,\overline{\xb})|^p-|u(0,\overline{\yb}) -u(0,\overline{\xb})|^p \\
\leq &Cp (|u(\ty,\overline{\yb}) -u(\tx,\overline{\xb}) |^{p-1}+|u(0,\overline{\yb}) -u(0,\overline{\xb}) |^{p-1})   \left| u(\ty,\overline{\yb}) -u(0,\overline{\yb})  -( u(\tx,\overline{\xb}) -  u(0,\overline{\xb}))\right|\\
\lesssim &\del \min ( |\yb -\xb|^{p-1} , \| u\|^{p-1}_{\infty}) \lesssim \del \min ( |\yb -\xb|^{p-1} , 1). 
\end{split}
\]
 On the other hand if $0\leq u(0,\overline{\yb}) -u(0,\overline{\xb})\leq 4 \del \tilde C$, then we have $|u(\ty,\overline{\yb}) -u(\tx,\overline{\xb})|\leq 6\del \tilde C$. Therefore, \eqref{eq:prop2.1A2} is true and this leads to
\begin{align*}
A_2\lesssim &\left|\del^{\,-2}\int_{(-\del,0)\times B^{d-1}(\bm 0,2M)} \int_{(-\del,0)\times B^{d-1}(\bm 0,2M)}  \frac{ \max(\del\verti{\yb-\xb}^{p-1}, \del^p)}{(|{\bm y}-{\bm x}|\vee\del)^{d+p-2}} \,d\yb d\xb\right|\\
&\hspace{3cm} +\left|\del^{\,-2}\int_{(-\del,0)\times B^{d-1}(\bm 0,M)} \int_{\mcR_{-\del}\backslash (-\del,0)\times B^{d-1}(\bm 0,2M)} \frac{\del }{(|{\bm y}-{\bm x}|\vee\del)^{d+p-2}} \,d\yb d\xb\right|\\
\lesssim &\,\del^{\,-1}\left(\int_{(-\del,0)\times B^{d-1}(\bm 0,2M)}   \int_{B^{d-1}(\bm 0,2M)}\int_{-\del}^0 \max\left(\frac{1}{(|\overline{\yb} - \overline{\xb}|\vee\del)^{d-1}},\frac{\del^{p-1}}{(|\overline{\yb} - \overline{\xb}|\vee\del)^{d+p-2}} \right) \,d\ty d\overline{\yb} d\xb\right.\\
&\hspace{3cm}
+\left.\int_{(-\del,0)\times B^{d-1}(\bm 0,M)} \int_{\R^{d-1}\backslash B^{d-1}(\overline{\bm x},M)}\int_{-\del}^0 \frac{1}{|\overline{\bm y}-\overline{\bm x}|^{d+p-2}} \,d\ty d\overline{\yb} d\xb\right)\\
\lesssim & \, \del(1-\log(\del))\xrightarrow{\del\to0} 0.
\end{align*}

{Lastly, for $A_3$ we note that
$$\verti{\overline{\yb}-\overline{\xb}}\leq\verti{\overline{\yb}-\overline{\xb}}\vee\del\leq \verti{\yb-\xb}\vee\del\leq \verti{\overline{\yb}-\overline{\xb}}+\del,$$
and therefore
$$\frac{1}{|\overline{\yb}-\overline{\xb}|^{d+p-2}}-\frac{1}{(\verti{\yb-\xb}\vee\del)^{d+p-2}}\leq \frac{1}{|\overline{\yb}-\overline{\xb}|^{d+p-2}}-\frac{1}{(|\overline{\yb}-\overline{\xb}|+\del)^{d+p-2}}.$$
With the fact that 
\begin{align*}
&\lim_{\del\to0}\left(\iint_{{\real^{2d-2}}}\dfrac{|u(0,\overline{\yb})-u(0,\overline{\xb})|^p}{|\overline{\yb}-\overline{\xb}|^{d+p-2}}d\overline{\yb} d\overline{\xb}-\iint_{{\real^{2d-2}}} \frac{|u(0,\overline{\yb})-u(0,\overline{\xb})|^p}{(|\overline{\bm y}-\overline{\bm x}|+\del)^{d+p-2}} \,d\overline{\yb} d\overline{\xb}\right)=0\,, \end{align*}
where the limits are achieved by the dominated convergence theorem, we then obtain
\begin{align*}
A_3=&\left|\int_{\real^{d-1}} \int_{\real^{d-1}}\dfrac{|u(0,\overline{\yb})-u(0,\overline{\xb})|^p}{|\overline{\yb}-\overline{\xb}|^{d+p-2}}d\overline{\yb} d\overline{\xb}-\del^{-2}\int_{\real^{d-1}}\int_{-\del}^0 \int_{\real^{d-1}}\int_{-\del}^0 \frac{|u(0,\overline{\yb})-u(0,\overline{\xb})|^p}{(|{\bm y}-{\bm x}|\vee\del)^{d+p-2}} \,d\tilde{y}d\overline{\yb}d\tilde{x} d\overline{\xb}\right|\\
=&\left|\del^{-2}\int_{\real^{d-1}}\int_{-\del}^0 \int_{\real^{d-1}}\int_{-\del}^0\dfrac{|u(0,\overline{\yb})-u(0,\overline{\xb})|^p}{|\overline{\yb}-\overline{\xb}|^{d+p-2}}-\frac{|u(0,\overline{\yb})-u(0,\overline{\xb})|^p}{(|{\bm y}-{\bm x}|\vee\del)^{d+p-2}} \,d\tilde{y}d\overline{\yb}d\tilde{x} d\overline{\xb}\right|\\
\leq&\left|\int_{\real^{d-1}} \int_{\real^{d-1}}\dfrac{|u(0,\overline{\yb})-u(0,\overline{\xb})|^p}{|\overline{\yb}-\overline{\xb}|^{d+p-2}}d\overline{\yb} d\overline{\xb}-\int_{\real^{d-1}} \int_{\real^{d-1}} \frac{|u(0,\overline{\yb})-u(0,\overline{\xb})|^p}{(|\overline{\bm y}-\overline{\bm x}|+\del)^{d+p-2}} \,d\overline{\yb} d\overline{\xb}\right|\overset{\del\rightarrow 0}{\longrightarrow} 0.   
\end{align*}}
\end{proof}

\subsection{Change of Variables and Scaling Identities}\label{sec:scale}

To further understand the trace theorem and nonlocal spaces of a given $\delta$ {and the connections with existing studies in the literature, we consider some scaling identities here. We recall the notation introduced in \eqref{stripenotation} so} that for any $\xb\in\mcR_{-1}^L$, we have $\del\xb\in \mcR_{-\del}^{L\del}$, which leads to the following scaling argument:
\begin{lem}\label{lem:scale}
Given $L>0$, for any {$u$ in ${\cS^{\,\beta}_\del(\mcR_{-\del}^L)}$ or
$\cT^{\,\beta}_\del(\mcR_{-\del})$}, let $v(\xb):=u(\del\xb)$, then {$v$ belongs to ${\cS^{\,\beta}_1(\mcR_{-1}^{L/\del})}$ or $\cT^{\,\beta}_1(\mcR_{-1})$ and we have respectively} 
\[\del^{d-p}\verti{v}^p_{\cS^{\,\beta}_1(\mcR_{-1}^{L/\del})}=\verti{u}^p_{\cS^{\,\beta}_\del(\mcR_{-\del}^{L})},\quad \del^{d-p}\verti{v}^p_{\cT^{\,\beta}_1(\mcR_{-1})}=\verti{u}^p_{\cT^{\,\beta}_\del(\mcR_{-\del})},\]
\[\del^{d}\vertii{v}^p_{L^p(\mcR_{-1}^{L/\del})}=\vertii{u}^p_{L^p(\mcR_{-\del}^{L})},\quad\del^{d}\vertii{v}^p_{L^p(\mcR_{-1})}=\vertii{u}^p_{L^p(\mcR_{-\del})}.\]
Moreover, the above results also hold for $L=\infty$:
\[\del^{d-p}\verti{v}^p_{\cS^{\,\beta}_1(\mcR_{-1}^{\infty})}=\verti{u}^p_{\cS^{\,\beta}_\del(\mcR_{-\del}^{\infty})},\quad \del^{d}\vertii{v}^p_{L^p(\mcR_{-1}^{\infty})}=\vertii{u}^p_{L^p(\mcR_{-\del}^{\infty})}.\]
\end{lem}
\begin{proof}
The proof is obtained by a change of variables. In particular, denoting $\wb=\yb/\del$ and $\zb=\xb/\del$, we have
{\begin{align}
\verti{u}^p_{\cS^{\,\beta}_\del(\mcR_{-\del}^{L})}=&\int_{\mcR_{-\del}^L}\int_{\mcR_{-\del}^L} \gamma^{\,\beta}_{\del}(|\yb-\xb|)|u(\yb)-u(\xb)|^pd\yb d\xb\nonumber\\
=&\del^{-d-p}\int_{\mcR_{-\del}^L}\int_{\mcR_{-\del}^L} \gamma^{\,\beta}_{1}\left(\frac{|\yb-\xb|}{\delta}\right)\left|v\left(\frac{\yb}{\del}\right)-v\left(\frac{\xb}{\del}\right)\right|^pd\yb d\xb \nonumber \\
=&\del^{d-p}\int_{\mcR_{-1}^{L/\del}}\int_{\mcR_{-1}^{L/\del}} \gamma^{\,\beta}_{1}\left(|\wb-\zb|\right)|v(\wb)-v(\zb)|^p d\wb d\zb=\del^{d-p}\verti{v}^p_{\cS^{\,\beta}_1(\mcR_{-1}^{L/\del})}.
\label{eq:scalenorm}
\end{align}
Similarly, for the trace norm we have
\begin{align*}
\verti{u}^p_{\cT^{\,\beta}_\del(\mcR_{-\del})}=&\del^{\,\beta-2}\int_{\mcR_{-\del}}\int_{\mcR_{-\del}} \frac{|u(\bm y) -u(\bm x)|^p}{(|\bm y-\bm x|{\vee} \del)^{d+p-2} (|\bm y -\bm x|\wedge\del)^{\,\beta}} \,d\yb d\xb\\
=&\del^{-d-p}\int_{\mcR_{-\del}}\int_{\mcR_{-\del}} \frac{\left| v\left(\frac{\bm y}{\del}\right) -v\left(\frac{\bm x}{\del}\right)\right|^p}{\left(\frac{|\bm y-\bm x|}{\del}{\vee}1 \right)^{d+p-2} \left(\frac{|\bm y -\bm x|}{\del}\wedge 1\right)^{\,\beta}} \,d\yb d\xb\\
=&\del^{d-p}\int_{\mcR_{-1}}\int_{\mcR_{-1}} \frac{(v(\bm w) -v(\bm z))^p}{(|\bm w-\bm z|{\vee} 1)^{d+p-2} (|\bm w -\bm z|\wedge 1)^{\,\beta}} \,d\wb d\zb=\del^{d-p}\verti{v}^p_{\cT^{\,\beta}_1(\mcR_{-1})}.
\end{align*}}
All other identities can be proved similarly. 
\end{proof}

\subsection{Equivalent Semi-norms}

We now introduce a lemma that allows us to compare the nonlocal spaces  $\cS_{\delta}^{\, \beta}$ with different sizes of $\delta$.   

\begin{lem}\label{lem:EnergyKernelEst}
Let $\alpha \in (0,1]$ and $U \subseteq \mathbb{R}^d$ a convex domain. There exists $C_1 = C_1(d,p)>0$ and $C_2= C_2(d, p, \alpha)>0$ such that for any $u \in \cS_\delta^\beta(U)$, 
\[
C_1|u|_{\cS_{\delta}^{\, \beta}(U)} \leq  |u|_{\cS_{\alpha\delta}^{\, \beta}(U)} \leq C_2 |u|_{\cS_{\delta}^{\, \beta}(U)}. 
\]
\end{lem}

\begin{proof}
First, note that 
\[
 |u|_{\cS_\delta^{\, \beta}(U)}^p = \int_{U}\int_{U} \gamma^{\,\beta}_{\delta}(|\by-\bx|)|u(\by)-u(\bx)|^p d \by d \bx .
\]
Choose an $m \in \mathbb{N}$ with $\dfrac{1}{m}<\alpha\leq \dfrac{1}{m-1}$. Notice that $m\geq 2$ for $\alpha\in (0,1]$. We split $u(\by)-u(\bx)$ into the following
\[
u(\by)-u(\bx)=\sum_{i=1}^{m} \left( u\left(\bx+\frac{i}{m}(\by-\bx)\right)-u\left(\bx+\frac{i-1}{m}(\by-\bx)\right)\right).
\]
Since $U$ is convex by assumption, for any $\bx\in U$ and $\by\in U$, we have $\bx+\frac{i}{m}(\by-\bx) \in U$ for each $0\le i \le m$ and so $u$ is well defined at these points.  Now, applying the inequality
$
\left| \sum_{i=1}^m a_i\right|^p \leq m^{p-1}\sum_{i=1}^m |a_i|^p
$,
 we have
\begin{align*}
    \int_{U}\int_{U}\gamma^{\,\beta}_{\delta}(|\by-\bx|)|u(\by)-u(\bx)|^p d \by d \bx& \le  m^{p-1}\sum_{i=1}^m\int_{U}\int_{U} \gamma_\delta^{\, \beta}(|\by-\bx|)\left|u\left(\bx+\frac{i}{m}(\by-\bx)\right)-u\left(\bx+\frac{i-1}{m}(\by-\bx)\right)\right|^p d \by d \bx.
    \end{align*}
 Notice that by the change of variables   $\bm w= \bx+\frac{i}{m}(\by-\bx)\in U$ and $\bz=\bx+\frac{i-1}{m}(\by-\bx)\in U$ we have $|\bm w-\bz|=  |\by-\bx|/m$ and the Jacobian matrix 
 \begin{equation}
 \label{eq:jabobina}
 \frac{\partial (\bm w, \bz)}{\partial (\by, \bx)} = 
 \begin{pmatrix}
\frac{i}{m} I_d & (1-\frac{i}{m} ) I_d \\
\frac{i-1}{m}  I_d & (1-\frac{i-1}{m}) I_d 
\end{pmatrix},
\end{equation}
 where $I_d\in \R^{d\times d}$ is the identity matrix.
 Thus $|\text{det}(\partial (\bm w, \bz)/\partial (\by, \bx))|= |\text{det}((\frac{i}{m}(1-\frac{i-1}{m}) -\frac{i-1}{m}(1-\frac{i}{m}) ) I_d)| =m^{-d}$ and then 
    \begin{align*}
  \int_{U}\int_{U}\gamma^{\,\beta}_{\delta}(|\by-\bx|)|u(\by)-u(\bx)|^p d \by d \bx  &\leq m^p \int_{U}\int_{U} \gamma_\delta^{\, \beta}(m|\bm w-\bz|)\left|u(\bm w)-u\left(\bz\right)\right|^p m^d d \bm w d \bz\\
 &= \int_{U}\int_{U} \frac{C_{d,p,\,\beta}}{(\del/m)^{d+p-\beta}}\frac{1}{|\bm w-\bz|^\beta} 1_{|\bm w-\bz|<\del/m}\left|u(\bm w)-u\left(\bz\right)\right|^p d \bm w d \bz
 \\
    &\leq \left(\frac{m}{m-1}\right)^{d+p-\beta} |u|^p_{\cS_{\alpha\delta}^{\, \beta}(U)} \leq 2^{d+p-\beta} |u|^p_{\cS_{\alpha\delta}^{\, \beta}(U)} \le {2^{d+p} |u|^p_{\cS_{\alpha\delta}^{\, \beta}(U)}},
\end{align*}
where we have used $1/m<\alpha\leq 1/(m-1)$ and $m\geq2$. So the left half of the inequality is true with $C_1= 2^{-d/p-1}$.

Lastly, the right half of the inequality is true with $C_2 = \alpha^{-d/p-1}$, since  
\[
\begin{split}
&|u|^p_{\cS_{\alpha\delta}^{\, \beta}(U)}  = \int_{U}\int_{U} \frac{C_{d,p,\,\beta}}{(\alpha\del)^{d+p-\beta}}\frac{1}{|\yb -\xb|^\beta} 1_{|\yb-\xb|<\alpha\del}\left|u(\yb)-u\left(\xb\right)\right|^p d\yb d\xb  \\
\leq &\alpha^{-d-p+\beta}\int_{U}\int_{U} \frac{C_{d,p,\,\beta}}{\del^{d+p-\beta}}\frac{1}{|\yb -\xb|^\beta} 1_{|\yb-\xb|<\del}\left|u(\yb)-u\left(\xb\right)\right|^p d\yb d\xb  \leq \alpha^{-d-p}  |u|^p_{\cS_{\delta}^{\, \beta}(U)}.
\end{split}
\]
\end{proof}

\subsection{Dyadic Cubes and Whitney Type Decomposition}\label{sec:whitney}

{The proof of Theorem \ref{mainthm_2} relies on extension results of Whitney type, the subject of which can be found in \cite{Stein1970}. 
Here we focus on defining Whitney type decompositions for the half space $\mcR^\infty:=\R^d_{+}$ and its subdomain $\mcR^L=(0,L)\times \R^{d-1}$.} 
For any $d\in\Z_+$, we define $\sQ_d$ the collection of dyadic cubes in $\R^d$, i.e., the cubes of the form {$Q=2^{-k}I(\mathbf{m})$
where $k\in\Z$ and
$I(\mathbf{m})=((0,1]^d+\mathbf{m})$ is the shifted unit cube for 
$\mathbf{m}\in\Z^d$}. Let $l(Q)$ denote the side length of the cube $Q\in \sQ_d$, and $\sQ_{d,k}$ the collection of cubes $Q\in \sQ_d$ with $l(Q)=2^{-k}$. 
For $\mcR^L= (0,L)\times \R^{d-1}$, we now define two types of decomposition of the domain $\mcR^L$ using the dyadic cubes for $L=2^m$ ($m\in \Z_+$), which will be useful in Section \ref{sec:inversetrace} to prove the inverse trace result. 
For any $k\in\Z$, we define
$\sW_k= \bigcup_{Q\in \sQ_{d-1,k}}(2^{-k}, 2^{-k+1}]\times Q $. 

\begin{itemize}
    \item Type I decomposition. Let $L= 2^m$ for some $m\in \Z_+\cup \{ 0\}$ and  $\overline{\sW_0}= \bigcup_{Q\in \sQ_{d-1,0}} (0,1]\times Q$, then 
    \begin{equation}\label{eq:decompositionI}
   \sW^I(\mcR^L):= \bigcup_{-m+1\leq k\leq 0,\, k\in \Z} \sW_k \cup  \overline{\sW_0}    
   \end{equation}
    \item Type II decomposition. Let $L= 2^m$ for some $m\in \Z_+\cup \{ 0\}$, then we define 
   \begin{equation} \label{eq:decompositionII}
    \sW^{II}(\mcR^L) := \bigcup_{-m+1\leq k,\, k\in\Z} \sW_k
   \end{equation}
\end{itemize}
Naturally, we will write $\sW^I(\mcR^\infty)=\bigcup_{k\in \Z_-\cup\{0\}} \sW_k \cup  \overline{\sW_0}$ and $\sW^{II}(\mcR^\infty)=\bigcup_{k\in \Z} \sW_k$. {Notice that $\sW^{II}(\mcR^\infty)$ coincides with the classical Whitney decomposition of the half space, where the length of each cube is proportional to the distance between the cube and the boundary of the domain. This type of decomposition is also used to prove the classical and fractional extension results \cite{dyda2019function,koskela2017traces}. The Type I decomposition, however, has a special set $\overline{\sW_0}$ which touches the boundary $\{0\}\times \R^{d-1}$ and it is used later to construct extension operator for the case $\beta<d$.}

\section{Nonlocal Trace Theorem}\label{sec:traceRd}

In this section we consider the trace theorem on half spaces and provide the proof for Theorem \ref{mainthm_1}. {We recall that the result stated corresponds to}
$\Om=\mcR^\infty=(0,\infty)\times\R^{d-1}$, $\Om_{\del}=\mcR_{-\del}=(-\del,0)\times\R^{d-1}$ and $\hat{\omg}=\mcR_{-\del}^{\infty}$. In particular, {with the help of the scaling arguments in Lemma \ref{lem:scale},} we first prove the results for the special case $\del=1$, then extend the results to general $\del$.
Since $C^1(\overline{\mcR_{-\del}^{\infty}})\cap \cS^{\,\beta}_\del(\mcR_{-\del}^{\infty})$ forms a dense set in $\cS^{\,\beta}_\del(\mcR_{-\del}^{\infty})$, it suffices to prove the conclusion for $u\in C^1(\overline{\mcR_{-\del}^{\infty}})\cap \cS^{\,\beta}_\del(\mcR_{-\del}^{\infty})$, {which is the case presented in the proofs.} 

When $\del=1$, we will prove the following theorem:
\begin{thm}[Trace theorem on half spaces when $\del=1$]\label{mainthm_1_d1}
Let $\beta\in[0,d+p)$,
then there exist a generic constants $C$ depending only on $d$ and $p$, such that for any $u\in \cS^{\,\beta}_1(\mcR_{-1}^{\infty})$ and any $\tau>0$,  
\begin{align}
 \label{eqn:L2result_flat_d1}\| u\|^p_{L^p(\mcR_{-1})}&\leq C |d+p-\beta|^{-1}\| u \|^{p-1}_{L^p(\mcR_{-1}^{\infty})}\verti{u}_{\cS^{\,\beta}_1(\mcR_{-1}^{\infty})}\leq C|d+p-\beta|^{-1}{\left(\tau^{-1} \| u \|^p_{L^p(\mcR_{-1}^{\infty})}+ \tau^{p-1} \verti{u}^p_{\cS^{\,\beta}_1(\mcR_{-1}^{\infty})}\right)},\\   
  \label{eqn:energynormresult_flat_d1}
  \verti{u}^p_{\cT^{\,\beta}_1(\mcR_{-1})}&\leq C|d+p-\beta|^{-1} \verti{u}^p_{\cS^{\,\beta}_1(\mcR_{-1}^{\infty})}.
\end{align}
\end{thm}

{To prove Theorem \ref{mainthm_1_d1}}, for any $\bm k=(k_1,k_2,\cdots,k_d)\in \mathbb{Z}^d$, we define the (hyper)cube associated with $\bm k$ by $I(\bm k)=(0,1)^{d+p-2}+{\bm k}= \prod_{i=1}^d (k_i, k_i+1)$. Now for any $\hb = (h_1, h_2, \cdots, h_d)\in \R^d$, we write $ [\hb] := ([h_1], [h_2], \cdots, [h_d])$, where 
\[
[h_i] = 
\left\{ 
\begin{aligned}
\lfloor h_i \rfloor \quad &  \text{if } h_i\geq 0  \\
 \lceil h_i \rceil \quad  &  \text{if } h_i< 0
\end{aligned}
\right. 
\]
for $i\in \{1,2, \cdots, d \}$. 
Then for all $\xb\in I(\bm k)$, we have $\xb +[\hb]\in I(\bm k +[\hb] )$. Notice that $I(\bm k +[\hb])$ has non-trivial intersections with the set $\{ \xb+\hb: \xb\in I(\bm k)\}$. Now we use $D(\bm k, [\hb])$ 
to denote the union of all (hyper)cubes in $\R^d$ that have non-trivial intersections with the diagonal line from the center of $I(\bm k)$ to the center of  $I(\bm k + [\hb])$. Then we have the following lemma. 
\begin{lem}
\label{lem:translation}
Let $u \in \cS_1^0(\R^{d})$, then
for any $\hb\in \R^d$, we have
\[
\frac{1}{\left( |\hb|+1 \right)^{p-1}}\int_{I(\bm k)} |u(\xb) -u(\xb+\hb)|^p d\xb \leq C \int_{D(\bm k,[\hb])} \int_{|\yb - \xb|<d+1} |u(\xb) - u(\yb)|^p d\yb d\xb , 
\]
where $C$ is a constant only dependent on $d$ and $p$. 
\end{lem}
\begin{proof}
Let $m$ be the number of (hyper)cubes in the set $D(\bm k, [\hb])$. Then we know that $m\leq |[\hb]|_{l_1} + 1\leq |\hb|_{l_1} + 1\leq \sqrt{d} |\hb| + 1$. We denote these (hyper)cubes by $I(\bm k+ \bm n^{(i)})$, where $\bm n^{(i)}\in \mathbb{Z}^d$ for each $i\in \{0,1,\cdots, m-1\}$ with $\bm n^{(0)}= \bf{0}$ and $\bm n^{(m-1)}= [\hb]$. Moreover, $|\bm n^{(i)} - \bm n^{(i+1)}|_{l_1} =1 $ for $i\in \{0,1,\cdots, m-2 \}$. Therefore, we can connect $\xb^{(0)}:= \xb$ and $\xb^{(m)}:= \xb + \hb$ by the points $\xb^{(i)}\in I(\bm k + \bm n^{(i)})$ for $i=1,2,\cdots, m-1$. Then
\[
|u(\xb) -u(\xb+\hb)|^p \leq 
\left| \sum_{i=0}^{m-1} \left(u(\xb^{(i)}) - u(\xb^{(i+1)}) \right) \right|^p \leq m^{p-1}  \sum_{i=0}^{m-1} |u(\xb^{(i)}) - u(\xb^{(i+1)})|^p. 
\]
Now integrate the above equation with respect to $\xb^{(i)}$ over $I(\bm k+ \bm n^{(i)})$ for each $i=0,1, \cdots, m-1$, we get 
\[
\begin{split}
&\int_{I(\bm k)} |u(\xb) - u(\xb+\hb)|^p dx  \\
\leq& m^{p-1} \left[ \sum_{i=0}^{m-2}  \int_{I(\bm k+\bm n^{(i)})} \int_{I(\bm k+\bm n^{(i+1)})} \left|u(\xb^{(i)}) - u(\xb^{(i+1)})\right|^p  d\xb^{(i+1)} d\xb^{(i)} + \int_{I(\bm k+\bm n^{(m-1)})}\int_{I(\bm k)}  \left| u(\xb^{(m-1)}) -u(\bx+\hb)\right|^p d\xb d\xb^{(m-1)} \right]  \\
 \leq& m^{p-1} \sum_{i=0}^{m-1} \int_{I(\bm k + \bm n^{(i)})} \int_{B(\xb,d+1)}  \left|u(\xb) - u(\yb)\right|^p d\yb d\xb  
 \leq  m^{p-1} \int_{D(\bm k,[\hb])} \int_{B(\xb,d+1)}  |u(\xb) - u(\yb)|^p d\yb d\xb .
\end{split}
\]
where $B(\xb,r)$ denotes the ball of radius $r$ in $\R^d$. The lemma is then a result of the above estimate and the fact that $m\leq \sqrt{d} |\hb| +1$.
\end{proof}

Lemma \ref{lem:translation}
has the following implication.
 \begin{cor}
 \label{cor:estimatebynewnorm}
 For any $m\in \Z_+$ and any  $u\in\cS_1^0(\mcR_{-1}^m)$, we have  
\begin{align}
\label{estimatebynewnorm-1}
  \int_{\R^{d-1}}  \int_{-1}^0 |u(\tilde{x}, \overline{\xb}) - u(\tilde{x}+m,  \overline{\xb}) |^p d\tilde{x} d\overline{\xb}  \leq C m^{p-1}  |u |^p_{\cS_1^0(\mcR_{-1}^m)}, 
\end{align}
where $C$ is a constant that depends only on $d$ and $p$.
\end{cor}
\begin{proof}
For $m\in \mathbb{Z}_+$, we 
use Lemma \ref{lem:translation}
with $\hb = (m, 0,\cdots, 0)$ to obtain 
\[
\begin{split}
\int_{\R^{d-1}}  \int_{-1}^0 |u(\tilde{x}, \overline{\xb}) - u(\tilde{x}+m,  \overline{\xb}) |^p d\tilde{x} d\overline{\xb}
&= \sum_{\overline{\bm k}\in \mathbb{Z}^{d-1}} \int_{I(\overline{\bm k})}  \int_{-1}^0 |u(\tilde{x}, \overline{\xb}) - u(\tilde{x}+m,  \overline{\xb}) |^p d\tilde{x} d\overline{\xb}  \\
&\leq C (m+1)^{p-1} \sum_{\overline{\bm k}\in \mathbb{Z}^{d-1}} \int_{(-1,m)\times I(\overline{\bm k})} \int_{|\yb -\xb|<d+1}|u(\xb) - u(\yb)|^p  d\yb d\xb  \\
& = C(m+1)^{p-1}\int_{\mcR_{-1}^m} \int_{|\yb -\xb|<d+1}|u(\xb) - u(\yb)|^p d\yb d\xb  \leq C m^{p-1} |u|^p_{\cS_1^0(\mcR_{-1}^m)}, 
\end{split}
\]
where we have use Lemma \ref{lem:EnergyKernelEst} in the last step of the above estimate.  \end{proof}

Lastly, we can show the $L^p$ estimate in Theorem \ref{mainthm_1_d1}. Note that 
\begin{equation}
\label{eq:generalb}
\begin{split}
| u |^p_{\cS^0_1(\mcR_{-1}^{\infty})}=C_{d,p,0}\int_{\mcR_{-1}^{\infty}}  \int_{\mcR_{-1}^{\infty}} 1_{\{|\bm y - \bm x|<1\}} \left| u(\bm y)-u(\bm x)\right|^p d\bm y d\bm x&\leq  \int_{\mcR_{-1}^{\infty}}  \int_{\mcR_{-1}^{\infty}} \dfrac{C_{d,p,0}1_{\{|\bm y - \bm x|<1\}}}{| \bm y-\bm x|^{\,\beta}}  \left| u(\bm y)-u(\bm x)\right|^p d\bm y d\bm x\\
&\leq C |d+p-\beta|^{-1}| u |^p_{\cS^{\,\beta}_1(\mcR_{-1}^{\infty})},
\end{split}
\end{equation}
where $C$ only depends on $d$ and $p$. 
Therefore it suffices to prove {the $L^p$ estimate}
with $\beta=0$ and then invoke the above inequality for general $\beta$. 
\begin{lem}[A nonlocal embedding lemma]\label{lem:L2d1}
For any $u\in \cS^0_1(\mcR_{-1}^{\infty})$ and $L>0$,
\begin{equation}
\label{eqn:L2}
\| u\|^p_{L^p(\mcR_{-1})}\leq C \left[ L^{-1}\| u\|^p_{L^p(\mcR_{-1}^L)} + L^{p-1} | u|^p_{\cS^0_1(\mcR_{-1}^L)} \right],
\end{equation}
where  $C$ is a constant independent of $u$ and $L$. 
\end{lem}

\begin{proof}
First, if $L<1$, then \eqref{eqn:L2} is trivially satisfied. Now we assume $L\geq 1$, then $1\leq \lfloor L \rfloor \leq L \leq 2\lfloor L \rfloor$. 
For any $m\in \Z_+$, we have 
\[
|u(\tilde{x}, \overline{\xb})|^p \leq 2^{p-1} \left[ | u(\tilde{x}, \overline{\xb}) - u(\tilde{x}+m, \overline{\xb}) |^p + |u(\tilde{x}+m, \overline{\xb})|^p \right]. 
\]
Therefore, combining the above inequality with Corollary \ref{cor:estimatebynewnorm}, we get
\[
\begin{split}
\| u\|^p_{L^p(\mcR_{-1})} =  \int_{\R^{d-1}}  \int_{-1}^0 |u(\tilde{x}, \overline{\xb})|^p   d\tilde{x} d\overline{\xb} & \leq 2^{p-1} \left[  \int_{\R^{d-1}}  \int_{-1}^0 |u(\tilde{x}, \overline{\xb}) - u(\tilde{x}+m,  \overline{\xb}) |^p d\tilde{x} d\overline{\xb} + \int_{\R^{d-1}}  \int_{-1}^0   |u(\tilde{x}+m, \overline{\xb})|^p  d\tilde{x} d\overline{\xb}  \right] \\ 
&\leq C \left[ m^{p-1} | u|^p_{\cS_1^0(\mcR_{-1}^m)}  +  \int_{\R^{d-1}}  \int_{m-1}^m   |u(\tilde{x}, \overline{\xb})|^p  d\tilde{x} d\overline{\xb} \right]. 
\end{split}
\]
Take a summation of the above inequality for $m$ from $1$ to $\lfloor L \rfloor$, we get 
\[
\lfloor L \rfloor \| u\|^p_{L^p(\mcR_{-1})}  \leq C \lfloor L \rfloor^p | u|^p_{\cS_1^0(\mcR_{-1}^{\lfloor L \rfloor})} + \|u \|^p_{L^p(\mcR_{-1}^{\lfloor L \rfloor})}, 
\]
which implies equation \eqref{eqn:L2}. 
\end{proof}

In the following we proceed to show the proof of \eqref{eqn:energynormresult_flat_d1}. {In a similar spirit to Corollary
\ref{cor:estimatebynewnorm} but with a more refined consideration than the application of a direct H\"{o}lder's inequality,  we first state an intermediate result.}

\begin{lem}\label{lem:botdiff_1}
{There exists a positive constant $C$ depending only on $d$  and $p$ such that for any positive integer $m$} and for any $
u\in\cS^{\,0}_1(\mcR_{-1}^{\infty})$,
the following estimate holds:
\begin{align*}
&\nonumber\int_{\mcR_{-1}}\left( \int_m^{m+1} |u(\tx,\overline{\xb})-u(\tz,\overline{\bm{x}})| d\tz  \right)^p d{\xb}\\
\leq &C |u |^p_{\cS_{1}^0(\mcR_{-1}^\infty)}
+C \left(\int_{-1}^{m}\left(\int_{(\lfloor\tilde{l}\rfloor,\lfloor\tilde{l}\rfloor+2)\times\real^{d-1}}\int_{B(\zb, d+1)\cap\mcR_{-1}^\infty}\verti{u(\wb)-u(\zb)}^p d\wb d\zb\right)^{1/p}d\tl\right)^p.
\end{align*}
\end{lem}

\begin{proof}
For $\tz\in {I}{(m)}:=(m,m+1)$, let us consider the covering of the path from $\tx$ to $\tz$ given by 
$I(i):=\left({i},{i+1}\right)$ for $i=0,\cdots,m$.
We take $\tz^{(i)}\in I(i)$ for $i\in \{ 0,1,\cdots, m-1\}$ and set $\tz^{(m)}:=\tz$. Then, since
\begin{align*}
|u(\tx,\overline{\xb})-u(\tz,\overline{\bm{x}})|
\le |u(\tx,\overline{\bm{x}})-u(\tz^{(0)},\overline{\bm{x}})|+
 \sum_{i=0}^{m-1}|u(\tz^{(i)},\overline{\bm{x}})-u(\tz^{(i+1)},\overline{\bm{x}})|,
\end{align*}
integrating the above inequality over $I{(i)}$ with respect to $\tz^{(i)}$ for each $i=0,\cdots,m$, and taking  both sides to the power of $p$ yields:
\begin{align*}
&\left(\int_{I{(m)}}
\nonumber|u(\tx,\overline{\xb})-u(\tz,\overline{\bm{x}})| d\tz \right)^p  \\
\le& 2^{p-1}\left(\int_{I{(0)}}|u(\tx,\overline{\xb})-u(\tz^{(0)},\overline{\bm{x}})| d\tz^{(0)} \right)^p +2^{p-1} \left(\sum_{i=0}^{m-1}\int_{I{(i)}}\int_{I{(i+1)}}|u(\tz^{(i)},\overline{\bm{x}})-u(\tz^{(i+1)},\overline{\bm{x}})|d\tz^{(i+1)}d\tz^{(i)}\right)^p.
\end{align*}
Now, we integrate the above inequality over $\mcR_{-1}$ with respect to ${\xb}$:
\begin{align}
\nonumber&\int_{\mcR_{-1}}\left( \int_{I{(m)}} |u(\tx,\overline{\xb})-u(\tz,\overline{\bm{x}})| d\tz \right)^p   d{\xb}
\leq2^{p-1}\int_{\mcR_{-1}}\left(\int_{I{(0)}}\verti{u(\tx,\overline{\bm{x}})-u(\tz^{(0)},\overline{\bm{x}})}d\tz^{(0)}\right)^p d{\xb}\\
&\qquad\qquad\qquad\quad +2^{p-1}\int_{\mcR_{-1}}\left(\sum_{i=0}^{m-1}\int_{I{(i)}}\int_{I{(i+1)}}\verti{u(\tz^{(i)},\overline{\bm{x}})-u(\tz^{(i+1)},\overline{\bm{x}})}d\tz^{(i+1)}d\tz^{(i)}\right)^p d{\xb}.\label{eqn:H13A1new}
\end{align}
For the first term above, 
{
\begin{align*}
\int_{\mcR_{-1}}\left(\int_{I{(0)}}\verti{u(\tx,\overline{\xb})-u(\tz^{(0)},\overline{\bm{x}})}d\tz^{(0)}\right)^p d{\xb}\leq
& \int_{\real^{d-1}}\int_{-1}^0\int_{0}^{1} \verti{u(\tx,\overline{\xb})-u(\tz^{(0)},\overline{\bm{x}})}^p d\tz^{(0)} d{\tx}d\overline{\xb}\\
=
& \int_{0}^{1} \int_{\real^{d-1}}\int_{-1}^0 \verti{u(\tx,\overline{\xb})-u(\tz^{(0)},\overline{\bm{x}})}^p 
d{\tx}d\overline{\xb}
d\tz^{(0)}\\
\leq & C |u |^p_{\cS_{1}^0(\mcR_{-1}^\infty)},
\end{align*}
where the last step is obtained by Corollary
\ref{cor:estimatebynewnorm}.
}
For the second term in \eqref{eqn:H13A1new} we use the Minkowski's integral inequality \cite{schep1995minkowski} and Lemma \ref{lem:translation}:
\begin{align*}
\nonumber&\int_{\mcR_{-1}}\left(\sum_{i=0}^{m-1}\int_{I{(i)}}\int_{I{(i+1)}}\verti{u(\tz^{(i)},\overline{\bm{x}})-u(\tz^{(i+1)},\overline{\bm{x}})}d\tz^{(i+1)}d\tz^{(i)}\right)^p d{\xb}\\
\leq&C\left(\sum_{i=0}^{m-1}\int_{I{(i)}}\left( \int_{\real^{d-1}}\left(\int_{I{(i+1)}}\verti{u(\tz^{(i)},\overline{\xb})-u(\tz^{(i+1)},\overline{\xb})}d\tz^{(i+1)}\right)^p d\overline{\xb} \right)^{1/p} d\tz^{(i)}\right)^p\\
\leq&C\left(\sum_{i=0}^{m-1}\int_{I{(i)}}\left( \sum_{\overline{\bm{k}}\in\mathbb{Z}^{d-1}}\int_{I(\overline{\bm{k}})}\int_{I{(i+1)}}\verti{u(\tz^{(i)},\overline{\xb})-u(\tz^{(i+1)},\overline{\xb})}^pd\tz^{(i+1)} d\overline{\xb} \right)^{1/p} d\tz^{(i)}\right)^p\\
\nonumber\leq&C\left(\int_{-1}^{m}\left( \sum_{\overline{\bm{k}}\in\mathbb{Z}^{d-1}}\int_{I(\overline{\bm{k}})}\int_{\lfloor\tl\rfloor+1}^{\lfloor\tl\rfloor+2}\verti{u(\tw+(\tl-\tw),\overline{\xb})-u(\tw,\overline{\xb})}^pd\tw d\overline{\xb}\right)^{1/p} d\tl\right)^p\\
\nonumber\leq&C\left(\int_{-1}^{m}\left(\sum_{\overline{\bm{k}}\in\mathbb{Z}^{d-1}}\int_{D((\lfloor\tilde{l}\rfloor+1,\overline{\bm{k}}),([\tl-\tw],\bm{0}))}\int_{B(\zb, d+1)\cap\mcR_{-1}^\infty}\verti{u(\wb)-u(\zb)}^pd\wb d\zb\right)^{1/p} d\tl\right)^p\\
\nonumber\leq&C\left(\int_{-1}^{m}\left(\sum_{\overline{\bm{k}}\in\mathbb{Z}^{d-1}}\int_{D((\lfloor\tilde{l}\rfloor+1,\overline{\bm{k}}),(-1,\bm{0}))}\int_{B(\zb, d+1)\cap\mcR_{-1}^\infty}\verti{u(\wb)-u(\zb)}^pd\wb d\zb \right)^{1/p} d\tl\right)^p\\
\nonumber\leq&C\left(\int_{-1}^{m}\left(\int_{(\lfloor\tilde{l}\rfloor,\lfloor\tilde{l}\rfloor+2)\times\real^{d-1}}\int_{B(\zb, d+1)\cap\mcR_{-1}^\infty}\verti{u(\wb)-u(\zb)}^pd\wb d\zb\right)^{1/p} d\tl\right)^p.
\end{align*}
In the above derivation, we have used the fact that $[\tl-\tw]=0$ or  $-1$ for $\tw\in({\lfloor\tl\rfloor+1},{\lfloor\tl\rfloor+2})$ and therefore, the {corresponding sets of (hyper)cube, as defined earlier, satisfy} $D((\lfloor\tilde{l}\rfloor+1,\overline{\bm{k}}),([\tl-\tw],\bm{0}))\subset D((\lfloor\tilde{l}\rfloor+1,\overline{\bm{k}}),(-1,\bm{0}))$.
\end{proof}

\begin{remark}
{We note that the terms in the inequality of the above lemma are well defined by H\"{o}lder's inequality. Indeed, it is easy to see that
$$
\left(\int_{-1}^{m}\left( \int_{(\lfloor\tilde{l}\rfloor,\lfloor\tilde{l}\rfloor+2)\times\real^{d-1}}\int_{B(\zb, d+1)\cap\mcR_{-1}^\infty}\verti{u(\wb)-u(\zb)}^p d\wb d\zb\right)^{1/p} d\tl\right)^p
\leq C m^{p-1}  \verti{u}^p_{\cS^{\,0}_1(\mcR_{-1}^{\infty})},
$$
for a positive constant $C$ depending only on $d$ and $p$. Our goal is, however, to show a much refined bound so that the dependence on $m$ in the above inequality can be dropped.}
\end{remark}

{Next, we derive a result that helps us to obtain an estimate related to the second term in the above lemma. For $\xb=(\tx,\overline{\xb})$, $\yb=(\ty,\overline{\yb})\in \mcR_{-1}$, we present some lemmas to bound the estimates on the nonlocal differences from $\xb$ to $(\tz,\overline{\xb})$ and the nonlocal difference from $(\tz,\overline{\xb})$ to $(\tz,\overline{\yb})$, respectively.
For any $\xb \in \mcR_{-1}$, we let   
$\bm k(\xb)$ be an integer lattice point associated with $\xb$ such that
$I(\bm k(\xb))$ be a (hyper)cube containing $\xb$. Note that the association may not be unique if
$\xb$ is on the boundary of some open (hyper)cube $I(\bm k)$, integer lattice point. In such a case, we may select any of the neighboring $I(\bm k)$ to be the associated  (hyper)cube. Naturally, if
$\xb$ is an integer lattice point itself, we can use the default choice $I(\xb)$. 
For $\overline{\hb}\in \real^{d-1}$, we let
 $D(\bm k(\xb), [(0, \overline{\hb} )])$ be the collection of (hyper)cubes associated with  $\bm k(\xb)$ and
 $\hb=(0, \overline{\hb} )$ defined previously.
 We denote $m=m(\xb,\overline{\hb})$ as the number of (hyper)cubes in  $D(\bm k(\xb), [(0, \overline{\hb} )])$. We again use ${I}{(m)}:=(m,m+1)$ for any integer $m$.}
 
 \begin{lem}\label{lem:botdiff_2_new}
{There exists a positive constant $C$  depending only on $d$ and $p$ such that for any $
u\in\cS^{\,0}_1(\mcR_{-1}^{\infty})$ {and any $\xb\in\mcR_{-1}$}, the following holds:}
\begin{equation}
\int_{\real^{d-1}}\left(\int_{-1}^{{m(\xb,\overline{\hb})}}\left( \int_{(\lfloor\tilde{l}\rfloor,\lfloor\tilde{l}\rfloor+2)\times\real^{d-1}}\int_{B(\zb, d+1)\cap\mcR_{-1}^\infty}\verti{u(\wb)-u(\zb)}^pd\wb d\zb\right)^{1/p} d\tl\right)^p \frac{d\overline{\hb}}{\left(\verti{\overline{\hb}}+1\right)^{d+p-2}}
\leq C \verti{u}^p_{\cS^{\,0}_1(\mcR_{-1}^{\infty})}.
\end{equation}
\end{lem}

\begin{proof}
{we
 write $\overline{\hb}$ in the spherical coordinate to get}
\begin{align*}
\nonumber&\int_{\real^{d-1}}\left(\int_{-1}^{
{m(\xb,\overline{\hb})}
}\left( \int_{(\lfloor\tilde{l}\rfloor,\lfloor\tilde{l}\rfloor+2)\times\real^{d-1}}\int_{B(\zb, d+1)\cap\mcR_{-1}^\infty}\verti{u(\wb)-u(\zb)}^pd\wb d\zb\right)^{1/p}d\tl\right)^p\frac{d\overline{\hb}}{\left(\verti{\overline{\hb}}+1\right)^{d+p-2}}\\
\nonumber\leq& C\int_{0}^{\infty}\left(\int_{-1}^{(r+1)\sqrt{d-1}-1}\left( \int_{(\lfloor\tilde{l}\rfloor,\lfloor\tilde{l}\rfloor+2)\times\real^{d-1}}\int_{B(\zb, d+1)\cap\mcR_{-1}^\infty}\verti{u(\wb)-u(\zb)}^pd\wb d\zb \right)^{1/p} d\tl\right)^pr^{d-2}\frac{dr}{\left(r+1\right)^{d+p-2}}\\
\nonumber\leq& C\int_{0}^{\infty}\left(\int_{0}^{\hat{r}}\left(\int_{(\lfloor\hat{l}\rfloor-1,\lfloor\hat{l}\rfloor+1)\times\real^{d-1}}\int_{B(\zb, d+1)\cap\mcR_{-1}^\infty}\verti{u(\wb)-u(\zb)}^pd\wb d\zb\right)^{1/p}d\hat{l}\right)^p\frac{d\hat{r}}{\hat{r}^p}\\
\nonumber\leq& C\int_{0}^{\infty} \hat{r}^{-p} \left(\int_{0}^{\hat{r}}
f(\hat{l})
d\hat{l} \right)^p
d\hat{r}
\end{align*}
where 
$$\hat{l}:=\tilde{l}+1,\; \hat{r}:=(r+1)\sqrt{d-1},\;
f(\hat{l}):=
\left( \int_{(\lfloor\hat{l}\rfloor-1,\lfloor\hat{l}\rfloor+1)\times\real^{d-1}}\int_{B(\zb, d+1)\cap\mcR_{-1}^\infty}\verti{u(\wb)-u(\zb)}^pd\wb d\zb \right)^{1/p},
$$
and we have used the fact that ${m(\xb,\overline{\hb})} \leq |\overline{\hb}|_{l_1}+1\leq \sqrt{d-1}|\overline{\hb}|+1\leq \sqrt{d-1}(|\overline{\hb}|+1)-1$.

Using the Hardy's inequality \cite{hardy1952inequalities,avkhadiev2006hardy,avkhadiev2014hardy}:
{$$\int_0^\infty x^{-p}\left(\int_{0}^x f(y)dy\right)^pdx\leq \left(\frac{p}{p-1}\right)^p\int_0^\infty (f(y))^pdy,$$}
we get
\begin{align}
\nonumber&\int_{\real^{d-1}}\left(\int_{-1}^{
{m(\xb,\overline{\hb})}
}\left(\int_{(\lfloor\tilde{l}\rfloor,\lfloor\tilde{l}\rfloor+2)\times\real^{d-1}}\int_{B(\zb, d+1)\cap\mcR_{-1}^\infty}\verti{u(\wb)-u(\zb)}^pd\wb d\zb\right)^{1/p}d\tl\right)^p\frac{d\overline{\hb}}{\left(\verti{\overline{\hb}}+1\right)^{d+p-2}}\\
\nonumber\leq& C\int_{0}^{\infty}\int_{(\lfloor\hat{l}\rfloor-1,\lfloor\hat{l}\rfloor+1)\times\real^{d-1}}\int_{B(\zb, d+1)\cap\mcR_{-1}^\infty}\verti{u(\wb)-u(\zb)}^pd\wb d\zb d\hat{l}\\
\leq&C\int_{\mcR_{-1}^\infty}\int_{B(\zb, d+1)\cap\mcR_{-1}^\infty}\verti{u(\wb)-u(\zb)}^pd\wb d\zb\leq C |{u}|^p_{\cS^0_{d+1}(\mcR_{-1}^\infty)}\leq C|{u}|^p_{\cS^0_{1}(\mcR_{-1}^\infty)}.\label{eqn:lemma39p2}
\end{align}
\end{proof}

With the above lemma, we can then bound the nonlocal difference from $\xb=(\tx,\overline{\xb})$ to $(\tz,\overline{\xb})$ with the trace semi-norm.
Such an estimate can be seen as the  norm of the nonlocal variations along the normal direction (with respect to the strip domain) being controlled by the nonlocal semin-norm for $\cS^{\,0}_1(\mcR_{-1}^{\infty})$.

\begin{lem}\label{lem:botdiff_2}
There exists a positive constant $C$ depending only on $d$ and $p$ such that {for any $
u\in\cS^{\,0}_1(\mcR_{-1}^{\infty})$, the following holds:}
\begin{equation}
\int_{\real^{d-1}}\int_{\mcR_{-1}} \left( \int_{{I(m(\xb,\overline{\hb}))}}|u(\tx,\overline{\xb})-u(\tz,\overline{\xb})| d\tz \right)^p d{\xb} \frac{d\overline{\hb}}{\left(\verti{\overline{\hb}}+1\right)^{d+p-2}} \leq C \verti{u}^p_{\cS^{\,0}_1(\mcR_{-1}^{\infty})}.
\end{equation}
\end{lem}

\begin{proof}
With Lemma \ref{lem:botdiff_1} we have
\begin{align*}
\nonumber&\int_{\real^{d-1}}\int_{\mcR_{-1}}\left( \int_{I{(m(\xb,\overline{\hb})
)}} |u(\tx,\overline{\xb})-u(\tz,\overline{\bm{x}})| d\tz  \right)^pd{\xb}\frac{d\overline{\hb}}{\left(\verti{\overline{\hb}}+1\right)^{d+p-2}}\\
\leq& C
|u |^p_{\cS_{1}^0(\mcR_{-1}^\infty)} 
+C \int_{\real^{d-1}}\left(\int_{-1}^{{m(\xb,\overline{\hb})}}\left(\int_{(\lfloor\tilde{l}\rfloor,\lfloor\tilde{l}\rfloor+2)\times\real^{d-1}}\int_{B(\zb, d+1)\cap\mcR_{-1}^\infty}\verti{u(\wb)-u(\zb)}^pd\wb d\zb\right)^{1/p} d\tl\right)^p\frac{d\overline{\hb}}{\left(\verti{\overline{\hb}}+1\right)^{d+p-2}}.
\end{align*}
{Now, combining with Lemma \ref{lem:botdiff_2_new}, the proof is complete.}
\end{proof}

\begin{remark}
Apply a similar argument to the nonlocal difference from $\yb=(\ty,\overline{\yb})$ to $(\tz,\overline{\yb})$, we can similarly obtain: 
\begin{equation}
\int_{\real^{d-1}}\int_{\mcR_{-1}}\left( \int_{{I{(m
(\yb,\overline{\hb})
)}}
}|u(\ty,\overline{\yb})-u(\tz,\overline{\yb})| d\tz \right)^p d {\yb} \frac{d\overline{\hb}}{\left(\verti{\overline{\hb}}+1\right)^{d+p-2}} \leq C \verti{u}^p_{\cS^{\,0}_1(\mcR_{-1}^{\infty})}
\end{equation}
where $C$ is a positive constant depending only on $d$ and $p$.
\end{remark}

We now proceed to investigate the nonlocal difference from $(\tz,\overline{\xb})$ to $(\tz,\overline{\yb})$ in the following lemma, 
which can be seen analogously as the  norm of the nonlocal tangential variations being controlled by the nonlocal semin-norm:
\begin{lem}\label{lem:pardiff}
There exists a positive constant $C$ depending on $d$ and $p$ such that for 
{any $u\in\cS^{\,0}_1(\mcR_{-1}^{\infty})$
}, the following estimate holds:
\begin{equation}
\int_{\real^{d-1}}\int_{\real^{d-1}} \int_{{I{(m(\xb,\overline{\hb}))}}}|u(\tz,\overline{\xb})-u(\tz,\overline{\xb}+\overline{\hb})|^p d\tz d\overline{\xb} \frac{d\overline{\hb}}{\left(\verti{\overline{\hb}}+1\right)^{d+p-2}} \leq C \verti{u}^p_{\cS^{\,0}_1(\mcR_{-1}^{\infty})}.
\end{equation}
\end{lem}
\begin{proof}
{For notation simplicity, we drop the dependence of $m$ on its argument in the derivation here.}
First notice that from Lemma \ref{lem:translation} we have 
\[
\begin{split}
 \int_{\R^{d-1}}
\int_{I{(m)}} |u(\tilde{z}, \overline{\xb}) - u(\tilde{z}, \overline{\xb}+\overline{\hb})|^p d\tilde{z} d\overline{\xb} & = \sum_{\overline{\bm k}\in \mathbb{Z}^{d-1}}\int_{I(\overline{\bm k})}
\int_{I{(m)}} |u(\tilde{z}, \overline{\xb}) - u(\tilde{z}, \overline{\xb}+\overline{\hb})|^p d\tilde{z} d\overline{\xb} \\
&\leq C (|\overline{\hb}|+1)^{p-1}\sum_{\overline{\bm k}\in \mathbb{Z}^{d-1}}\int_{D((m, \overline{\bm k}), (0,[\overline{\hb}]))}
\int_{B(\zb,d+1)\cap \mcR_{-1}^\infty} |u(\wb) - u(\zb)|^p d\wb d\zb \\
&\leq C (|\overline{\hb}|+1)^p \int_{I{(m)}\times \R^{d-1}}
\int_{B(\zb,d+1)\cap \mcR_{-1}^\infty} |u(\wb) - u(\zb)|^p d\wb d\zb.
\end{split}
\]
The last inequality in the above estimate is due to the fact that there are at most $|[\overline{\hb}]|_{l_1} + 1 \leq \sqrt{d-1}|\overline{\hb}| + 1 $ (hyper)cubes in the set $D((m, \overline{\bm k}), (0,[\overline{\hb}]))$. Therefore, 
\[
\begin{split}
&\int_{\R^{d-1}}\int_{\R^{d-1}}
\int_{I{(m)}} |u(\tilde{z}, \overline{\xb}) - u(\tilde{z}, \overline{\xb}+\overline{\hb})|^p d\tilde{z} d\overline{\xb} \frac{d\overline{\hb}}{\left(\verti{\overline{\hb}}+1\right)^{d+p-2}} \\
\leq& C \int_{\R^{d-1}} \frac{d\overline{\hb}}{\left(\verti{\overline{\hb}}+1\right)^{d-2}} \int_{I{(m)}\times \R^{d-1}}
\int_{B(\zb,d+1)\cap \mcR_{-1}^\infty} |u(\wb) - u(\zb)|^p d\wb d\zb \\
\leq& C \int_0^\infty \frac{r^{d-2}}{(r+1)^{d-2}}dr \int_{I{(m)}\times \R^{d-1}}
\int_{B(\zb,d+1)\cap \mcR_{-1}^\infty} |u(\wb) - u(\zb)|^p d\wb d\zb  \leq C |u|^p_{\cS_1^0(\mcR_{-1}^\infty)}, 
\end{split}
\]
where we have also used Lemma \ref{lem:EnergyKernelEst} in the last inequality.
\end{proof}

We now have the following lemma for the trace semi-norm:
\begin{lem}\label{lem:Td1}
 There exist a positive constant $C$ depending only $d$ and $p$ such that for any $u\in \cS^{\,\beta}_1(\mcR_{-1}^{\infty})$,
\begin{equation}
\label{eqn:T}
\verti{u}^p_{\cT^{\,\beta}_1(\mcR_{-1})}\leq C |d+p-\beta|^{-1}\verti{u}^p_{\cS^{\,\beta}_1(\mcR_{-1}^{\infty})}\,.   
\end{equation}
\end{lem}
\begin{proof}
We first note that
\begin{align*}
&\int_{\mcR_{-1}}\int_{\mcR_{-1}\cap B(\xb,1)} \frac{|u(\bm y) -u(\bm x)|^p}{(|\bm y-\bm x|\vee 1)^{d+p-2} (|\bm y -\bm x|\wedge1)^{\,\beta}} \,d\yb d\xb=\int_{\mcR_{-1}}\int_{\mcR_{-1}\cap B(\xb,1)} \frac{|u(\bm y) -u(\bm x)|^p}{|\bm y -\bm x|^{\,\beta}} \,d\yb d\xb\leq C |d+p-\beta|^{-1}|{u}|^p_{\cS^{\,\beta}_{1}(\mcR_{-1}^\infty)}.
\end{align*}
Moreover, we note that $\verti{\overline{\yb}-\overline{\xb}}+1\leq 2\verti{\yb-\xb} =2 (\verti{\yb-\xb}\vee 1)$ for $\yb\in\mcR_{-1}\cap (\mcR_{-1}\backslash B(\xb,1))$, it then suffices to show
\begin{align*}
\underset{\mcR_{-1}}{\iint}\underset{{\mcR_{-1}\backslash B(\xb,1)}}{\iint} \frac{|u(\ty,\overline{\bm y}) -u(\tx,\overline{\bm x})|^p}{(|\overline{\bm y}-\overline{\bm x}|+1)^{d+p-2} } \,d\ty d\overline{\yb}d\tx d\overline{\xb}&\leq \underset{\mcR_{-1}}{\iint} \underset{\mcR_{-1}}{\iint}\frac{|u(\ty,\overline{\bm y}) -u(\tx,\overline{\bm x})|^p}{(|\overline{\bm y}-\overline{\bm x}|+1)^{d+p-2} } \,d\ty d\overline{\yb}d\tx d\overline{\xb} \leq C |u|^p_{\cS_1^0(\mcR_{-1}^\infty)}\,,
\end{align*}
since $|u|^p_{\cS_1^0(\mcR_{-1}^\infty)}\leq C|d+p-\beta|^{-1} |{u}|^p_{\cS^{\,\beta}_{1}(\mcR_{-1}^\infty)}$ for any $\beta\in [0,d+p)$ where $C$ only depends on $d$ and $p$. Taking the path from $\xb$ to $(\tz,\overline{\xb})$, $(\tz,\overline{\yb})$ and then finally $\yb$, we have
\begin{align}
\nonumber|u(\tx,\overline{\xb})-u(\ty,\overline{\yb})| \le
|u(\tx,\overline{\xb})-u(\tz,\overline{\xb})|+
|u(\ty,\overline{\yb})-u(\tz,\overline{\yb})|
+|u(\tz,\overline{\xb})-u(\tz,\overline{\yb})|.
\end{align}
Denoting $\overline{\hb}:=\overline{\yb}-\overline{\xb}$ and
integrating the above over $I{(m(\xb,\overline{\hb}))}$ with respect to $\tz$ and taking both hand sides to the power of $p$ yields:
\begin{align}
\nonumber|u(\tx,\overline{\xb})-u(\ty,\overline{\yb})|^p \le & 3^{p-1} \left( \int_{I{(m(\xb,\overline{\hb}))}}
|u(\tx,\overline{\xb})-u(\tz,\overline{\xb})| d\tz \right)^p+3^{p-1}   \left(  \int_{I{(m(\xb,\overline{\hb}))}}|u(\ty,\overline{\yb})-u(\tz,\overline{\yb})| d\tz \right)^p\\
\nonumber&+3^{p-1}   \left(  \int_{I{(m(\xb,\overline{\hb}))}}|u(\tz,\overline{\xb})-u(\tz,\overline{\yb})|d\tz \right)^p.
\label{eqn:H12A0}
\end{align}
{Considering a fixed $\overline{\hb}=\overline{\yb}-\overline{\xb}$, } we integrate the above inequality over $(-1,0)$ with respect to $\tx$, $\ty$, respectively, then integrate over $\real^{d-1}$ with respect to $\overline{\xb}$:
\begin{align}
\nonumber
&\int_{\real^{d-1}} \int_{-1}^0 \int_{-1}^0
|u(\tx,\overline{\xb})-u(\ty,\overline{\yb})|^p  d\ty d \tx d\overline{\xb}\\
\nonumber\le&3^{p-1}  \int_{\real^{d-1}}   \int_{-1}^0\left( \int_{I{(m(\xb,\overline{\hb}))}}|u(\tx,\overline{\xb})-u(\tz,\overline{\xb})| d\tz \right)^p d\tx d\overline{\xb} + 3^{p-1}  \int_{\real^{d-1}}  \int_{-1}^0 \left(  \int_{I{(m(\xb,\overline{\hb}))}}|u(\ty,\overline{\yb})-u(\tz,\overline{\yb})| d\tz \right)^p d\ty d\overline{\xb}\\
\nonumber & + 3^{p-1}   \int_{\real^{d-1}} \left(  \int_{I{(m(\xb,\overline{\hb}))}}|u(\tz,\overline{\xb})-u(\tz,\overline{\yb})|d\tz \right)^p d\overline{\xb}.
\end{align}
Multiplying the above inequalities with $\left(\verti{\overline{\hb}}+1\right)^{d+p-2}$ and integrating with respect to $\overline{\hb}$ over $\real^{d-1}$ yield:
\begin{align*}
&\int_{\real^{d-1}}\int_{\real^{d-1}}\int_{-1}^{0}\int_{-1}^{0}\dfrac{\verti{u(\tx,\overline{\xb})-u(\ty,\overline{\yb})}^p}{\left(\verti{\overline{\yb}-\overline{\xb}}+1\right)^{d+p-2}} d\ty d\tx d\overline{\xb}d\overline{\yb}\\
\leq&C\int_{\real^{d-1}}\int_{\real^{d-1}}   \int_{-1}^0\left( \int_{I{(m(\xb,\overline{\hb}))}}|u(\tx,\overline{\xb})-u(\tz,\overline{\xb})| d\tz \right)^p d\tx d\overline{\xb}\frac{d\overline{\hb}}{\left(\verti{\overline{\hb}}+1\right)^{d+p-2}}\\
&+C\int_{\real^{d-1}}\int_{\real^{d-1}}  \int_{-1}^0 \left(  \int_{I{(m(\xb,\overline{\hb}))}}|u(\ty,\overline{\xb}+\overline{\hb})-u(\tz,\overline{\xb}+\overline{\hb})| d\tz \right)^p d\ty d\overline{\xb}\frac{d\overline{\hb}}{\left(\verti{\overline{\hb}}+1\right)^{d+p-2}}\\
&+C\int_{\real^{d-1}}\int_{\real^{d-1}} \left(  \int_{I{(m(\xb,\overline{\hb}))}}|u(\tz,\overline{\xb})-u(\tz,\overline{\yb})|d\tz \right)^p d\overline{\xb}\frac{d\overline{\hb}}{\left(\verti{\overline{\hb}}+1\right)^{d+p-2}}\leq C \verti{u}^p_{\cS^{\,0}_1(\mcR_{-1}^{\infty})},
\end{align*}
where the last inequality {follows} immediate from Lemmas \ref{lem:botdiff_2}-\ref{lem:pardiff}.
\end{proof}

We now show the proof of Theorem \ref{mainthm_1_d1}  and Theorem \ref{mainthm_1}. 

\noindent{\it Proof of Theorem \ref{mainthm_1_d1} and Theorem \ref{mainthm_1}.} 
From Lemma \ref{lem:L2d1}, we have
\[
\| u\|^p_{L^p(\mcR_{-1})}\leq C L^{-1}\| u \|^p_{L^p(\mcR_{-1}^\infty)}+CL^{p-1}\verti{u}^p_{\cS^0_1(\mcR_{-1}^\infty)}\leq C |d+p-\beta|^{-1}(L^{-1}\| u \|^p_{L^p(\mcR_{-1}^\infty)}+L^{p-1}\verti{u}^p_{\cS^{\,\beta}_1(\mcR_{-1}^\infty)})\,,
\]
for any $L>0$. 
By taking $L=\| u \|_{L^p(\mcR_{-1}^{\infty})}/\verti{u}_{\cS^{\,\beta}_1(\mcR_{-1}^{\infty})}$ for $u\neq 0$ in the above inequality, we obtain \eqref{eqn:L2result_flat_d1}. 
\eqref{eqn:energynormresult_flat_d1} is an immediate result in Lemma \ref{lem:Td1}. 

The proof of the general nonlocal trace Theorem \ref{mainthm_1} on half spaces then follows from Theorem \ref{mainthm_1_d1} and the scaling argument in Lemma \ref{lem:scale}: for any $u\in \cS^{\,\beta}_\del(\mcR_{-\del}^\infty)$ and let $v(\xb)=u(\del\xb)$, then $v\in {\cS^{\,\beta}_1(\mcR_{-1}^{\infty})}$ and
\begin{align*}
\frac{1}{\del}\|  u\|^p_{L^p(\mcR_{-\del})}&=\del^{d-1}\|  v\|^p_{L^p(\mcR_{-1})}\leq C |d+p-\beta|^{-1}\del^{d-1}\| v \|^{p-1}_{L^p(\mcR_{-1}^\infty)}\verti{v}_{\cS^{\,\beta}_1(\mcR_{-1}^\infty)}\\
&=C |d+p-\beta|^{-1}\| u \|^{p-1}_{L^p(\mcR_{-\del}^\infty)}\verti{u}_{\cS^{\,\beta}_\del(\mcR_{-\del}^\infty)}\leq C |d+p-\beta|^{-1}\left(\| u \|^p_{L^p(\mcR_{-\del}^\infty)}+\verti{u}^p_{\cS^{\,\beta}_\del(\mcR_{-\del}^\infty)}\right).\\
| u|^p_{\cT^{\,\beta}_\del(\mcR_{-\del})}&=\del^{d-p}\verti{v}^p_{\cT^{\,\beta}_1(\mcR_{-1})}\leq C|d+p-\beta|^{-1} \del^{d-p}\verti{v}^p_{\cS^{\,\beta}_1(\mcR_{-1}^\infty)}=C|d+p-\beta|^{-1}\verti{u}^p_{\cS^{\,\beta}_\del(\mcR_{-\del}^\infty)}.
\end{align*}

\qed

\begin{remark}
{We have proved in Lemma \ref{lem:L2d1} that $\vertii{u}_{L^p(\mcR_{-1})}$ is bounded by the $L^p$ norm and the $\cS_1^{\,\beta}$ semi-norm on a stripe domain ${\mcR_{-1}^L}$, which {might appear as a stronger statement that implies the result on a half space ${\mcR_{-1}^\infty}$ in Theorem \ref{mainthm_1_d1}. Though, we note that the result on a stripe domain ${\mcR_{-1}^L}$ can also be a consequence of Theorem \ref{mainthm_1_d1}. 
Likewise,} a bound on the ${\cT_1^{\,\beta}(\mcR_{-1})}$ semi-norm can also be obtained on a stripe domain ${\mcR_{-1}^L}$:
$$\verti{u}^p_{\cT^{\,\beta}_1(\mcR_{-1})}\leq C|d+p-\beta|^{-1}\left(|u|^p_{\cS_1^{\,\beta}(\mcR_{-1}^{L})} + L^{-p}\|u\|_{L^p(\mcR_{-1}^{L})}^p\right).$$
{To verify the above conclusions,
 let us derive the results on the stripe domain from  \eqref{eqn:L2result_flat_d1} and \eqref{eqn:energynormresult_flat_d1}.}
We take a smooth cutoff function $\phi(x)\in C_c^{\infty}(\real)$ such that supp$(\phi)\subset [-1,1/2]$ and $\phi(x)=1$ for $x\leq 0$. Denoting $\tilde{\phi}(\xb):\real^d\rightarrow \real$ such that $\tilde{\phi}(\tilde{x},\overline{\xb})=\phi(\tilde{x}/L)$ for $\tilde{x}>0$ and $\tilde{\phi}(\tilde{x},\overline{\xb})=\phi(\tilde{x})$ for $\tilde{x}<0$, then we note that there exists a generic constant $C$ independent of $L$ and $\verti{\tilde{\phi}}_{C^0}\leq C$, $\verti{\tilde{\phi}'}_{C^0}\leq CL^{-1}$. Substituting $\tilde{\phi}u$ into \eqref{eqn:energynormresult_flat_d1}, we have
\[
\begin{split}
 &|u|^p_{\cT_1^{\,\beta}(\mcR_{-1})}=|\tilde{\phi}u|^p_{\cT_1^{\,\beta}(\mcR_{-1})}\leq C |d+p-\beta|^{-1} |\tilde{\phi}u|^p_{\cS_1^{\,\beta}(\mcR_{-1}^\infty)} \\
 \leq & C|d+p-\beta|^{-1} \left(\int_{\mcR_{-1}^{L}}\int_{\mcR_{-1}^{L}}+ 2\int_{\mcR_{-1}^{L/2}}\int_{\mcR^\infty\backslash\mcR^{L}}\right) \gamma^{\,\beta}_1(|\bx-\by|)|\tilde{\phi}(\by)u(\by)-\tilde{\phi}(\bx)u(\bx)|^pd\by d\bx   \\
 \leq & C|d+p-\beta|^{-1} \int_{\mcR_{-1}^{L}}\int_{\mcR_{-1}^{L}}\gamma^{\,\beta}_1(|\bx-\by|)|\tilde{\phi}(\by)|^p|u(\by)-u(\bx)|^p+\gamma^{\,\beta}_1(|\bx-\by|)|u(\bx)|^p|\tilde{\phi}(\by)-\tilde{\phi}(\bx)|^p d\by d\bx  \\
 &+ C|d+p-\beta|^{-1} \int_{\mcR_{-1}^{L/2}} |\tilde{\phi}(\bx)u(\bx)|^p \int_{L/2<|\by-\bx|<1}  \gamma^{\,\beta}_1(|\bx-\by|) d\by d\bx \\
 \leq&C|d+p-\beta|^{-1}\left(\|\tilde{\phi}\|^{p}_{C^0}|u|^p_{\cS_1^{\,\beta}(\mcR_{-1}^{L})} +\|\tilde{\phi}'\|^{p}_{C^0}\int_{\mcR_{-1}^{L}}|u(\bx)|^p\int_{\mcR_{-1}^{L} \cap B(\bx,1)}\gamma^{\,\beta}_1(|\bx-\by|)|\bx-\mathbf{y}|^pd\by d\bx\right.\\
 &\left.+ \| \tilde{\phi}\|_{C^0}^p \| u\|_{L^p(\mcR_{-1}^{L/2})}^p  \int_{L/2<|\bz|<1} \frac{1}{|\bz|^\beta} d\bz\right)\\
  \leq & C|d+p-\beta|^{-1}\left( |u|^p_{\cS_1^{\,\beta}(\mcR_{-1}^{L})} +{L^{-p}}\|u\|_{L^p(\mcR_{-1}^{L})}^p  + \| \tilde{\phi}\|_{C^0}^p \| u\|_{L^p(\mcR_{-1}^{L/2})}^p  \int_{L/2<|\bz|<1} \frac{1}{|\bz|^\beta} d\bz\right).
\end{split}
\] 
Notice that 
\[
\int_{L/2<|\bz|<1} \frac{1}{|\bz|^\beta} d\bz=0 \quad\text{ if } L\geq2,\]
otherwise 
\[
\int_{L/2<|\bz|<1} \frac{1}{|\bz|^\beta} d\bz \leq C L^{d-\beta} \leq C L^{-p}
\]
since $d-\beta > -p$. 
So we have $|u|^p_{\cT_1^{\,\beta}(\mcR_{-1})}\leq C|d+p-\beta|^{-1}\left( |u|^p_{\cS_1^{\,\beta}(\mcR_{-1}^{L})} +L^{-p}\|u\|_{L^p(\mcR_{-1}^{L})}^p\right) $ for all $L>0$.
} Similarly, substituting  $\tilde{\phi}u$ into \eqref{eqn:L2result_flat_d1} for $\beta=0$ yields
\begin{align*}
    &\| u\|^p_{L^p(\mcR_{-1})}=\|\tilde{\phi} u\|^p_{L^p(\mcR_{-1})}\leq C \| \tilde{\phi}u \|^{p-1}_{L^p(\mcR_{-1}^{\infty})}\verti{\tilde{\phi}u}_{\cS^{0}_1(\mcR_{-1}^{\infty})}\leq C\|u\|^{p-1}_{L^p(\mcR_{-1}^{L})}\left(|u|_{\cS_1^{0}(\mcR_{-1}^{L})}+L^{-1}\|u\|_{L^p(\mcR_{-1}^{L})}\right)\\
    &\leq CL^{-1}\| u \|^p_{L^p(\mcR_{-1}^L)}+C \| u \|^{p-1}_{L^p(\mcR_{-1}^L)}\verti{u}_{\cS^0_1(\mcR_{-1}^L)}\leq CL^{-1}\| u \|^p_{L^p(\mcR_{-1}^L)}+C L^{p-1}\verti{u}^p_{\cS^0_1(\mcR_{-1}^L)}.
\end{align*}
\end{remark}

\section{Nonlocal Inverse Trace Theorem}
\label{sec:inversetrace}

For $u: \mcR_{-1}\to\R$ and $L=2^m$ ($m\in\Z_+\cup\{0\}$), we now define the extension operator $E^L u: \mcR_{-1}^\infty\to\R$. 
Notice that the kernel $\ga^{\,\beta}_1$ is defined in \eqref{def:frackernel} with the two cases $\beta\in [0,d)$ and $\beta\in (d,d+p)$. We have the following two cases for the definition of $E^L$.\\
Case 1: $\beta\in [0,d)$. We define a partition of unity for $\mcR^L$ according to the decomposition $\sW^I(\mcR^L)$ defined in \eqref{eq:decompositionI}. For any $W\in \sW^I(\mcR^L)$, let $\phi_W^I: \R^d_+ \to [0,1]$ be a smooth function associated with $W$ such that $\phi_W^I$ is bounded below uniformly on $W$, Lip$\phi_W^I \lesssim 1/l(W)$ and supp$(\phi_W^I)$ is contained in an $l(W)/4$-neighborhood of $W$. Moreover, $\sum_{W\in \cW^I(\mcR^L)} \phi_W^I \equiv 1$ on $\mcR^L$. Notice that $\{ \phi^I_W\}$ should also depend on $L$ and here we drop the $L$ dependence for simplicity of notations.  The extension operator is then defined as
\begin{equation}\label{def:extensionop_1}
 E^L u (\bm x) = 
\left\{
\begin{aligned}
\sum_{W\in \sW^I(\mcR^L)} a_W^I \phi_W^I(\bm x) \quad &\bm x\in \mcR^\infty \,, \\
u(\bm x) \quad  &\bm x\in \mcR_{-1}\,,
\end{aligned}
\right.   
\end{equation}
where $$a_W^I:= \left(\int_{\sM_1(W)} u \right) \big/ |\sM_1(W)|,$$ 
and the map $\sM_1$ is defined for any $W=(a, b]\times Q \in \sW^I(\mcR^L) $ as 
\begin{equation} \label{def:map_M1}
\sM_1(W)= (-1, 0)\times Q\,. \end{equation}
Case 2: $\beta\in (d,d+p)$. We similarly define $\{ \phi_W^{II} \}_{W\in \sW^{II}(\mcR^L)}$ as  a partition of unity for $\mcR^L$ according to  $\sW^{II}(\mcR^L)$ defined in \eqref{eq:decompositionII}. More specifically, for any $W\in \sW^{II}(\mcR^L)$, $\phi_W^{II}: \R^d_+ \to [0,1]$ is a smooth function bounded below uniformly on $W$, Lip$\phi_W^{II} \lesssim 1/l(W)$ and supp$(\phi_W^{II})$ is contained in an $l(W)/4$-neighborhood of $W$. Moreover, $\sum_{W\in \cW^{II}(\mcR^L)} \phi_W^{II} \equiv 1$ on $\mcR^L$.    Then the extension operator is given by 
\begin{equation} \label{def:extensionop_2}
 E^L u (\bm x) = 
\left\{
\begin{aligned}
\sum_{W\in \sW^{II}(\mcR^L)} a_W^{II} \phi_W^{II}(\bm x) \quad &\bm x\in\mcR^\infty \,, \\
u(\bm x) \quad  &\bm x\in \mcR_{-1} \,,
\end{aligned}
\right.  
\end{equation}
where 
$$a_W^{II}:= \left(\int_{\sM_2(W)} u \right) \big/ |\sM_2(W)|,
$$
and the map $\sM_2$ for any $W=(a, b]\times Q \in \sW^{II}(\mcR^L)$ as 
\begin{equation} \label{def:map_M2}
\sM_2(W)= 
\left\{ 
\begin{aligned}
& (-b, -a)\times Q  \quad &\text{if } b\leq 1 \,,\\
& (-1, 0) \times Q \quad  &\text{if } b> 1 \,.
\end{aligned}
\right.     
\end{equation}
Notice that for such $W$ in $\sW^{II}(\mcR^L)$, if $b>1$, then we have $a\geq 1$ and $|b-a|\geq 1$. 


\begin{remark}
The two types of extensions in \eqref{def:extensionop_1} and \eqref{def:extensionop_2} work for $\beta\in [0,d)$ and $\beta\in(d,d+p)$ respectively. Notice that if $u\in C(\overline{\mcR_1})$, then the extension in \eqref{def:extensionop_2} gives a continuously function across the boundary $\partial\mcR^\infty$ to have necessary regularity. On the other hand, the extended function in \eqref{def:extensionop_1} is discontinuous across $\partial\mcR^\infty$.
Such an extension is fine in this case, since  $\cS_1^{\,\beta}(\mcR_{-1}^\infty)$ is equivalent to $L^p(\mcR_{-1}^\infty)$
for $\beta\in [0,d)$
, and it accepts discontinuous functions. 
We also note that the map in \eqref{def:map_M2} characterizes two regimes --  the ``fractional regime'', where any cube in $\mcR^1$ is mapped to its symmetric reflection in $\mcR_{-1}$, and the ``classical regime'', where any cube in $\mcR^\infty\backslash\mcR^1$ is mapped to a (hyper)rectangle in $\mcR_{-1}$. Related discussions on extension operators for the fractional and classical Sobolev spaces using Whitney decompositions can be found in \cite{dyda2019function} and \cite{koskela2017traces}.
\end{remark}

\begin{thm}
\label{thm:extensionL}
For any $L=2^m$ ($m\in \Z_+$), let $E^L$ be the extension operator defined in \eqref{def:extensionop_1} for $\beta\in[0,d)$ or in \eqref{def:extensionop_2} for $\beta\in(d,d+p)$, then  $E^L:\cT_{1}^{\,\beta}(\mcR_{-1}) \to \cS_1^{\,\beta}(\mcR_{-1}^\infty)$ and 
\begin{align}
 \label{eq:EL1}&\|E^L u\|^p_{L^p(\mcR_{-1}^\infty)} \leq  C L\| u\|^p_{L^p(\mcR_{-1})}   \\
 \label{eq:EL2}& |E^L u|^p_{\cS_{1}^{\,\beta}(\mcR_{-1}^\infty)} \leq C \left( L^{-(p-1)}\| u\|^p_{L^p(\mcR_{-1})} +|\beta-d|^{-1} |u|^p_{\cT_1^{\,\beta}(\mcR_{-1})} \right)
\end{align}
where $C$ is a constant independent of $L$, $\beta$ and $u\in \cT_1^{\,\beta}(\mcR_{-1})$. 
\end{thm}
\begin{proof}
 We first take the case $\beta\in [0,d)$, where $E^L$ is defined in \eqref{def:extensionop_1}. Notice that by the construction of $\{\phi_W^I\}_{W\in \sW^{I}(\mcR^L)}$ the support of each $\phi_W^I$ overlaps with only a finite number of the supports of other function.
  For the $L^p$ estimate, we have 
\[
\begin{split}
\int_{\mcR^\infty} |E^L u|^p &= \int_{\mcR^\infty} \left| \sum_{W\in \sW^{I}(\mcR^L)} a_W^I \phi_W^I(\bm x) \right|^p d\bm x
\lesssim \int_{\mcR^\infty} \sum_{W\in\sW^{I}(\mcR^L)} \left|a_W^I\right|^p | \phi_W^I |^p(\bm x) d\bm x \\
& \lesssim  \sum_{W\in\sW^{I}(\mcR^L)} |W| \left| a_W^I\right|^p \leq 
\sum_{W\in\sW^{I}(\mcR^L)}  \frac{|W|}{|\sM_1(W)|} \int_{\sM_1(W)} |u(\bm x)|^p d\bm x \\
& \lesssim \sum_{k=0}^{m} 2^k \int_{\mcR_{-1}} |u(\bm x)|^p  d\bm x \lesssim 2^{m} \int_{\mcR_{-1}} |u(\bm x)|^p  d\bm x \lesssim L \| u\|^p_{L^p(\mcR_{-1})}\,.
\end{split}
\]
Thus \eqref{eq:EL1} is true. 

Now to estimate $|E^L u|_{\cS_1^{\,\beta}(\mcR_{-1}^\infty)}$, we first note that 
\[
\int_{\mcR_{-1}}\int_{\mcR_{-1}} \ga^{\,\beta}_1(|\bm y -\bm x|) |E^L u(\bm y) -E^L u(\bm x)|^p d\bm y d\bm x = \int_{\mcR_{-1}}\int_{\mcR_{-1}} \ga^{\,\beta}_1(|\bm y -\bm x|) |u(\bm y) -u(\bm x)|^p d\bm y d\bm x \lesssim |u|^p_{\cT_{1}^{\,\beta}(\mcR_{-1})}\,.
\]
So we only need to estimate 
\begin{equation}\label{eqn:I1}
I.1=\int_{\mcR_{-1}} \int_{\mcR^\infty}  \ga^{\,\beta}_1(|\bm y -\bm x|) |E^L u(\bm y) -E^L u(\bm x)|^p  d\bm y d\bm x \,,  \end{equation}
and 
\begin{equation}\label{eqn:I2}
I.2=\int_{\mcR^\infty} \int_{\mcR^\infty}  \ga^{\,\beta}_1(|\bm y -\bm x|) |E^L u(\bm y) -E^L u(\bm x)|^p  d\bm y d\bm x \,.\end{equation}

Notice that for any $\bm y \in W_1 \in \sW^{I}(\mcR^L)$, 
\[
\begin{split}
E^L u(\bm y) - u(\bm x) = &\sum_{W\in \sW^I(\mcR^L)} a_W^I \phi_W^I(\bm y) - u(\bm x) = \sum_{W\in \sN(W_1)\cap\sW^I(\mcR^L)} \left(a_W^I-a_{W_1}^I\right) \phi_W^I(\bm y) + \left(a_{W_1}^I- u(\bm x)\right)\,,    
\end{split}
\]
where $\sN(W_1)\subset \sW^{I}(\mcR^\infty)$ denotes the collections of all the cubes that have nontrivial overlaps with the $l(W_1)/4$-neighborhood of $W_1$.
We then have the estimate 
\[
\begin{split}
I.1& = \sum_{W_1\in \overline{\sW_0}} \int_{\mcR_{-1}} \int_{W_1}  \ga^{\,\beta}_1(|\bm y -\bm x|) |E^L u(\bm y) -E^L u(\bm x)|^p  d\bm y d\bm x \\ 
& \lesssim \sum_{W_1\in \overline{\sW_0}}\sum_{W\in \sN(W_1)} \int_{\mcR_{-1}} \int_{W_1}  \ga^{\,\beta}_1 (|\bm y -\bm x|)  \left| a_{W}^I-a_{W_1}^I\right|^p d\bm y d\bm x +  \sum_{W_1\in \overline{\sW_0}}\int_{\mcR_{-1}} \int_{W_1}  \ga^{\,\beta}_1(|\bm y -\bm x|)  \left| a_{W_1}^I- u(\bm x)\right|^p d\bm y d\bm x \\
& =: I.1.a + I.1.b \,.
\end{split}
\]
We now estimate $I.1.a$. Note that
\begin{equation} \label{eq:estimate_a}
\begin{split}
a_{W}^I-a_{W_1}^I &=\frac{1}{|\sM_1(W)|} \int_{\sM_1(W)} u(\by') d\by'  - \frac{1}{|\sM_1(W_1)|} \int_{\sM_1(W_1)} u(\bx') d\bx' \\
&=\frac{1}{|\sM_1(W)||\sM_1(W_1)|} \int_{\sM_1(W_1)}\int_{\sM_1(W)}  (u(\by') -  u(\bx')) d\by' d\bx' \,,
\end{split}    
\end{equation}
we then have
\[
\begin{split}
\left| a_{W}^I-a_{W_1}^I\right|^p
&\leq \frac{1}{|\sM_1(W)||\sM_1(W_1)|}\int_{\sM_1(W_1)}\int_{\sM_1(W)}  |u(\by') -  u(\bx')|^p d\by' d\bx' \\
&\lesssim \frac{1}{|\sM_1(W)||\sM_1(W_1)|}\int_{\sM_1(W_1)}\int_{\sM_1(W)}  \frac{|u(\by') -  u(\bx')|^p}{(|\by'-\bx'|\vee 1)^{d+p-2} (|\by' -\bx'|\wedge 1)^{\,\beta}} d\by' d\bx'\,,
\end{split}
\]
where we have used $|\by'-\bx'|\lesssim 1$ in the last inequality as a result of $W_1\in \overline{\sW_0}$ and $W\in\sN(W_1)$. Notice that $\ga^{\,\beta}_1(|\by-\bx|) =C_{d,p,\,\beta}1_{\{|\by-\bx|<1\}}/|\by-\bx|^{\,\beta}$ by \eqref{def:frackernel}, so
\[
\begin{split}
I.1.a\lesssim \sum_{W_1\in \overline{\sW_0}}\sum_{W\in \sN(W_1)}& \left(\int_{\mcR_{-1}} \int_{W_1}  \frac{C_{d,p,\,\beta}1_{\{|\by-\bx|<1\}}}{ |\by-\bx|^{\,\beta}} \frac{1}{|\sM_1(W)||\sM_1(W_1)|} d\by d\bx \right) \\
& \cdot \left(\int_{\sM_1(W_1)}\int_{\sM_1(W)}  \frac{|u(\by') -  u(\bx')|^p}{(|\by'-\bx'|\vee 1)^{d+p-2} (|\by' -\bx'|\wedge1)^{\,\beta}} d\by' d\bx'\right) \,.   
\end{split}
\]
For any $\by\in W_1\in \overline{\sW_0}$,
 \[
 \int_{\mcR_{-1}} \frac{C_{d,p,\,\beta}1_{\{|\by-\bx|<1\}}}{|\by-\bx|^{\,\beta}} d\bx \leq \int_{|\bm s|<1} \frac{C_{d,p,\,\beta}}{|\bm s|^{\,\beta}} d\bm s \lesssim \frac{1}{d-\beta} \,.
 \]
 Notice that $1/(d-\beta)$ is the constant that appears in \eqref{eq:EL2} and it blows up as $\beta\to d$ so we have to take a fixed $\beta\in[0,d)$. 
Therefore, 
\[
\int_{\mcR_{-1}} \int_{W_1}  \frac{C_{d,p,\,\beta}1_{\{|\by-\bx|<1\}}}{ |\by-\bx|^{\,\beta}} \frac{1}{|\sM_1(W)||\sM_1(W_1)|} d\by d\bx \lesssim  \frac{|\beta-d|^{-1}|W_1|}{|\sM_1(W)||\sM_1(W_1)|}  \lesssim 
|\beta-d|^{-1}\,,
\]
where we have used the fact that $|W_1| \approx |\sM_1(W)|\approx |\sM_1(W_1)| \approx 1$ for any $W_1\in \overline{\sW_0}$ and $W\in \sN(W_1)$. So 
\[
I.1.a \lesssim |\beta-d|^{-1} \sum_{W_1\in \overline{\sW_0}}\sum_{W\in \sN(W_1)}\int_{\sM_1(W_1)}\int_{\sM_1(W)}  \frac{|u(\by') -  u(\bx')|^p}{(|\by'-\bx'|\vee1)^{d+p-2} (|\by' -\bx'|\wedge1)^{\,\beta}} d\by' d\bx' \lesssim |\beta-d|^{-1}|u|^p_{\cT_1^{\,\beta}(\mcR_{-1})}\,.
\]
Now since
\[
\begin{split}
\left| a_{W_1}^I- u(\bm x)\right|^p &= \left| \frac{1}{|\sM_1(W_1)|} \int_{\sM_1(W_1)} (u(\by')- u(\bx)) d\by' \right|^p \leq \frac{1}{|\sM_1(W_1)|} \int_{\sM_1(W_1)} |u(\by')- u(\bx)|^p d\by'  \\
&\lesssim \frac{1}{|\sM_1(W_1)|}\int_{\sM_1(W_1)}  \frac{|u(\by')- u(\bx)|^p}{(|\by'-\bx|\vee1)^{d+p-2} (|\by' -\bx|\wedge1)^{\,\beta}} d\by'\,,
\end{split}
\]
as a result of the fact that $\by\in W_1$ and $|\by-\bx|\lesssim1$, we have
\[
\begin{split}
I.1.b &\lesssim \sum_{W_1\in \overline{\sW_0}} \int_{\mcR_{-1}} \left( \int_{W_1}  \frac{1_{\{|\by-\bx|<1\}}}{ |\by-\bx|^{\,\beta}} \frac{1}{|\sM_1(W_1)|} d\by \int_{\sM_1(W_1)}  \frac{|u(\by') -  u(\bx)|^p}{(|\by'-\bx|\vee1)^{d+p-2} (|\by' -\bx|\wedge1)^{\,\beta}} d\by'  \right)d\bx \\
&\lesssim  |\beta-d|^{-1}\sum_{W_1\in \overline{\sW_0}} \int_{\mcR_{-1}} \int_{\sM_1(W_1)} \frac{|u(\by') -  u(\bx)|^p}{(|\by'-\bx|\vee1)^{d+p-2} (|\by' -\bx|\wedge1)^{\,\beta}} d\by'  d\bx \lesssim |\beta-d|^{-1}|u|^p_{\cT_1^{\,\beta}(\mcR_{-1})} \,.
\end{split}
\]

Together we have shown that $I.1$ in \eqref{eqn:I1} is bounded by $|u|^p_{\cT_1^{\,\beta}(\mcR_{-1})}$.

To estimate $I.2$ in \eqref{eqn:I2}, we first note that 
\begin{equation}\label{eq:I.2}
 I.2 = \sum_{W_1\in\sW^I(\mcR^\infty)} \sum_{W_2\in \sN(W_1)}\int_{W_1}\int_{W_2} \ga^{\,\beta}_1(|\bm y -\bm x|) |E^L u(\bm y) -E^L u(\bm x)|^p  d\by d\bx\,,
\end{equation}
since all cubes in $\sW^I(\mcR^\infty)$ have length greater than or equal to $1$ and $\ga^{\,\beta}_1(|\by-\bx|)=0$ if $|\by-\bx|>1$. 
Now suppose
 $\xb\in W_1 \in \sW_k$ for $k\leq -m+1$, and $y\in W_2\in \sN(W_1)$ then 
\[
E^L u(\by)-E^L u(\bx) = \sum_{W\in (\sN(W_1)\cup \sN(W_2))\cap \sW^I(\mcR^L)} \left[ \left( a_W^I -a_{W_1}^I\right) \left(\phi_W^I(\by) - \phi_W^I(\bx)\right) + a_{W_1}^I \left(\phi_W^I(\by) - \phi_W^I(\bx)\right)  \right],
\]
On the other hand if
$\bx\in W_1\in \sW^{I}(\mcR^L)\backslash \sW_{-m+1}$ and $\yb\in W_2 \in \sN(W_1)$, then we know that both $\xb$ and $\yb$ are in $\mcR^L$, and therefore
\[
\sum_{W\in \sW^I(\mcR^L)} \phi_W^I(\by) = \sum_{W\in \sW^I(\mcR^L)} \phi_W^I(\bx)= 1.
\]
In turn, we have
\[
E^L u(\by)-E^L u(\bx) = \sum_{W\in \sW^{I}(\mcR^L)} a_W^I \left(\phi_W^I(\by) - \phi_W^I(\bx)\right) = \sum_{W\in (\sN(W_1)\cup \sN(W_2))\cap \sW^I(\mcR^L)} \left( a_W^I -a_{W_1}^I\right) \left(\phi_W^I(\by) - \phi_W^I(\bx)\right) \,. 
\]
Taking into account the two cases, we can show
\[
I.2 \lesssim I.2.a + I.2.b
\]
where 
\[
I.2.a  = \sum_{W_1\in\sW^I(\mcR^\infty)} \sum_{W_2\in \sN(W_1)} \int_{W_1}\int_{W_2} \ga^{\,\beta}_1(|\bm y -\bm x|)\left| \sum_{W\in (\sN(W_1)\cup \sN(W_2))\cap \sW^I(\mcR^L)} \left( a_W^I -a_{W_1}^I\right) \left(\phi_W^I(\by) - \phi_W^I(\bx)\right)\right|^p   d\by d\bx ,
\]
and 
\[
I.2.b = \sum_{k\leq -m+1}\sum_{W_1\in \sW_k} \sum_{W_2\in \sN(W_1)} \int_{W_1}\int_{W_2} \ga^{\,\beta}_1(|\bm y -\bm x|) \left| \sum_{W\in (\sN(W_1)\cup \sN(W_2))\cap \sW^I(\mcR^L)}a_{W_1}^I \left(\phi_W^I(\by) - \phi_W^I(\bx)\right) \right|^p d\by d\bx .
\]
We first estimate $I.2.b$. Since the number of sets in $\sN(W_1)\cup \sN(W_2)$ is uniformly bounded by a constant for any $W_2\in \sN(W_1)$, and Lip$\phi_W^I\lesssim 1/l(W)$, we have 
\[
\left|\sum_{W\in (\sN(W_1)\cup \sN(W_2))\cap \sW^I(\mcR^L)} 
a_{W_1}^I \left(\phi_W^I(\by) - \phi_W^I(\bx)\right)\right|^p \lesssim  \left| a_{W_1}^I\right|^p \frac{|\yb-\xb|^p}{|l(W_1)|^p} \lesssim \frac{|\yb-\xb|^p}{|\sM_1(W_1)| |l(W_1)|^p}  \int_{\sM_1(W_1)} |u(\zb)|^p d\zb ,
\] 
where we have also used $l(W_1)\approx l(W_2) \approx l(W)$ for any $W\in \sN(W_1)\cup \sN(W_2)$ and $W_2\in \sN(W_1)$.
Therefore 
\[
\begin{split}
I.2.b &\lesssim
\sum_{k\leq -m+1} \sum_{W_1 \in \sW_k} \sum_{W_2\in \sN(W_1)} \int_{W_1} \int_{W_2}  \frac{\ga^{\,\beta}_1(|\bm y -\bm x|) |\yb-\xb|^p d\yb d\xb }{|\sM_1(W_1)| |l(W_1)|^p}  \int_{\sM_1(W_1)} |u(\zb)|^p d\zb  \\
& \lesssim  \sum_{k\leq -m+1} \sum_{W_1 \in \sW_k} \frac{|W_1|}{|\sM_1(W_1)| |l(W_1)|^p } \int_{\sM_1(W_1)} |u(\zb)|^p d\zb \\
&\lesssim \sum_{k\leq -m+1} 2^{k(p-1)} \int_{\mcR_{-1}} |u(\zb)|^p d\zb \lesssim L^{-(p-1)} \| u\|^p_{L^p(\mcR_{-1})}.
\end{split}
\]
For I.2.a, we first notice that
\[
\left| \sum_{W\in (\sN(W_1)\cup \sN(W_2))\cap \sW^I(\mcR^L)} \left( a_W^I -a_{W_1}^I\right)  \left(\phi_W^I(\by) - \phi_W^I(\bx)\right) \right|^p \lesssim 
\sum_{W\in \sN(W_1)\cup \sN(W_2)} \left| a_W^I -a_{W_1}^I\right|^p  \frac{|\by -\bx|^p}{l(W)^p}\,, 
\]
where we have used the fact that $\sN(W_1)\cup \sN(W_2)$ contains only finite number of sets and  Lip$\phi_W^I\lesssim 1/l(W)$. 
 Now from \eqref{eq:estimate_a}, we obtain 
\[
\begin{split}
\left| a_{W}^I-a_{W_1}^I\right|^p
&\leq \frac{1}{|\sM_1(W)||\sM_1(W_1)|}\int_{\sM_1(W_1)}\int_{\sM_1(W)}  |u(\by') -  u(\bx')|^p d\by' d\bx' \\
&\lesssim \frac{1}{|\sM_1(W)||\sM_1(W_1)|}\int_{\sM_1(W_1)}\int_{\sM_1(W)}  \frac{|u(\by') -  u(\bx')|^p}{ (|\by' -\bx'|\wedge1)^{\,\beta}} d\by' d\bx'\,.
\end{split}
\]
Notice also that 
\[
\int_{W_1} \int_{W_2}\frac{C_{d,p,\,\beta}1_{\{|\by-\bx|<1\}}}{|\by-\bx|^{\beta-p}} d\by d\bx \leq \int_{W_1} \int_{\R^d}\frac{C_{d,p,\,\beta} 1_{\{|\by-\bx|<1\}}}{|\by-\bx|^{\beta-p}} d\by d\bx  \lesssim |W_1| \,.
\]
So we have 
\[
\begin{split}
I.2.a&\lesssim \sum_{W_1\in\sW^I(\mcR^\infty)} \sum_{W_2\in\sN(W_1)} \sum_{W\in\sN(W_1)\cup\sN(W_2)}\left(\int_{W_1} \int_{W_2}  \frac{C_{d,p,\,\beta}1_{\{|\by-\bx|<1\}}}{ |\by-\bx|^{\beta-p} l(W)^p} \frac{1}{|\sM_1(W)||\sM_1(W_1)|} d\by d\bx \right) \\
&\hspace{7cm} \cdot \left(  \int_{\sM_1(W_1)}\int_{\sM_1(W)}  \frac{|u(\by') -  u(\bx')|^p}{ (|\by' -\bx'|\wedge1)^{\,\beta}} d\by' d\bx' \right)    \\
&\lesssim \sum_{W_1\in\sW^I(\mcR^\infty)} \sum_{W_2\in\sN(W_1)}  \sum_{W\in\sN(W_1)\cup\sN(W_2)}\frac{|W_1|}{l(W)^p|\sM_1(W)||\sM_1(W_1)|} \int_{\sM_1(W_1)}\int_{\sM_1(W)}  \frac{|u(\by') -  u(\bx')|^p}{ (|\by' -\bx'|\wedge1)^{\,\beta}} d\by' d\bx' \\
&\lesssim \sum_{W_1\in\sW^I(\mcR^\infty)} \sum_{W_2\in\sN(W_1)}  \sum_{W\in\sN(W_1)\cup\sN(W_2)} \frac{l(W_1)^{d+p-2}}{l(W_1)^p l(W_1)^{2d-2}} \int_{\sM_1(W_1)}\int_{\sM_1(W)} \frac{|u(\by') -  u(\bx')|^p}{ (|\by' -\bx'|\wedge1)^{\,\beta}} d\by' d\bx' \\
&\lesssim  \left(\sum_{k\leq 0,k\in\Z} \sum_{W_1\in \sW_k}+\sum_{W_1\in \overline{\sW_0},\, k=0}\right) \sum_{W_2\in\sN(W_1)}  \sum_{W\in\sN(W_1)\cup\sN(W_2)} ( 2^{-k})^{-(d+p-2)}   \int_{\sM_1(W_1)}\int_{\sM_1(W)} \frac{|u(\by') -  u(\bx')|^p}{ (|\by' -\bx'|\wedge1)^{\,\beta}} d\by' d\bx' \\
&\lesssim  \left(\sum_{k\leq 0,k\in\Z} \sum_{W_1\in \sW_k}+\sum_{W_1\in \overline{\sW_0},\, k=0}\right) 2^{k(d+p-2)}   \int_{\sM_1(W_1)}\int_{B(\bx', C\cdot 2^{-k})\cap \mcR_{-1}} \frac{|u(\by') -  u(\bx')|^p}{ (|\by' -\bx'|\wedge1)^{\,\beta}} d\by' d\bx' \,,
\end{split}  
\]
for some $C>0$, and $B(\bx', r)$ denotes the ball of radius $r$ centered at $\bx'$. Now chose $C_1=C+1$, then $(C_1-C) \cdot 2^{-k} \geq 1$ for all $k\leq 0$. Therefore
\[
\begin{split}
 I.2.a &\lesssim \left(\sum_{k\leq 0,k\in\Z} \sum_{W_1\in \sW_k}+\sum_{W_1\in \overline{\sW_0},\, k=0}\right)  \int_{C_1\cdot 2^{-k}}^{C_1\cdot 2^{-k+1}}\int_{\sM_1(W_1)}\int_{B(\bx', C\cdot 2^{-k})\cap \mcR_{-1}} \frac{|u(\by') -  u(\bx')|^p}{ (|\by' -\bx'|\wedge1)^{\,\beta}} d\by' d\bx'  \frac{dt}{t^{d+p-1}} \\   
 &\lesssim \left(\sum_{k\leq 0,k\in\Z} \sum_{W_1\in \sW_k}+\sum_{W_1\in \overline{\sW_0},\, k=0}\right)  \int_{C_1\cdot 2^{-k}}^{C_1\cdot 2^{-k+1}}\int_{\sM_1(W_1)}\int_{B(\bx', C_1\cdot 2^{-k} -1 )\cap \mcR_{-1}} \frac{|u(\by') -  u(\bx')|^p}{ (|\by' -\bx'|\wedge1)^{\,\beta}} d\by' d\bx'  \frac{dt}{t^{d+p-1}} \\  
 &\lesssim \left(\sum_{k\leq 0,k\in\Z} \sum_{W_1\in \sW_k}+\sum_{W_1\in \overline{\sW_0},\, k=0}\right)  \int_{C_1\cdot 2^{-k}}^{C_1\cdot 2^{-k+1}}\int_{\sM_1(W_1)}\int_{B(\bx', t -1 )\cap \mcR_{-1}} \frac{|u(\by') -  u(\bx')|^p}{ (|\by' -\bx'|\wedge1)^{\,\beta}} d\by' d\bx'  \frac{dt}{t^{d+p-1}} \\
 &\lesssim  \int_{0}^{\infty}\int_{\mcR_{-1}}\int_{B(\bx', t -1)\cap \mcR_{-1}} \frac{|u(\by') -  u(\bx')|^p}{ (|\by' -\bx'|\wedge
 1)^{\,\beta}} d\by' d\bx' \frac{dt}{t^{d+p-1}}   \\
 &\leq   \int_{\mcR_{-1}}\int_{ \mcR_{-1}} \frac{|u(\by') -  u(\bx')|^p}{ (|\by' -\bx'|\wedge1)^{\,\beta}} \int_{|\by'-\bx'|<t-1} \frac{dt}{t^{d+p-1}} d\by' d\bx' \\
 &\lesssim  \int_{\mcR_{-1}}\int_{ \mcR_{-1}} \frac{|u(\by') -  u(\bx')|^p}{ (|\by'-\bx'|+1)^{d+p-2}(|\by' -\bx'|\wedge1)^{\,\beta}}  d\by' d\bx' \lesssim |u|^p_{\cT_1^{\,\beta}(\mcR_{-1})}\,.
\end{split}
\]
Together we have shown \eqref{eq:EL2} for the case $\beta\in [0,d) $.

Now for $\beta\in (d,d+p)$, we take the extension operator $E^L$ defined in \eqref{def:extensionop_2}. The $L^p$ estimate of $E^L u$ can be shown similarly as in the first case. For the estimate of $|E^L u|_{\cS_1^{\,\beta}(\mcR_{-1}^\infty)}$, similar to the first case considered earlier, it is not hard to see that we only need to estimate \begin{equation}
II.1=\int_{\mcR_{-1}} \int_{\mcR^\infty}  \ga^{\,\beta}_1(|\bm y -\bm x|) |E^L u(\bm y) -E^L u(\bm x)|^p  d\bm y d\bm x \,,  \end{equation}
and 
\begin{equation}
II.2=\int_{\mcR^\infty} \int_{\mcR^\infty}  \ga^{\,\beta}_1(|\bm y -\bm x|) |E^L u(\bm y) -E^L u(\bm x)|^p  d\bm y d\bm x \,.\end{equation}
Similar to the first case, we can split $II.1$ into two parts. 
\[
\begin{split}
II.1& = \sum_{k=1}^\infty\sum_{W_1\in\sW_k} \int_{\mcR_{-1}} \int_{W_1}  \ga^{\,\beta}_1(|\bm y -\bm x|) |E^L u(\bm y) -E^L u(\bm x)|^p  d\bm y d\bm x \\ 
& \lesssim \sum_{k=1}^\infty\sum_{W_1\in\sW_k} \sum_{W\in \sN(W_1)} \int_{\mcR_{-1}} \int_{W_1}  \ga^{\,\beta}_1(|\bm y -\bm x|)  \left|a_{W}^{II}-a_{W_1}^{II}\right|^p d\bm y d\bm x +  \sum_{k=1}^\infty\sum_{W_1\in\sW_k} \int_{\mcR_{-1}} \int_{W_1}  \ga^{\,\beta}_1(|\bm y -\bm x|)  \left|a_{W_1}^{II}- u(\bm x)\right|^p d\bm y d\bm x \\
& =: II.1.a + II.1.b \,.
\end{split}
\]
From a similar equation to \eqref{eq:estimate_a}, we have 
\begin{equation}\label{eq:estimate_aII}
 \left|a_{W}^{II}-a_{W_1}^{II}\right|^p
\lesssim \frac{ l(W_1)^{\,\beta}}{|\sM_2(W)||\sM_2(W_1)|}\int_{\sM_2(W_1)}\int_{\sM_2(W)}  \frac{|u(\by') -  u(\bx')|^p}{(|\by'-\bx'|\vee1)^{d+p-2} (|\by' -\bx'|\wedge1)^{\,\beta}} d\by' d\bx'\,,   
\end{equation}
where we have used $W\in \sN(W_1)$ and $l(W_1)\leq 1$.
So 
\[
\begin{split}
II.1.a\lesssim \sum_{k=1}^\infty\sum_{W_1\in\sW_k}\sum_{W\in \sN(W_1)}& \left(\int_{\mcR_{-1}} \int_{W_1}  \frac{C_{d,p,\,\beta}1_{\{|\by-\bx|<1\}}}{ |\by-\bx|^{\,\beta}} \frac{ l(W_1)^{\,\beta}}{|\sM_2(W)||\sM_2(W_1)|} d\by d\bx \right) \\
& \cdot \left(\int_{\sM_2(W_1)}\int_{\sM_2(W)}  \frac{|u(\by') -  u(\bx')|^p}{(|\by'-\bx'|\vee 1)^{d+p-2} (|\by' -\bx'|\wedge1)^{\,\beta}} d\by' d\bx'\right) \,.   
\end{split}
\]
For any $\by\in W_1 \in \sW_k$ and  $\beta\in (d,d+p)$,
\[
\int_{\mcR_{-1}} \frac{C_{d,p,\,\beta}1_{\{|\by-\bx|<1\}}}{|\by-\bx|^{\,\beta}} d\bx  \lesssim \frac{1}{\beta-d} l(W_1)^{d-\beta} \,. 
\] 
Notice again that the constant $1/(\beta-d)$ blows up as $\beta\to d$. Then
\[
\int_{\mcR_{-1}} \int_{W_1}  \frac{C_{d,p,\,\beta}1_{\{|\by-\bx|<1\}}}{ |\by-\bx|^{\,\beta}} \frac{ l(W_1)^{\,\beta}}{|\sM_2(W)||\sM_2(W_1)|} d\by d\bx \lesssim |\beta-d|^{-1} l(W_1)^{d-\beta}\cdot \frac{ l(W_1)^{\,\beta}}{|\sM_2(W)||\sM_2(W_1)|} |W_1| \lesssim |\beta-d|^{-1}\,,
\]
where we have used the fact that $|\sM_2(W)|\approx|\sM_2(W_1)|\approx |W_1| \approx l(W_1)^{d+p-2}$. Using this estimate, one can show 
\[
II.1.a\lesssim |\beta-d|^{-1} |u|^p_{\cT_1^{\,\beta}(\mcR_{-1})}.
\]  
To estimate $II.1.b$, we first define a decomposition of $\mcR_{-1}$
\begin{equation}
\sW(\mcR_{-1}) =  \left\{\sM_2(W): W\in \sW_k, k\in \Z_+ \right\}.    \end{equation}
Then 
\begin{equation}
\label{eq:estimate_II1b}
II.1.b = \sum_{W_2\in \sW(\mcR_{-1})}\sum_{k=1}^\infty\sum_{W_1\in\sW_k} \int_{W_2} \int_{W_1}  \ga^{\,\beta}_1(|\bm y -\bm x|)  \left|a_{W_1}^{II}- u(\bm x)\right|^p d\bm y d\bm x \,.   
\end{equation}
Notice that for any $\bx\in W_2 \in \sW(\mcR_{-1})$ and $\by'\in \sM_2(W_1)$ for $W_1\in \sW_k, k\in \Z_+$, we have 
\[
|\by' -\bx| \lesssim \dist(W_1,W_2) + l(W_1) + l(W_2) \lesssim \dist(W_1,W_2)\,. 
\]
Moreover, we must have $\dist(W_1, W_2)< 1$ for the double integral in \eqref{eq:estimate_II1b} to be non-zero. 
We thus have the estimate  
\[
 \left|a_{W_1}^{II}- u(\bm x)\right|^p \leq \frac{1}{|\sM_2(W_1)|} \int_{\sM_2(W_1)} |u(\by') -u(\bx)|^p d\by'    \lesssim \frac{ \dist(W_1, W_2)^{\,\beta}}{|\sM_2(W_1)|}\int_{\sM_2(W_1)} \frac{|u(\by') -u(\bx)|^p}{(|\by' -\bx|\vee1)^{d+p-2} (|\by' -\bx|\wedge1)^{\,\beta} } d\by' \,.
\]
Therefore 
\[
\begin{split}
&II.1.b \lesssim \sum_{W_2\in \sW(\mcR_{-1})} \sum_{W_1\in \sW_k,\,k\in \Z_+} \int_{W_2} \left( \int_{W_1}  \frac{C_{d,p,\,\beta}1_{\{|\by-\bx|<1\}}}{ |\by-\bx|^{\,\beta}} \frac{ \dist(W_1, W_2)^{\,\beta}}{|\sM_2(W_1)|} d\by \int_{\sM_2(W_1)}  \frac{|u(\by') -  u(\bx)|^p}{(|\by'-\bx|\vee1)^{d+p-2} (|\by' -\bx|\wedge1)^{\,\beta}} d\by'  \right)d\bx \\
&\lesssim \sum_{W_2\in \sW(\mcR_{-1})} \sum_{W_1\in \sW_k,\,k\in \Z_+}  \int_{W_2} \left(  \int_{W_1}  \frac{C_{d,p,\,\beta}}{\dist(W_1, W_2)^{\,\beta}} \frac{\dist(W_1, W_2)^{\,\beta}}{|\sM_2(W_1)|}  d\by \int_{\sM_2(W_1)}  \frac{|u(\by') -  u(\bx)|^p}{(|\by'-\bx|\vee1)^{d+p-2} (|\by' -\bx|\wedge1)^{\,\beta}} d\by'  \right)d\bx\\
&\lesssim  \sum_{W_2\in \sW(\mcR_{-1})} \sum_{W_1\in \sW_k,\,k\in \Z_+}  \int_{W_2} \int_{{\sM_2(W_1)}} \frac{|u(\by') -  u(\bx)|^p}{(|\by'-\bx|\vee1)^{d+p-2} (|\by' -\bx|\wedge1)^{\,\beta}} d\by'  d\bx \lesssim |u|^p_{\cT_1^{\,\beta}(\mcR_{-1})} \,.
\end{split}
\]
Now for $II.2$, we can first write   
\begin{equation}\label{eq:II.2}
 II.2 = \sum_{W_1\in\sW^{II}(\mcR^\infty)} \sum_{W_2\in\sW^{II}(\mcR^\infty)}\int_{W_1}\int_{W_2} \ga^{\,\beta}_1(|\bm y -\bm x|) |E^L u(\bm y) -E^L u(\bm x)|^p  d\by d\bx\,.  
\end{equation}
Observe that 
\[
E^L u(\by)-E^L u(\bx) = \sum_{W\in (\sN(W_1)\cup \sN(W_2))\cap \sW^{II}(\mcR^L)} \left[ \left( a_W^{II} -a_{W_1}^{II}\right) \left(\phi_W^{II}(\by) - \phi_W^{II}(\bx)\right) + a_{W_1}^{II} \left(\phi_W^{II}(\by) - \phi_W^{II}(\bx)\right)  \right], 
\]
where the second part in the above equation is only nonzero for $\xb\in W_1 \in  \sW_{k}$ for $k\leq -m+1$ (and therefore $\yb\in W_2 \in \sN(W_1)$ in this case because of the nonlocal interaction length).
Similarly as before, we have
\[
II.2 \lesssim II.2.a + II.2.b
\]
where 
\[
II.2.a  = \sum_{W_1\in\sW^{II}(\mcR^\infty)} \sum_{W_2\in \sW^{II}(\mcR^\infty)} \int_{W_1}\int_{W_2} \ga^{\,\beta}_1(|\bm y -\bm x|)\left| \sum_{W\in (\sN(W_1)\cup \sN(W_2))\cap \sW^{II}(\mcR^L)} \left( a_W^{II} -a_{W_1}^{II}\right) \left(\phi_W^{II}(\by) - \phi_W^{II}(\bx)\right)\right|^p   d\by d\bx 
\]
and 
\[
II.2.b = \sum_{k\leq -m+1}\sum_{W_1\in \sW_k} \sum_{W_2\in \sN(W_1)} \int_{W_1}\int_{W_2} \ga^{\,\beta}_1(|\bm y -\bm x|) \left| \sum_{W\in (\sN(W_1)\cup \sN(W_2))\cap \sW^{II}(\mcR^L)}a_{W_1}^{II} \left(\phi_W^{II}(\by) - \phi_W^{II}(\bx)\right) \right|^p d\by d\bx. 
\]
It is easy to see that $II.2.b$ can be estimated similarly as $I.2.b$ and we have $II.2.b\lesssim  L^{-(p-1)} \| u\|^p_{L^p(\mcR_{-1})}$.
Now for the estimate of $II.2.a$, we have two different cases where $W_2\in \sN(W_1)$ and $W_2\notin \sN(W_1)$. For the case $W_2\in \sN(W_1)$, the estimate follows similarly to the estimate of $I.2.a$ which is omitted here. For the case $W_2\notin \sN(W_1)$, 
we proceed by noticing that if $l(W_1)\geq 2$ and $W_2\notin \sN(W_1)$, then we must have $\dist(W_1, W_2) \geq 1$ so that the double integral in \eqref{eq:II.2} becomes zero. Therefore we only need to consider $l(W_1)\leq 1$ (i.e., $W_1\in \sW_k$ for $k\geq 0$) in this case. Notice that
\[
\begin{split}
\left| \sum_{W\in (\sN(W_1)\cup \sN(W_2))\cap \sW^{II}(\mcR^L)}\left( a_W^{II} -a_{W_1}^{II}\right) \left(\phi_W^{II}(\by) - \phi_W^{II}(\bx)\right)\right|^p &\lesssim \sum_{W\in \sN(W_1)\cup \sN(W_2)} \left| a_W^{II} -a_{W_1}^{II}\right|^p  \left|\phi_W^{II}(\by) - \phi_W^{II}(\bx)\right|^p \\
&\lesssim 
\sum_{W\in \sN(W_1)\cup \sN(W_2)} \left| a_W^{II} -a_{W_1}^{II}\right|^p \,, 
\end{split}
\]
for $\bx\in W_1$ and $\by\in W_2$. For $W\in \sN(W_1)$, we use \eqref{eq:estimate_aII} and 
\[
\sum_{W_2\notin \sN(W_1)}\int_{W_1} \int_{W_2}  \frac{C_{d,p,\,\beta}1_{\{|\by-\bx|<1\}}}{|\by-\bx|^{\,\beta}} d\by d\bx \leq \int_{W_1} \int_{|\by-\bx|\geq  l(W_1)/2}   \frac{C_{d,p,\,\beta}}{|\by-\bx|^{\,\beta}} d\by d\bx  \lesssim |\beta-d|^{-1}|W_1| l(W_1)^{d-\beta}
\]
to get 
\[
\begin{split}
& \sum_{W_1\in\sW^{II}(\mcR^\infty)} \sum_{W_2\in\sW^{II}(\mcR^\infty)\backslash\sN(W_1)}\int_{W_1}\int_{W_2} \ga^{\,\beta}_1(|\bm y -\bm x|) \sum_{W\in \sN(W_1)} \left| a_W^{II} -a_{W_1}^{II}\right|^p  d\by d\bx   \\
\lesssim & \sum_{\substack{W_1\in\sW_k \\ k\geq 0}} \sum_{W_2\in\sW^{II}(\mcR^\infty)\backslash\sN(W_1)}  \sum_{W\in \sN(W_1)} \left(\int_{W_1}\int_{W_2} \frac{C_{d,p,\,\beta}1_{\{|\by-\bx|<1\}}}{|\by-\bx|^{\,\beta}}  \frac{ l(W_1)^{\,\beta}}{|\sM_2(W)||\sM_2(W_1)|}d\by d\bx \right) \\
& \hspace{6cm}\cdot\left(\int_{\sM_2(W_1)}\int_{\sM_2(W)}  \frac{|u(\by') -  u(\bx')|^p}{(|\by'-\bx'|\vee1)^{d+p-2} (|\by' -\bx'|\wedge1)^{\,\beta}} d\by' d\bx'\right) \\
\lesssim & |\beta-d|^{-1}\sum_{\substack{W_1\in\sW_k \\ k\geq 0}}  \sum_{W\in \sN(W_1)}  \int_{\sM_2(W_1)}\int_{\sM_2(W)}  \frac{|u(\by') -  u(\bx')|^p}{(|\by'-\bx'|\vee1)^{d+p-2} (|\by' -\bx'|\wedge1)^{\,\beta}} d\by' d\bx' \lesssim |\beta-d|^{-1}|u|^p_{\cT_1^{\,\beta}(\mcR_{-1})}\,.
\end{split}
\]
On the other hand, if $W\in \sN(W_2)$, then $|W|\approx |\sM(W)|\approx|W_2|$, and we can use
\[
\begin{split}
\left| a_{W}^{II}-a_{W_1}^{II}\right|^p 
&\leq \frac{1}{|\sM_1(W)||\sM_1(W_1)|} \int_{\sM_1(W_1)}\int_{\sM_1(W)}  |u(\by') -  u(\bx')|^p d\by' d\bx' \\
&\lesssim \frac{ \dist(W_1,W_2)^{\,\beta}}{|\sM_2(W)||\sM_2(W_1)|}\int_{\sM_2(W_1)}\int_{\sM_2(W)}  \frac{|u(\by') -  u(\bx')|^p}{(|\by'-\bx'|\vee1)^{d+p-2} (|\by' -\bx'|\wedge1)^{\,\beta}} d\by' d\bx' 
\end{split}
\]
and 
\[
\int_{W_1} \int_{W_2}  \frac{1_{\{|\by-\bx|<1\}}}{|\by-\bx|^{\,\beta}} d\by d\bx \leq  \frac{|W_1| |W_2|}{\dist(W_1,W_2)^{\,\beta}}
\]
to arrive at 
\[
\begin{split}
& \sum_{W_1\in\sW^{II}(\mcR^\infty)} \sum_{W_2\in\sW^{II}(\mcR^\infty)\backslash\sN(W_1)}\int_{W_1}\int_{W_2} \ga^{\,\beta}_1(|\bm y -\bm x|) \sum_{W\in \sN(W_2)} \left| a_W^{II} -a_{W_1}^{II}\right|^p  d\by d\bx   \\
\lesssim & \sum_{\substack{W_1\in\sW_k \\ k\geq 0}} \sum_{W_2\in\sW^{II}(\mcR^\infty)\backslash\sN(W_1)}  \sum_{W\in \sN(W_2)} \left(\int_{W_1}\int_{W_2} \frac{C_{d,p,\,\beta}1_{\{|\by-\bx|<1\}}}{|\by-\bx|^{\,\beta}}  \frac{ \dist(W_1,W_2)^{\,\beta}}{|\sM_2(W)||\sM_2(W_1)|}d\by d\bx \right) \\
& \hspace{6cm}\cdot\left(\int_{\sM_2(W_1)}\int_{\sM_2(W)}  \frac{|u(\by') -  u(\bx')|^p}{(|\by'-\bx'|\vee1)^{d+p-2} (|\by' -\bx'|\wedge1)^{\,\beta}} d\by' d\bx'\right) \\
\lesssim &\sum_{\substack{W_1\in\sW_k \\ k\geq 0}} \sum_{W_2\in\sW^{II}(\mcR^\infty)\backslash\sN(W_1)}  \sum_{W\in \sN(W_2)} \int_{\sM_2(W_1)}\int_{\sM_2(W)}  \frac{|u(\by') -  u(\bx')|^p}{(|\by'-\bx'|\vee1)^{d+p-2} (|\by' -\bx'|\wedge1)^{\,\beta}} d\by' d\bx' \lesssim |u|^p_{\cT_1^{\,\beta}(\mcR_{-1})}.
\end{split}
\]
Together we have shown $II.2.a\lesssim |\beta-d|^{-1}|u|^p_{\cT_1^{\,\beta}(\mcR_{-1})}$ and as a result the estimate \eqref{eq:EL2} is proved for $\beta\in (d,d+p)$. 
\end{proof}

\noindent {\it Proof of Theorem \ref{mainthm_2}.} 
For any $u\in \cT^{\,\beta}_\del (\mcR_{-\del})$, define $v(x)=u(\del x)$, then from Lemma \ref{lem:scale} we know $|v|^p_{\cT_1^{\,\beta}(\mcR_{-1})} = \del^{p-d} |u|^p_{\cT_1^{\,\beta}(\mcR_{-\del})}$ and $\| v\|^p_{L^p(\mcR_{-1})}=\del^{-d} \| u\|^p_{L^p(\mcR_{-\del})}$. From Theorem \ref{thm:extensionL} we know that 
\[
\| E^L v\|^p_{L^p(\mcR_{-1}^\infty)} \leq C  L \| v\|^p_{L^p(\mcR_{-1})} , \text{ and } |E^L v|^p_{\cS_{1}^{\,\beta}(\mcR^\infty_{-1})} \leq C\left( L^{-(p-1)} \| v\|^p_{L^p(\mcR_{-1})}  +  |\beta-d|^{-1}|v|^p_{\cT_{1}^{\,\beta}(\mcR^\infty_{-1})} \right)
\]
for any $L=2^m$ ($m\in \Z_+$).
Now choose $m = \lceil \log_2(M/\del) \rceil$ and $L=2^m$, and define the extension operator $E: \cT^{\,\beta}_\del(\mcR_{-\del}) \to \cS^{\,\beta}_\del(\mcR_{-\del}^\infty)$ as 
\[
Eu (\bx) = (E^L v)(\bx/\del)\,.
\]
Then 
\[
\| Eu\|^p_{L^p(\mcR_{-\del}^\infty)} = \del^d \| E^L v\|^p_{L^p(\mcR_{-1}^\infty)}\leq C \del^d L \| v\|_{L^p(\mcR_{-1})} = C L \| u\|^p_{L^p(\mcR_{-\del})} \lesssim \frac{1}{\del} \| u\|^p_{L^p(\mcR_{-\del})} \,,
\]
and 
\[
\begin{split}
| Eu|^p_{\cS_{\del}^{\,\beta}(\mcR_{-\del}^\infty)} = \del^{d-p} | E^L v|^p_{\cS_{1}^{\,\beta}(\mcR_{-1}^\infty)}\leq &C \left( \del^{d-p} L^{-(p-1)} \| u\|^p_{L^p(\mcR_{-\del})} + |\beta-d|^{-1}\del^{d-p}  | v|^p_{\cT_1^{\,\beta}(\mcR_{-1})} \right) \\ \lesssim &C\left(   \frac{1}{\del} \| u\|^p_{L^p(\mcR_{-\del})}+ |\beta-d|^{-1}| u|^p_{\cT^{\,\beta}_\del(\mcR_{-\del})} \right)\,.   
\end{split}
\]
Thus Theorem \ref{mainthm_2} is shown. 
\qed 

\section{Extension to General Lipschitz Domains }\label{sec:generaldomain}

In this section we will extend the trace theorems for the half plane to a general Lipschitz domain.  
We take the strategy
to first generalize to a ``special'' Lipschitz domain before showing the fully general case.  

\subsection{Some technical lemmas}
We will first list some 
lemmas which are used to show the transformations 
from the special Lipschitz domain to the half-space are continuous, with detailed proofs elaborated in \ref{app:1}.

\begin{lem}\label{lem:psiLip}
Let $\varphi:\mathbb{R}^{d-1} \to [0,\infty)$ be a Lipschitz function with Lipschitz constant $L$, $\Phi(\overline{\bx})=(\varphi(\overline{\bm{x}}),\overline{\bm{x}})$, and define the function $\psi:\mathbb{R}^{d-1} \to [-\delta,\infty)$ where $\psi(\overline{\bx}):=\min\{\tilde{x}\in \mathbb{R}: \text{dist}((\tilde{x},\overline{\bm{x}}),\Phi(\mathbb{R}^{d-1}))=\delta, \ \psi(\overline{\bx}) < \varphi(\overline{\bx})\}$ for all $\overline{\bm{x}} \in \mathbb{R}^{d-1}$.  Then $\psi$ is a Lipschitz function with the same Lipschitz constant $L$, which is independent of $\delta$.
\end{lem}

\begin{lem}\label{lem:phipsisupnorm}
  Let $\varphi$ and $\psi$ be as defined in Lemma \ref{lem:psiLip}.  Then  $\delta\leq |\varphi(\overline{\bx})-\psi(\overline{\bx})|\leq\delta\sqrt{L^2+1}$ for any $\overline{\bx}\in \R^{d-1}$.
\end{lem}

\begin{lem}\label{lem:projectiondist}
Let $\varphi$ and $\psi$ be defined as in the Lemma \ref{lem:psiLip}.  Then for $\bm{x}=(\tilde{x},\overline{\xb}),\bm{y}=(\tilde{y},\overline{\yb})\in \mcR_{-\delta}$, if
\begin{equation}
\label{eq:changeofvariable_x} 
\begin{aligned}
\tilde{x}'&=\left(1+\frac{\tilde{x}}{\delta}\right)\varphi(\overline{\bx})-\frac{\tilde{x}}{\delta}\psi(\overline{\bx})\\
 \tilde{y}'&=\left(1+\frac{\tilde{y}}{\delta}\right)\varphi(\overline{\by})-\frac{\tilde{y}}{\delta}\psi(\overline{\by})
\end{aligned}
\end{equation}
 and $\bx'=(\tilde{x}',\overline{\bx})$, $\by'=(\tilde{y}',\overline{\by})$, then
 \[
 K_L'|\bx-\by| \le|\bx'-\by'| \le K_L|\bx-\by|
 \]
 for some positive constants $K_L\geq 1$ and $K_L' \le 1$, {independent of $\delta$}. 
\end{lem}

\begin{lem}\label{lem:KernelEstimate}
{Let $\varphi$ and $\psi$ be defined as in the Lemma \ref{lem:psiLip}.} Let $M=\max(L+1+\sqrt{L^2+1}, K_L)$ with $K_L$ as in Lemma \ref{lem:projectiondist}, {$\xb=(\tx, \overline{\xb})$, $\yb=(\ty, \overline{\yb})$}, $\bm w = (\tilde{y}+\varphi(\overline{\bm{y}}),\overline{\bm{y}})$ and $\bz =(\tilde{x}+\varphi(\overline{\bm{x}}),\overline{\bm{x}})$.  Then we have the following kernel inequalities:
\begin{enumerate}[(a)]
\item For $\bx,\by \in \mcR^\infty$, $\gamma_{\delta/M}^{\,\beta}(|\bx-\by|) \le M^{d+p}\gamma_{\delta}^{\,\beta}(|\bm z - \bm w|){.}$
\item For $\bx \in \mcR_{-\delta}$ and $\by \in \mcR^\infty$, $\gamma_{\delta/M}^{\,\beta}(|\bx-\by|) \le M^{d+p}\gamma_{\delta}^{\,\beta}(|\bm x' - \bm w|){.}$
\item For $\bx,\by \in \mcR_{-\delta}$, $\gamma_{\delta/M}^{\,\beta}(|\bx-\by|) \le M^{d+p}\gamma_{\delta}^{\,\beta}(|\bm x' - \bm y'|){.}$
\end{enumerate}
where $\bx'$ and $\by'$ are defined as in \eqref{eq:changeofvariable_x}.
\end{lem}

\begin{lem}\label{lem:InverseKernelEstimate}
{Let $\varphi$ and $\psi$ be defined as in the Lemma \ref{lem:psiLip}.}
Let $M=\max(L+2, K'_L)$ with $K'_L$ as in Lemma \ref{lem:projectiondist}, $\bx'=(\tilde{x}',\overline{\bx})$,  $\by'=(\tilde{y}',\overline{\by})$, $\bm w = (\tilde{y}'-\varphi(\overline{\bm{y}}),\overline{\bm{y}})$ and $\bz =(\tilde{x}'-\varphi(\overline{\bm{x}}),\overline{\bm{x}})$.  Then we have the following kernel inequalities:
\begin{enumerate}[(a)]
\item For $\bx',\by' \in \Omega$, $\gamma_{\delta}^{\,\beta}(|\bx'-\by'|) \le M^{d+p}\gamma_{M\delta}^{\,\beta}(|\bm z - \bm w|)${.}
\item For $\bx' \in  \Omega_\delta$ and $\by '\in  \Omega$, $\gamma_{\delta}^{\,\beta}(|\bx'-\by'|) \le M^{d+p}\gamma_{
M\delta}^{\,\beta}(|\bm x - \bm w|)${.}
\item For $\bx',\by '\in  \Omega_\delta$, $\gamma_{M\delta}^{\,\beta}(|\bx'-\by'|) \le M^{d+p}\gamma_{M\delta}^{\,\beta}(|\bm x - \bm y|)${.}
\end{enumerate}
{where $\bx=(\tilde{x},\overline{\bx})$,  $\by=(\tilde{y},\overline{\by})$ such that $\tilde{x}$, $\tilde{y}$ and $\tilde{x}'$, $\tilde{y}'$ satisfy \eqref{eq:changeofvariable_x}. More specifically,
\begin{equation}
\label{eq:changeofvariable_x'} 
\begin{aligned}
\tilde{x}:=\frac{\delta(\tilde{x}'-\varphi(\overline{\xb}))}{\varphi(\overline{\xb})-\psi(\overline{\xb})},\quad \tilde{y}:=\frac{\delta(\tilde{y}'-\varphi(\overline{\yb}))}{\varphi(\overline{\yb})-\psi(\overline{\yb})}.
\end{aligned}
\end{equation}}
\end{lem}

\subsection{Extension to special Lipshitz domains}

We now first present a nonlocal trace theorem to generalize the earlier result shown on a half space with a flat boundary on one side to the case of any infinite domain whose boundary is defined by a Lipshitz graph.
\begin{thm}\label{thm:speclip}
Consider a special Lipschitz domain, $\Omega$, where there is a Lipschitz function $\varphi:\mathbb{R}^{d-1}\to \mathbb{R}$ with Lipschitz constant $L$ such that
\[
\Omega=\{\bm{x}=(\tilde{x},\overline{\bm{x}}) \in \mathbb{R}^d: \varphi(\overline{\bm{x}})<\tilde{x}, \ \overline{\bm{x}} \in \mathbb{R}^{d-1}\}
\]
with a $\delta$ collar boundary $\Omega_\delta$
\[
\Omega_\delta=\{\bm{x}=(\tilde{x},\overline{\bm{x}}) \in \mathbb{R}^d: \psi(\overline{\bm{x}})<\tilde{x}<\varphi(\overline{\bm{x}}), \ \overline{\bm{x}} \in \mathbb{R}^{d-1}\}.
\]
Here $\psi$ is as defined in Lemma \ref{lem:psiLip}.  
Then the nonlocal trace theorem holds on this domain, i.e.,
\[
\|u\|_{\cT_{\delta}^{\,\beta}(\Omega_\delta)}\le C|d+p-\beta|^{-1/p}\|u\|_{\cS_\delta^{\,\beta}(\hat\Omega)}
\]
for some constant $C$ which is independent of $\delta$,  $\beta$ and $u \in \cS_\delta^{\,\beta}(\hat\Omega)$.
\end{thm}
\begin{proof}
Define the operators $P_{\varphi,\psi}:L^p(\Omega) \to L^p(\mcR^\infty)$  and $G_{\varphi,\psi}:L^p(\Omega_\delta) \to L^p(\mcR_{-\delta})$ 
\begin{align}
\label{eq:projP}(P_{\varphi,\psi} u)(\bm{x})&=u(\tilde{x}+\varphi(\overline{\bm{x}}), \overline{\bm{x}}),\\
\label{eq:projG}(G_{\varphi,\psi} u)(\bm{x})&=u\left(\left(1+\frac{\tilde{x}}{\delta}\right)\varphi(\overline{\bm{x}})-\frac{\tilde{x}}{\delta}\psi(\overline{\bm{x}}),\overline{\bm{x}}\right).
\end{align} 
Then the operator $S_{\varphi,\psi}:L^p(\hat\Om)\to L^p(\mcR_{-\del}^\infty)$ is defined as
\begin{align}
\label{eq:projS}
S_{\varphi,\psi}u(\bm{x})=\begin{cases}
    P_{\varphi,\psi}u (\bm{x}) & \bm{x} \in \mcR^\infty,\\
    G_{\varphi,\psi}u(\bm{x}) & \bm{x} \in \mcR_{-\del}.
\end{cases}
\end{align}
We will show that $G_{\varphi,\psi}$ is a bounded operator from $\cT^{\,\beta}_\del(\Om_\del)$ to $\cT^{\,\beta}_\del(\mcR_{-\del})$ with a bounded inverse $G_{\varphi,\psi}^{-1}$. Moreover, $S_{\varphi,\psi}$ is a bounded operator from $\cS_{\delta}^{\,\beta}(\hat\Omega)$ to $\cS_{\delta}^{\,\beta}(\mcR^\infty_{-\del})$.

To show $G_{\varphi,\psi}$ is a bounded operator from $\cT^{\,\beta}_\del(\Om_\del)$ to $\cT^{\,\beta}_\del(\mcR_{-\del})$, we let $\bx'=(\tilde{x}',\overline{\xb})$,  $\by'=(\tilde{y}',\overline{\yb})$ and the constant $K_L\geq1$ as defined in Lemma \ref{lem:projectiondist}. From definition, we have $G_{\varphi,\psi} (u)(\xb)=u(\xb')$, $d \bx' = \frac{|\varphi(\overline{\bx})- \psi(\overline{\bx})|}{\del} d\bx$ and $d \by' = \frac{|\varphi(\overline{\by})- \psi(\overline{\by})|}{\del} d\by$.  Then,
\begin{equation}\label{eqn:GLpEstimate}
     \|G_{\varphi,\psi} u\|^p_{L^p(\mcR_{-\delta})}=\frac{1}{\delta}\int_{\mcR_{-\delta}}|G_{\varphi,\psi}(u)(\bx)|^p d \bx=\frac{1}{\delta}\int_{\mcR_{-\delta}}|u(\bx')|^p d \bx \le \frac{1}{\delta}\frac{\delta}{\inf |\varphi-\psi|}\int_{\Omega_{\delta}}|u(\bx')|^p d \bx'  \le \frac{1}{\delta}\|u\|^p_{L^p(\Omega_\delta)}
\end{equation}
and
\begin{align*}
    |G_{\varphi,\psi} u|^p_{\cT^{\,\beta}_{\delta}(\mcR_{-\delta})}&=\frac{1}{\delta}\int_{\mcR_{-\delta}}|G_{\varphi,\psi}(u)(\bx)|^p d \bx+\delta^{\,\beta-2}\int_{\mcR_{-\delta}}\int_{\mcR_{-\delta}}\frac{|G_{\varphi,\psi}(u)(\by)-G_{\varphi,\psi}(u)(\bx)|^p}{(|\by-\bx|\vee\delta)^{d+p-2}(|\by-\bx|\wedge\delta)^{\beta}}d\by d\bx\\
    &\le \delta^{\,\beta-2}\int_{\mcR_{-\delta}}\int_{\mcR_{-\delta}}\frac{|u(\by')-u(\bx')|^p}{(\frac{|\by'-\bx'|}{K_L}\vee\delta)^{d+p-2}(\frac{|\by'-\bx'|}{K_L}\wedge\delta)^{\beta}}d\by d\bx\\
    &\le \delta^{\,\beta-2}K_L^{d+\beta+p-2}\frac{\delta^p}{\inf |\varphi-\psi|^p}\int_{\Omega_\delta}\int_{\Omega_\delta}\frac{|u(\by')-u(\bx')|^p}{(|\by'-\bx'|\vee\delta)^{d+p-2}(|\by'-\bx'|\wedge\delta)^{\beta}}d\by' d\bx'\\
    &\le K_L^{2d+2p-2}\frac{\delta^p}{\inf |\varphi-\psi|^p}\delta^{\,\beta-2}\int_{\Omega_\delta}\int_{\Omega_\delta}\frac{|u(\by')-u(\bx')|^p}{(|\by'-\bx'|\vee\delta)^{d+p-2}(|\by'-\bx'|\wedge\delta)^{\beta}}d\by' d\bx'\\
    &\le K_L^{2d+2p-2}\verti{u}^p_{\cT^{\,\beta}_\delta(\Omega_{\delta})}
\end{align*}
where we have used Lemma \ref{lem:phipsisupnorm} for the last estimate as well as \eqref{eqn:GLpEstimate}. Together these estimates show \newline $\|G_{\varphi,\psi} u\|^p_{\cT^{\,\beta}_{\delta}(\mcR_{-\delta})}\le C\|G_{\varphi,\psi} u\|^p_{\cT^{\,\beta}_{\delta}(\Omega_\delta)}$ where the constant $C$ is independent of $\beta$ and $\delta$.

To show that $G_{\varphi, \psi}^{-1}$ exists and it is a bounded operator from  $\cT^{\,\beta}_\del(\mcR_{-\del})$ to  $\cT^{\,\beta}_\del(\Om_\del)$,  we note first that one can get \eqref{eq:changeofvariable_x'} from 
\eqref{eq:changeofvariable_x}. 
It is then easy to check that $G_{\varphi,\psi}^{-1}$ can be defined by $G_{\varphi,\psi}^{-1}u(\bx')= u(\bx)$ for $\bx'=(\tilde{x}',\overline{\bx}) \in \Om_\del$ and $\bx=(\tilde{x},\overline{\bx})\in \mcR_{-\del}$ where $\tilde{x}$ is given by \eqref{eq:changeofvariable_x'}. Using a change of variable estimate given by Lemma \ref{lem:projectiondist} we can then deduce the continuity of $G_{\varphi,\psi}^{-1}$ the same way as we have done for $G_{\varphi,\psi}$ with details omitted. 


Now to show that $S_{\varphi,\psi}$ is a bounded operator from $\cS^{\,\beta}_\del(\hat\Om)$ to $\cS^{\,\beta}_\del(\mcR_{-\del}^\infty)$, we note that 
$\|S_{\varphi,\psi}u \|_{L^p(\mcR_{-\del}^\infty)} = \|u \|_{L^p(\hat\Om)}$, 
and by applying Lemma \ref{lem:EnergyKernelEst} with $U=\mcR_{-\del}^\infty$ and $\alpha= 1/M$ where $M>1$ is defined in Lemma \ref{lem:KernelEstimate}, we have
\begin{align*}
|S_{\varphi,\psi}(u)|^p_{\cS^{\,\beta}_{\delta}(\mcR_{-\delta}^\infty)} &\le C \int_{\mcR_{-\delta}^\infty}\int_{\mcR_{-\delta}^\infty}\gamma^{\,\beta}_{\delta/M}(|\by-\bx|)|S_{\varphi,\psi}(u)(\by)-S_{\varphi,\psi}(u)(\bx)|^p d\by d\bx\\
&\le C \underbrace{\int_{\mcR^\infty}\int_{\mcR^\infty}\gamma^{\,\beta}_{\delta/M}(|\by-\bx|)|u(\tilde{y}+\varphi(\overline{\by}), \overline{\by})-u(\tilde{x}+\varphi(\overline{\bx}), \overline{\bx})|^pd\by d\bx}_{I_1}\\
&+C\underbrace{\int_{\mcR_{-\delta}}\int_{\mcR^\infty}\gamma^{\,\beta}_{\delta/M}(|\by-\bx|)\left|u(\tilde{y}+\varphi(\overline{\by}), \overline{\by})-u\left(\left(1+\frac{\tilde{x}}{\delta}\right)\varphi(\overline{\bm{x}})-\frac{\tilde{x}}{\delta}\psi(\overline{\bm{x}}),\overline{\bm{x}}\right)\right|^pd\by d\bx}_{I_2}\\
&+C\underbrace{\int_{\mcR_{-\delta}}\int_{\mcR_{-\delta}}\gamma^{\,\beta}_{\delta/M}(|\by-\bx|)\left|u\left(\left(1+\frac{\tilde{y}}{\delta}\right)\varphi(\overline{\bm{y}})-\frac{\tilde{y}}{\delta}\psi(\overline{\bm{y}}),\overline{\bm{y}}\right)-u\left(\left(1+\frac{\tilde{x}}{\delta}\right)\varphi(\overline{\bm{x}})-\frac{\tilde{x}}{\delta}\psi(\overline{\bm{x}}),\overline{\bm{x}}\right)\right|^pd\by d\bx}_{I_3}.
\end{align*}
{Here the constant $C$ is independent of $\beta$ and $\delta$.}

Applying Lemma \ref{lem:KernelEstimate} part (a) to  $I_1$ and then using a change of variable, we have
\[
I_1\le M^{d+p} \int_{\Om}\int_{\Om} \gamma^{\,\beta}_\delta(|\bm{w}-\bm{z}|)|u(\bm{w})-u(\bm{z})|^pd\bm{w}d\bm{z} \le C |u|^p_{\cS^{\,\beta}_\del(\hat\Om)}.
\]
Similarly applying Lemma \ref{lem:KernelEstimate} to $I_2$ and $I_3$ along with the proper change of variables, we finally have  $|S_{\varphi,\psi}(u)|^p_{\cS^{\,\beta}_{\delta}(\mcR_{-\delta}^\infty)} \le C|u|^p_{\cS^{\,\beta}_{\delta}(\hat\Om)}$ where $C$ is independent of $\delta$ and $\beta$.  Thus, we have the continuity of $S_{\varphi,\psi}$.  Using proven properties of $G_{\varphi,\psi}$ and $S_{\varphi,\psi}$ along with the trace theorem for the half-plane, we have
 \begin{align*}
     \|u\|_{\cT_{\delta}^\beta(\Om_{\delta})} &=\|G^{-1}_{\varphi,\psi}(G_{\varphi,\psi}(u))\|_{\cT^{\,\beta}_{\delta}(\Om_{\delta})}\le C\|G_{\varphi,\psi}(u)\|_{\cT^{\,\beta}_{\delta}(\mcR_{-\delta})}= C\|S_{\varphi,\psi}(u)\|_{\cT^{\,\beta}_{\delta}(\mcR_{-\delta})}\\
     &\le C|d+p-\beta|^{-1/p}\|S_{\varphi,\psi}(u)\|_{\cS^{\,\beta}_{\delta}(\mcR_{-\delta}^\infty)}\le C|d+p-\beta|^{-1/p}\|u\|_{\cS^{\,\beta}_{\delta}(\hat\Om)}.
 \end{align*}
\end{proof} 

Using the transformation operators $G_{\varphi,\psi}$ and $S_{\varphi,\psi}$, we can also generalize the inverse nonlocal trace theorem to the special Lipschitz domain.  
\begin{thm}\label{thm:invSpecLip}
Let $\phi$, $\psi$, $\Omega$ and $\Omega_\delta$ be as defined in Theorem \ref{thm:speclip} and the extension operator $E:\cT^{\,\beta}_\delta(\mcR_{-\delta}) \to \cS^\beta_\delta(\mcR_{-\delta}^\infty)$ as defined in Theorem  \ref{mainthm_2}.  Let $G_{\varphi,\psi}$ be defined in \eqref{eq:projG} and $S_{\varphi,\psi}$ in \eqref{eq:projS}. 
Then we can define an extension operator $\tilde{E}:\cT_\delta^\beta(\Omega_\delta) \to \cS^\beta_\delta(\hat{\Omega})$ where $\tilde{E}=S^{-1}_{\varphi,\psi} E G_{\varphi,\psi}$ and
\[
\|\tilde{E} u\|_{\cS_\delta^\beta(\hat\Omega)} \le C |d-\beta|^{-1/p} \|u\|_{\cT_\delta^\beta(\Omega_\delta)}.
\]
Here $C$ is a constant independent of $\delta$, $\beta$ and $u\in \cT_\delta^\beta(\Omega_\delta)$.
\end{thm}
\begin{proof}
Notice that the inverse operator  $S^{-1}_{\varphi,\psi}: L^p(\mcR_{-\del}^\infty)\to L^p(\hat\Om) $ can be defined as 
\[
S^{-1}_{\varphi,\psi} u (\bx') =
\left\{
\begin{aligned}&
u\left(\frac{\del(\tilde{x}'-\varphi(\overline{\bx}))}{\varphi(\overline{\bx})-\psi(\overline{\bx})}, \overline{\bx}\right), & \quad  \bx'=(\tilde{x}',\overline{\bx}) \in \Omega_\delta,\\
  &  u(\tilde{x}'-\varphi(\overline{\bx}),\overline{\bx}), & \quad \bx' = (\tilde{x}',\overline{\bx}) \in \Omega.
\end{aligned}
\right.   
\]
Now we wish to show that $S^{-1}_{\varphi,\psi}$ is a bounded operator from $\cS^{\,\beta}_\del(\mcR_{-\del}^\infty)$ to  $\cS^{\,\beta}_\del(\hat\Om)$.
Using Lemma \ref{lem:InverseKernelEstimate} and the $M$ defined there along with the appropriate changes of the variables, we can show similarly as we have done in Theorem \ref{thm:speclip}  that 
\begin{align*}
|S^{-1}_{\varphi,\psi}(u)|^p_{\cS^{\,\beta}_{\delta}(\hat\Om)} 
\leq C  \int_{\mcR^\infty_{-\del}}\int_{\mcR^\infty_{-\del}}\gamma^{\,\beta}_{M\delta}(|\by-\bx|)|u(\by)-u(\bx)|^p d\by d\bx\leq C \int_{\mcR^\infty_{-\del}}\int_{\mcR^\infty_{-\del}}\gamma^{\,\beta}_{\delta}(|\by-\bx|)|u(\by)-u(\bx)|^p d\by d\bx ,
\end{align*}
where we have used Lemma \ref{lem:EnergyKernelEst} with $U=\mcR_{-\del}^\infty$ and $\alpha=1/M$ in the last step. The continuity of $S^{-1}_{\varphi,\psi}$ is thus shown. 

Finally, using the continuity properties of $S_{\varphi,\psi}^{-1}$ , $G_{\varphi,\psi}$ and $E$  we have
\[
\begin{split}
&\|\tilde{E} u\|_{\cS^{\,\beta}_\delta(\hat\Omega)}=\|S^{-1}_{\varphi,\psi} E G_{\varphi,\psi} u\|_{\cS^{\,\beta}_\delta(\hat\Omega)} \le C\|E G_{\varphi,\psi} u\|_{\cS^{\,\beta}_{\delta}(\mcR_{-\delta}^\infty)} \\
\le &C |d-\beta|^{-1/p}\| G_{\varphi,\psi} u\|_{\cT^{\,\beta}_{\delta}(\mcR_{-\delta})} \le C |d-\beta|^{-1/p}\| u\|_{\cT^{\,\beta}_{\delta}(\Omega_\delta)}.
\end{split}
\]

\end{proof}

\subsection{Extension to more general Lipshitz domains}

The extension to more general Lipshitz domains can be otained by using the partition of unity technique.
We first decompose the boundary collar region into finitely many balls so that we can locally view the boundary as multiple special Lipschitz domains.  From there we apply Theorem \ref{thm:speclip} to each part and join the estimates together with a partition of unity. The detailed derivation is given as follows.

First, for all discussions given in this subsection, we state some assumptions on the domains and define the necessary spaces and functions.

We consider a general bounded simply connected Lipschitz domain $\Omega$, which naturally makes $\overline{\Omega_{\delta}}$ a compact set for any finite $\delta>0$. Since $\Omega$ has a Lipschitz boundary (see, e.g., \cite{grisvard1985elliptic} Def 1.2.1.1), there exist $N$ local coordinate systems $\bx^i=(x_1^i,x_2^i,\cdots,x_d^i)$ for $1\le i\le N$, a collection of balls $\{B(\bx_i,r_i)\}_{i=1}^N$, and Lipschitz functions $\varphi_i: \real^{d-1} \to \real$ for some $N \in \mathbb{N}$, $\bx_i \in \partial \Omega, \ r_i>0$ such that $\partial \Omega \subseteq\bigcup_{i=1}^N B(\bx_i,r_i)$, 
\begin{align*}
\{\bx^i=(\tilde{x}^i,\overline{\bx}^i) \in B(\bx_i,r_i)  :  \varphi_i(\overline{\bx}^i) < \tilde{x}^i\} = \Omega \cap B(\bx_i,r_i),\\
\{\bx^i=(\tilde{x}^i,\overline{\bx}^i) \in B(\bx_i,r_i)  :  \varphi_i(\overline{\bx}^i) = \tilde{x}^i\} = \partial\Omega \cap B(\bx_i,r_i).
\end{align*}
Notice that in the above definition, $\{ \bx_i\}_{i=1}^N$ are $N$ fixed points on $\partial\Om$. 
Letting $\delta_0:=\frac{1}{2}\dist(\real^d \setminus\bigcup_{i=1}^N B(\bx_i,r_i),\partial\Omega)$, we have $\Omega_{\delta_0}\subseteq\bigcup_{i=1}^N B(\bx_i,r_i)$ and also a positive number $\epsilon<\delta_0$ such that $\Omega_{\delta_0}\subseteq\bigcup_{i=1}^N B(\bx_i,r_i-2\epsilon)$. Defining $\psi_i^0(\overline{\bx}^{i})$ from $\varphi_i$ as in Lemma \ref{lem:psiLip} we have
\begin{align*}
(\Omega \cup \Omega_{\delta_0})\cap B(\bx_i,r_i)&:=\Omega^i_0 \cap B(\bx_i,r_i)\\
\Omega_{\delta_0}\cap B(\bx_i,r_i)& :=\Omega_{\delta_0}^i \cap B(\bx_i,r_i)
\end{align*}
where $\Omega^i_0 =\{\bx^i \in \mathbb{R}^d:\psi_i^0(\overline{\bx}^{i})<\tilde{x}^{i}\}$ and $\Omega^i_{\delta_0} =\{\bx^i  \in \mathbb{R}^d:\psi_i^0(\overline{\bx}^{i})<\tilde{x}^{i}<\varphi_i(\overline{\bx}^{i})\}$.
Notice here the Lipschitz constants of $\varphi_i$ and $\psi_i^0$ 
depend on the domain $\Omega$, $\delta_0$ and the collection of balls
$\{B(\bx_i,r_i)\}$ only and 
are thus independent of $\delta$ and $\beta$. 
Then given $\delta\in(0,\epsilon)$ 
we have $\Omega_\delta \subset \Omega_{\delta_0} \subset\bigcup_{i=1}^N B(\bx_i,r_i-2\epsilon)$ along with functions $\psi_{i}(\overline{\bx}^{i})$ so that $\hat{\Omega}\cap B(\bx_i,r_i)=\Omega^i \cap B(\bx_i,r_i)$ and $
\Omega_{\delta}\cap B(\bx_i,r_i)=\Omega_{\delta}^i \cap B(\bx_i,r_i)$ where
\[
\Omega^i =\{{\bx^i}  \in \mathbb{R}^d:\psi_i(\overline{\bx}^{i})<\tilde{x}^{i}\} \ \text{and} \  \Omega^i_{\delta} =\{{\bx^i}  \in \mathbb{R}^d:\psi_i(\overline{\bx}^{i})<\tilde{x}^{i}<\varphi_i(\overline{\bx}^{i})\}.
\]
Define functions $\{ \lambda_i\}_{i=1}^N$ such that
\begin{enumerate}
    \item $\lambda_i \in C_c^\infty(B(\bx_i,r_i-\epsilon)), 1 \le i \le N$,
    \item $0 \le \lambda_i \le 1$ and $ \lambda_i\equiv 1$ on  $B(\bx_i,r_i-2\epsilon)$.
\end{enumerate}
Since $\Om_\del$ is covered by $\{B(\bx_i,r_i-2\epsilon)\}_{i=1}^N$, we have $1 \leq \sum_{i=1}^N \lambda_i(\bx) \leq C$ for a fixed constant $C>0$ depending only on the maximum number of overlapped balls in the set $\{ B(\xb_i,r_i-\epsilon)\}_{i=1}^N$. We also define \[
\widetilde{\lambda_i}(\bx) = \lambda_i(\bx)/ \sum_{j=1}^N \lambda_j^2(\bx)
\]
for each $i\in \{1,\cdots, N \}$. Then since $\sum_{j=1}^N \lambda_j^2(\bx)$ is uniformly bounded above and below, we also have $\widetilde{\lambda_i} \in C_c^\infty(B(\bx_i,r_i-\epsilon))$.

Furthermore, we would like to define the extension operator on the general domain.  By Theorem \ref{thm:invSpecLip}, there exists an extension operator $E^i: \cT^{\,\beta}_\del(\Om^i_\del) \to \cS^{\,\beta}_\del(\Om^i) $ such that $E^i(\lambda_i u)(\bx) =(\lambda_i u)(\bx) $ for $x\in \Om^i_\del$ and
 \[
 \| E^i(\lambda_i u)\|_{\cS^{\,\beta}_\del(\Om^i)} \leq C|d-\beta|^{-1/p} \| \lambda_i u\|_{\cT^{\,\beta}_\del(\Om_\del^i)}\,.
 \]
For any $\xb \in \hat\Om \setminus\Om^i$, we assume $E^i(\lambda_i u)(\xb) =0$. Then we can define the extension operator $E: L^p(\Om_\del) \to L^p(\hat\Om)$ by 
 \begin{equation}
 \label{def:extension_generalLip}
 Eu(\bx) = \sum_{i=1}^N \widetilde{\lambda_i}(\bx) E^i(\lambda_i u)(\bx) \,, 
 \end{equation}
for any $\bx\in \hat\Om$.  Before continuing we show some useful estimates which will together swiftly prove Theorem \ref{mainthm_1_general}.
\begin{lem}\label{lem:GenTraceESt1}
Let $\Omega$ be a simply connected Lipschitz domain with an interaction domain $\Omega_\delta$ for $0<\delta< \epsilon$ where $\epsilon$, $\Omega^i$, and $\Omega_\delta^i$,
as well as $\lambda_i$, are defined as above. For any $u \in \cS_\delta^{\,\beta}(\hat\Om)$,  $\| u\|^p_{\cT_\delta(\Omega_{\delta})} \le C \sum_{i=1}^N\|\lambda_i u\|^p_{\cT_\delta(\Omega_{\delta}^i)}$,
where $C$ is a constant independent of $\delta$ and $\beta$.
\end{lem}
\begin{proof}
First note for $\bx \in \Omega_{\delta}$,
\[
|u(\bx)|\le\sum_{i=1}^N |(\lambda_i u)(\bx)|
\]
and hence
\[
\|u\|_{\cT^{\,\beta}_\delta(\Omega_{\delta})}\le \sum_{i=1}^N\|\lambda_i u\|_{\cT^{\,\beta}_\delta(\Omega_{\delta})}.
\]
Notice that we can extend $\lambda_i$ by zero so that $\lambda_i u$ can be viewed as a function on $\Om^i$.  Using the fact that $\lambda_i\in C_c^\infty(B(\bx_i,r_i))$ we have
\begin{align}
    |\lambda_i u|^p_{\cT^{\,\beta}_\delta(\Omega_{\delta})}&=\delta^{\,\beta-2}\left(\int_{\Omega_\delta\cap B(\bx_i,r_i)}\int_{\Omega_\delta\cap B(\bx_i,r_i)}+2\int_{\Omega_\delta\cap B(\bx_i,r_i)}\int_{\Omega_\delta\setminus B(\bx_i,r_i)}\right)\frac{|\lambda_i u(\by)-\lambda_i u(\bx)|^p}{(|\by-\bx|\vee\delta)^{d+p-2}(|\by-\bx|\wedge \delta)^{\beta}} d \by d \bx \nonumber\\
    =&\delta^{\,\beta-2}\left(\int_{\Omega_\delta^i\cap B(\bx_i,r_i)}\int_{\Omega_\delta^i\cap B(\bx_i,r_i)}+2\int_{\Omega_\delta^i\cap B(\bx_i,r_i)}\int_{\Omega_\delta\setminus B(\bx_i,r_i)}\right)\frac{|\lambda_i u(\by)-\lambda_i u(\bx)|^p}{(|\by-\bx|\vee\delta)^{d+p-2}(|\by-\bx|\wedge \delta)^{\beta}} d \by d \bx \nonumber \\
    \le &|\lambda_i u|^p_{\cT_\delta(\Omega_{\delta}^i)}+2^p\delta^{\,\beta-2}\int_{\Omega_\delta^i\cap B(\bx_i,r_i-\epsilon)}\int_{\Omega_\delta\setminus B(\bx_i,r_i)}\frac{|\lambda_i u(\bx)|^p}{(|\by-\bx|\vee\delta)^{d+p-2}(|\by-\bx|\wedge \delta)^{\beta}} d \by d \bx. \label{eqn:GeneralLipSemiBound}
\end{align}
Notice that for each $\bx \in \Omega_\delta\cap B(\bx_i,r_i-\epsilon)$, and $\delta \in (0,\epsilon)$ we have 
\begin{align}
    &\delta^{\,\beta-2}\int_{\Omega_\delta\setminus B(\bx_i,r_i)}\frac{1}{(|\by-\bx|\vee\delta)^{d+p-2}(|\by-\bx|\wedge \delta)^{\beta}} d \by=\delta^{-2}\int_{\Omega_\delta\setminus B(\bx_i,r_i)}\frac{1}{|\by-\bx|^{d+p-2}} d \by \nonumber \\
   &\le \delta^{-2}\int_{\Omega_\delta \cap\{\by \in \mathbb{R}^d: |\by-\bx|>\epsilon\}}\frac{1}{|\by-\bx|^{d+p-2}} d \by\le\delta^{-2}\int_{\tilde{x}-\delta/2}^{\tilde{x}+\delta/2}\int_{|\overline{\by}-\overline{\bx}|>\epsilon-\delta/2}\frac{1}{(|\tilde{y}-\tilde{x}|+|\overline{\by}-\overline{\bx}|)^{d+p-2}} d \overline{\by} d \tilde{y} \label{eqn:GeneralLipSemiBound2}\\
   &\le \delta^{-1} \int_{|\overline{\by}|>{\epsilon/2}}\frac{1}{|\overline{\by}|^{d+p-2}} d \overline{\by}\le C\delta^{-1} \nonumber
\end{align}
where $C$ depends on $\epsilon$, $d$ and $p$.  Therefore, the integral term in \eqref{eqn:GeneralLipSemiBound} is bounded by a multiple of $\frac{1}{\delta}\|\lambda_iu\|^p_{L^p(\Omega_\delta^i)}$. 
Moreover, using the compact support of $\lambda_i$,
\begin{equation}\label{eqn:GeneralLipL2Bound}
\|\lambda_iu\|^p_{L^p(\Omega_\delta)}=\|\lambda_iu\|^p_{L^p(\Omega_\delta\cap B(\bx_i,r_i))}=\|\lambda_iu\|^p_{L^p(\Omega_\delta^i)}.
\end{equation}
Therefore \eqref{eqn:GeneralLipL2Bound} along with the estimates of \eqref{eqn:GeneralLipSemiBound} and \eqref{eqn:GeneralLipSemiBound2}
\[
 \|\lambda_i u\|^p_{\cT_\delta(\Omega_{\delta})} \le C \|\lambda_i u\|^p_{\cT_\delta(\Omega_{\delta}^i)}
\]
where the constant $C$ is independent of $\beta$ and $\delta$.
\end{proof}

\begin{lem}\label{lem:GenTraceEst2}
Let $\Omega$ be a simply connected Lipschitz domain with an interaction domain $\Omega_\delta$ for $0<\delta< \epsilon$ where $\epsilon$, $\Omega^i$, and $\Omega_\delta^i$ are defined as above, along with $\lambda_i$ and $\widetilde{\lambda_i}$ for each $i$. 
For any $u \in \cS_\delta^{\,\beta}(\hat\Om)$, we have
\[
\|\lambda_i u\|^p_{\cS_\delta^{\,\beta}(\Omega^i)} \le C \| u\|^p_{\cS_\delta^{\,\beta}(\hat\Omega \cap B(\xb_i ,r_i))}, \quad\text{ and }\quad \|\widetilde{\lambda_i} u\|^p_{\cS_\delta^{\,\beta}(\Omega^i)} \le C \| u\|^p_{\cS_\delta^{\,\beta}(\hat\Omega \cap B(\xb_i ,r_i))},
\]
where $C$ is independent of $\beta$ and $\delta$.
\end{lem}
\begin{proof}
First, note that $\|\lambda_i u\|_{
L^p(\Omega^i)} =\|\lambda_i u\|_{L^p(\hat{\Omega}\cap B(\bx_i,r_i))}$ and also if $\by \in \Omega^i\setminus B(\bx_i,r_i)$ and $\bx \in \Omega^i\cap B(\bx_i,r_i-\epsilon)$, then 
\[
|\by-\bx|\ge |\by-\bx_i|-|\bx-\bx_i|> r_i-(r_i-\epsilon)=\epsilon > \delta.
\]
Moreover, using the compact support of $\lambda_i$,
\begin{align*}
    |\lambda_i u|^p_{\cS_\delta^{\,\beta}(\Omega^i)} &= \left(\int_{\Omega^i\cap B(\bx_i,r_i)}\int_{\Omega^i\cap B(\bx_i,r_i)}+2\int_{\Omega^i\cap B(\bx_i,r_i)}\int_{\Omega^i\setminus B(\bx_i,r_i)}\right)\gamma^{\,\beta}_\delta(|\bx-\by|)|\lambda_i u(\bx)-\lambda_i u(\by)|^pd\by d\bx\\
     &\le |\lambda_i u|^p_{\cS_\delta^{\,\beta}(\Omega^i\cap B(\bx_i,r_i))}+2\int_{\Omega^i\cap B(\bx_i,r_i-\epsilon)}\int_{\Omega^i\cap\{\by \in\mathbb{R}^d: |\by-\bx|>\epsilon\}}\gamma^{\,\beta}_\delta(|\bx-\by|)|\lambda_iu(\bx)|^pd\by d \bx=|\lambda_i u|^p_{\cS_\delta^{\,\beta}(\hat{\Omega}\cap B(\bx_i,r_i))}
\end{align*}
where the last equality is because $\gamma_\delta^\beta(|\bx-\by|)=0$ since $|\bx-\by|>\epsilon>\delta$.
Then for each $1 \le i \le N$,
\begin{align}\label{eqn:lambdaseminormest1}
&|\lambda_i u|_{\cS_\delta^{\,\beta}(\hat{\Omega}\cap B(\bx_i,r_i))}^p= \int_{\hat{\Omega} \cap B(\bx_i,r_i)}\int_{\hat{\Omega} \cap B(\bx_i,r_i)}\gamma^{\,\beta}_\delta(|\bx-\by|)|\lambda_i(\by)u(\by)-\lambda_i(\bx)u(\bx)|^pd\by d\bx\nonumber\\
 &\le  2^{p-1}\int_{\hat{\Omega} \cap B(\bx_i,r_i)}\int_{\hat{\Omega} \cap B(\bx_i,r_i)}\gamma^{\,\beta}_\delta(|\bx-\by|)|\lambda_i(\by)|^p |u(\by)-u(\bx)|^p+\gamma^{\,\beta}_\delta(|\bx-\by|)|u(\bx)|^p|\lambda_i(\by)-\lambda_i(\bx)|^p d\by d\bx\nonumber\\
&\le 2^{p-1}|u|^p_{\cS_\delta^{\,\beta}(\hat{\Omega}\cap B(\bx_i,r_i))} +2^{p-1}\int_{\hat{\Omega} \cap B(\bx_i,r_i)}|u(\bx)|^p\int_{\hat{\Omega} \cap B(\bx_i,r_i)}\gamma^{\,\beta}_{\delta}(|\bx-\by|)|\lambda_i(\by)-\lambda_i(\bx)|^pd\by d\bx\nonumber\\
&\le 2^{p-1}|u|^p_{\cS_\delta^{\,\beta}(\hat{\Omega}\cap B(\bx_i,r_i))} +2^{p-1} \| \lambda_i\|^p_{C^1}\int_{\hat{\Omega} \cap B(\bx_i,r_i)}|u(\bx)|^p\int_{\hat{\Omega} \cap B(\bx_i,r_i)}\gamma^{\,\beta}_{\delta}(|\bx-\by|)|\by-\bx|^pd\by d\bx\nonumber\\
& \leq  C \| u\|^p_{\cS_\delta^{\,\beta}(\hat{\Omega}\cap B(\bx_i,r_i))} ,
\end{align}
where $C$ is independent of $\beta$ and $\delta$.
Also, since $|\lambda_i| \le 1$,
\[
\|\lambda_i u\|_{L^p(\hat\Om\cap B(\bx_i,r_i))}^p \le \|u\|_{L^p(\hat\Om\cap B(\bx_i,r_i))}^p \le \|u\|_{L^p(\hat\Om)}^p.
\]
The estimates for $\widetilde{\lambda_i}$ can be done similarly.
\end{proof}

\begin{lem}\label{lem:GenInvTraceEst2}
 Let $\Omega$ be a simply connected Lipschitz domain with an interaction domain $\Omega_\delta$ for $0<\delta< \epsilon$ where $\epsilon$, $\Omega^i$, and $\Omega_\delta^i$,
 along with 
   $\lambda_i$,
  are defined as above. For any $u \in \cT_\delta^{\,\beta}(\Om_\delta)$, we have
  \[
  \| \lambda_i u (\bx) \|^p_{\cT^\beta_\del(\Om^i_\del)} \le C \left( \| u \|^p_{\cT_\del^{\,\beta}(\Om_\del)} + |d+p-\beta|^{-1}\| u\|^p_{L^p(\Om_\del)}\right),
  \]
  where $C$ is independent of $\delta$ and $\beta$. 
\end{lem}
\begin{proof} Since $\lambda_i$ is supported in $B(\bx_i, r_i-\epsilon)$ and $\lambda_i\leq 1$, it is obvious that $\| \lambda_i u (\bx) \|_{L^p(\Om^i_\del)} =\| \lambda_i u (\bx) \|_{L^p(\Om_\del)}\leq \|u (\bx) \|_{L^p(\Om_\del)}$. Now 
\begin{equation}
\label{eq1:GenInvTraceEst2}
\begin{split}
|\lambda_i u|^p_{\cT^{\,\beta}_\del(\Om^i_\del)} &= \del^{\,\beta-2}\left(\int_{\Om^i_\del \cap B(\bx_i, r_i)}\int_{\Om^i_\del \cap B(\bx_i, r_i)} + 2 \int_{\Om^i_\del \cap B(\bx_i, r_i)} \int_{\Om_\del^i \backslash B(\bx_i, r_i)}  \right) \frac{|\lambda_i u (\by)-\lambda_i u (\bx) |^p}{(|\by-\bx|\vee\del)^{d+p-2}(|\by-\bx|\wedge\del)^{\,\beta}} d\by d\bx   \\
&= |\lambda_i u|^p_{\cT^{\,\beta}_\del(\Om^i_\del\cap B(\xb_i, r_i))} + 2 \int_{\Om^i_\del \cap B(\bx_i, r_i-\epsilon)} \int_{\Om_\del^i \backslash B(\bx_i, r_i)}   \frac{|\lambda_i u (\bx) |^p}{(|\by-\bx|\vee\del)^{d+p-2}(|\by-\bx|\wedge\del)^{\,\beta}} d\by d\bx 
\end{split}
\end{equation}
For the first term in the last line, we have
\[
\begin{split}
|\lambda_i u|^p_{\cT^{\,\beta}_\del(\Om^i_\del\cap B(\xb_i, r_i))} &\leq 2^{p-1} \del^{\,\beta-2}\int_{\Om^i_\del \cap B(\bx_i, r_i)}\int_{\Om^i_\del \cap B(\bx_i, r_i)}  \frac{|\lambda_i(\by)|^p |u(\by)-u(\bx)|^p+|u(\bx)|^p|\lambda_i(\by)-\lambda_i(\bx)|^p}{(|\by-\bx|\vee\del)^{d+p-2}(|\by-\bx|\wedge\del)^{\,\beta}} d\by d\bx \\
 \leq& 2^{p-1} | u|^p_{\cT^{\,\beta}_\del(\Om_\del)} + 2^{p-1} \del^{\,\beta-2}\| \lambda_i\|_{C^1}^p  \int_{\Om_\del \cap B(\bx_i, r_i)}|u(\bx)|^p \int_{\Om^i_\del \cap B(\bx_i, r_i)}  \frac{|\by-\bx|^p}{(|\by-\bx|\vee\del)^{d+p-2}(|\by-\bx|\wedge\del)^{\,\beta}} d\by d\bx. 
\end{split}
\]
Now for any $\xb\in \Om^i_\del \cap B(\bx_i, r_i)$, 
\[
\begin{split}
&\del^{\,\beta-2}
 \int_{\Om^i_\del \cap B(\bx_i, r_i)}  \frac{|\by-\bx|^p}{(|\by-\bx|\vee\del)^{d+p-2}(|\by-\bx|\wedge\del)^{\,\beta}} d\by \\
 \leq& \del^{\,\beta-d-p}  \int_{|\yb-\xb|<\del}\frac{1}{|\yb-\xb|^{\,\beta-p}}  d\by + \del^{-2}  \int_{\{ \yb\in \Om_\del^i: \del<|\yb-\xb|<2r_i\}}\frac{1}{|\yb-\xb|^{d-2}}  d\by \\
  \leq& C |d+p-\beta|^{-1} + \del^{-2}  \int_{ \left\{ \by\in\R^d:  \psi_i(\overline{\by})<\tilde{y}<\varphi_i(\overline{\by}),\, |\by-\bx|<2r_i\right\}}\frac{1}{|\yb - \xb|^{d-2}}  d\by \,,
 \end{split}
\]
where $C$ only depends on $d$. Note that in the last step above, we have adopted the $i$-th local coordinate system defined before to represent the points $\bx$ and $\by$ without adding labels to them. Define $G^i: \mcR_{-\del}\to \{ \by\in\R^{d}:  \psi_i(\overline{\by})<\tilde{y}<\varphi_i(\overline{\by}) \}$ such that
for any $\bx=(\tilde{x}, \overline{\bx})$, with the same $i$-th local coordinate system representation
\begin{equation}
\label{eq2:GenInvTraceEst2}
G^i \bx = \left(\left(1+\frac{\tilde{x}}{\delta}\right)\varphi_i(\overline{\bx})-\frac{\tilde{x}}{\delta}\psi_i(\overline{\bx}), \overline{\bx}\right)\,.
\end{equation}
We can see from Lemma \ref{lem:projectiondist} that $K_L^\prime |\bx-\bz|\leq | G^i \bx - G^i \bz| \leq K_L |\bx -\bz| $, so $\left|\frac{\partial G^i}{\partial \bx}\right|\leq K_L$. 
By the change of variable,
\[
\begin{split}
&\del^{-2}  \int_{ \left\{ \by\in\R^d:  \psi_i(\overline{\by})<\tilde{y}<\varphi_i(\overline{\by}),\, |\by-\bx|<2r_i\right\}}\frac{1}{|\yb -\xb|^{d-2}}  d\by
\leq K_L \del^{-2} \int_{\left\{\yb\in\mcR_{-\del}: |G^i \yb-\xb|<2r_i\right\}} \frac{1}{|G^i \yb -\xb|^{d-2}} d\yb \\
\leq & \frac{K_L}{(K_L^\prime)^{d-2}} \del^{-2} \int_{\left\{\yb\in\mcR_{-\del}: |\yb - (G^i)^{-1}\xb|<2r_i/K_L^\prime\right\}} \frac{1}{| \yb - (G^i)^{-1}\xb|^{d-2}} d\yb
\end{split}. 
\]
Let $\bm{w} = (G^i)^{-1}\bx = (\tilde{w}, \overline{\bm{w}})$ and $
\mcR_{-\del} - \bm{w} := \{ \yb-\bm{w}: \yb\in \mcR_{-\del}\} = (-\del-\tilde{w}, -\tilde{w} )\times \R^{d-1}
$, the last line can be estimated by
\[
\begin{split}
&\frac{K_L}{(K_L^\prime)^{d-2}} \del^{-2} \int_{\left\{\yb\in\mcR_{-\del} -\bm{w}: |\yb|<2r_i/K_L^\prime\right\}} \frac{1}{|\yb|^{d-2}} d\yb \leq \frac{K_L}{(K_L^\prime)^{d-2}} \del^{-2} \int_{\left\{\yb\in\mcR_{-\del} -\bm{w}: |\overline{\yb}|<2r_i/K_L^\prime\right\}} \frac{1}{|\overline{\yb}|^{d-2}} d\yb \\
&= \frac{K_L}{(K_L^\prime)^{d-2}} \del^{-2} \int_{-\del-\tilde{w}}^{-\tilde{w}} d\tilde{y} \int_{ |\overline{\yb}|<2r_i/K_L^\prime} \frac{1}{|\overline{\yb}|^{d-2}} d\overline{\yb}  \leq C \del^{-1}\,,
\end{split}
\]
where $C$ is independent of $\del$ and $\beta$. Collecting the above estimates, we have
\[
|\lambda_i u|^p_{\cT^{\,\beta}_\del(\Om^i_\del\cap B(\xb_i, r_i))}\leq C \left( \| u \|^p_{\cT_\del^{\,\beta}(\Om_\del)} + |d+p-\beta|^{-1}\| u\|^p_{L^p(\Om_\del)}\right). 
\]

Now for the second term in \eqref{eq1:GenInvTraceEst2}, we notice that if $\bx\in B(\bx_i, r_i-\epsilon)$ and $\by\in \R^d \setminus B(\bx_i, r_i)$, then $|\by-\bx|\geq \epsilon>\del$, therefore for each $\bx\in B(\bx_i, r_i-\epsilon)$ 
\[
\begin{split}
&\del^{\,\beta-2} \int_{\Om_\del^i \backslash B(\bx_i, r_i)} \frac{1}{(|\by-\bx|\vee\del)^{d+p-2}(|\by-\bx|\wedge\del)^{\,\beta}} d\by = \del^{-2} \int_{\Om_\del^i \backslash B(\bx_i, r_i)} \frac{1}{|\by-\bx|^{d+p-2}} d\by\\
\leq& \del^{-2} \int_{\left\{ \by\in\R^d:  \psi_i(\overline{\by})<\tilde{y}<\varphi_i(\overline{\by}),\, |\by-\bx|>\epsilon \right\}} \frac{1}{|\by-\bx|^{d+p-2}} d\by \,.
\end{split}
\]
Note that in the last line, we have again adopted the $i$-th local coordinate system. Let $G^i$ be defined in \eqref{eq2:GenInvTraceEst2} and $\bm{w}= (G^i)^{-1}\bx$. Then by the same reasoning as above, 
\[
\begin{split}
&\del^{-2}  \int_{ \left\{ \by\in\R^d:  \psi_i(\overline{\by})<\tilde{y}<\varphi_i(\overline{\by}),\, |\by-\bx|>\epsilon\right\}}\frac{1}{|\yb -\xb|^{d+p-2}}  d\by  \le K_L \del^{-2} \int_{\left\{ \by\in \mcR_{-\del},\, |G^i \by - \bx|>\epsilon \right\}} \frac{1}{|G^i \yb -\xb|^{d+p-2}} d\by  \\
 \le  & \frac{K_L}{(K_L^\prime)^{d+p-2}}\del^{-2} \int_{\left\{ \by\in \mcR_{-\del} -\bm{w},\, |\by|>\epsilon/K_L \right\}} \frac{1}{|\by|^{d+p-2}} d\by. 
 \end{split}
 \]
 Since $\{ \by: |\by|>\epsilon/K_L\}\subset \{ \by=(\tilde{y},\overline{\by}): |\overline{\by}|>\sqrt{3}\epsilon/(2K_L)\} \cup \{ \by=(\tilde{y},\overline{\by}): |\tilde{y}|>\epsilon/(2K_L)\} $, the above quantity can be bounded by
 \[
 \begin{split}
& \frac{K_L}{(K_L^\prime)^{d+p-2}}\del^{-2}\left( \int_{\left\{ \by\in\mcR_{-\del}-\bm{w}, |\tilde{y}|>\epsilon/(2K_L) \right\}}\frac{1}{|\by|^{d+p-2}} d\by + \int_{\left\{ \by\in\mcR_{-\del}-\bm{w}, |\overline{\by}|>\sqrt{3}\epsilon/(2K_L) \right\}}\frac{1}{|\by|^{d+p-2}} d\by \right) \\
 \leq & C\del^{-2}\left( 
\del \int_{ \R^{d-1}}\frac{1}{(|\overline{\by}|+\epsilon/(2 K_L))^{d+p-2}} d\overline{\by} + \del  \int_{|\overline{\by}|>\frac{\sqrt{3}\epsilon}{2K_L}} \frac{1}{|\overline{\by}|^{d+p-2}} d\overline{\by} \right) \leq  C \del^{-1}\,,
\end{split}
\]
where $C$ is a constant that depends on $\epsilon$, $K_L$, $d$ and $p$. 
Combining the estimates, we get
\[
|\lambda_i u|^p_{\cT^{\,\beta}_\del(\Om^i_\del)} \le  C\left( \| u \|^p_{\cT_\del^{\,\beta}(\Om_\del)} + |d+p-\beta|^{-1}\| u\|^p_{L^p(\Om_\del)}\right)\,.
\]
\end{proof}

\begin{lem}\label{lem:GenInvTraceEst1}
Let $\Omega$ be a simply connected Lipschitz domain with an interaction domain $\Omega_\delta$ for $0<\delta< \epsilon$ where $\epsilon$, $\Omega^i$, and $\Omega_\delta^i$ are defined as above. Assume each $\lambda_i$ as above and $E$ defined by \eqref{def:extension_generalLip}, then for any $u \in \cT_\delta^{\,\beta}(\Om_\delta)$,  $\| Eu\|_{\cS^{\,\beta}_\del(\hat\Om)} \le C |d-\beta|^{-1/p} \sum_{i=1}^N\|\lambda_i u (\bx) \|_{\cT^{\,\beta}_\del(\Om^i_\del)}$ where $C$ is independent of $\delta$ and $\beta$.
 \end{lem}
 \begin{proof}
 Notice that $Eu|_{\Om_\del} = u$ since, for $\bx\in \Om_\del$, $\lambda_i(\bx) E^i(\lambda_i u)(\bx) = \lambda_i^2(\bx) u(\bx)$ for all $i$. 
 By \eqref{def:extension_generalLip}, we observe that
\[
\|Eu\|_{\cS^{\,\beta}_\del(\hat\Om)} \leq \sum_{i=1}^N \|\widetilde{\lambda_i} E^i(\lambda_i u)  \|_{\cS^{\,\beta}_\del(\hat\Om)}. 
\]
For each $i\in \{1,\cdots,N\}$, we have 
\[
\begin{split}
 &|\widetilde{\lambda_i} E^i(\lambda_i u)  |^p_{\cS^{\,\beta}_\del(\hat\Om)} \\
 =& 
 \left(\int_{\hat\Omega\cap B(\bx_i,r_i)}\int_{\hat\Omega\cap B(\bx_i,r_i)}+2\int_{\hat\Omega\cap B(\bx_i,r_i)}\int_{\hat\Omega\setminus B(\bx_i,r_i)}\right)\gamma^{\,\beta}_\delta(|\bx-\by|)\left|\left(\widetilde{\lambda_i} E^i(\lambda_i u)\right)(\bx)-\left(\widetilde{\lambda_i} E^i(\lambda_i u)\right)(\by)\right|^pd\by d\bx   \\
 = &| \widetilde{\lambda_i} E^i(\lambda_i u)  |^p_{\cS^{\,\beta}_\del(\Omega^i)} + 2\int_{\hat\Omega\cap B(\bx_i,r_i-\epsilon)}\int_{\hat\Omega\cap\{\by \in\mathbb{R}^d: |\by-\bx|>\epsilon\}}\gamma^{\,\beta}_\delta(|\bx-\by|)\left|\left(\widetilde{\lambda_i} E^i(\lambda_i u)\right)(\bx)-\left(\widetilde{\lambda_i} E^i(\lambda_i u)\right)(\by)\right|^pd\by d\bx \\
 = &| \widetilde{\lambda_i} E^i(\lambda_i u)  |^p_{\cS^{\,\beta}_\del(\Omega^i)} \leq C | E^i(\lambda_i u)  |^p_{\cS^{\,\beta}_\del(\hat\Omega\cap B(\xb_i, r_i))}  \leq C | E^i(\lambda_i u)  |^p_{\cS^{\,\beta}_\del(\Omega^i)} \\
 \leq& C |\beta-d|^{-1} \| \lambda_i u\|^p_{\cT^{\,\beta}_\del(\Om^i_\del)}\,,
\end{split}
\]
where we have used Lemma \ref{lem:GenTraceEst2} and Theorem \ref{thm:invSpecLip}. 
Moreover, since $\widetilde{\lambda_i} \leq 1$, we have
\[
\|\widetilde{\lambda_i} E^i(\lambda_i u)  \|^p_{L^p(\hat\Om)}  = \|\widetilde{\lambda_i} E^i(\lambda_i u)  \|^p_{L^p(\hat\Om\cap B(\bx_i,r_i))}  \leq  \|E^i(\lambda_i u)  \|^p_{L^p(\Om^i)} \leq C \| \lambda_i u\|^p_{\cT^{\,\beta}_\del(\Om^i_\del)}.  
\]
Therefore, we have $\| Eu\|_{\cS^{\,\beta}_\del(\hat\Om)}\leq C |d-\beta|^{-1/p}
 \sum_{i=1}^N\|\lambda_i u (\bx) \|_{\cT^{\,\beta}_\del(\Om^i_\del)} $.
 \end{proof}
 
\noindent\textit{Proof of Theorem 1.3.} 
To show general trace theorem we use Lemmas \ref{lem:GenTraceESt1} and \ref{lem:GenTraceEst2} along with Theorem \ref{thm:speclip} to obtain the estimate
\begin{align*}
\| u\|_{\cT^{\,\beta}_\del(\Om_\del)}&\le C\sum_{i=1}^N\|\lambda_i u\|_{\cT_\delta^{\,\beta}(\Omega^i_\delta)}\le C|d+p-\beta|^{-1/p}\sum_{i=1}^N\|\lambda_i u\|_{\cS_\delta^{\,\beta}(\Omega^i)} \\
&\le C|d+p-\beta|^{-1/p}\sum_{i=1}^{N}\|u\|_{\cS_\delta^{\, \beta}(\hat\Omega\cap B(\bx_i,r_i))} \le  C|d+p-\beta|^{-1/p}\|u\|_{\cS_\delta^{\,\beta}(\hat{\Omega})}\,,
\end{align*}
where $C$ is independent of $\beta$ and $\delta$.  For the general inverse trace theorem, let $E$ be defined by \eqref{def:extension_generalLip}, then from Lemmas \ref{lem:GenInvTraceEst2} and \ref{lem:GenInvTraceEst1} we have 
\[
\| Eu\|_{\cS^{\,\beta}_\del(\hat\Om)}\le C|d-\beta|^{-1/p}\sum_{i=1}^N\|\lambda_i u \|_{\cT^{\,\beta}_\del(\Om^i_\del)} \le C|d-\beta|^{-1/p} \left( \| u \|_{\cT_\del^{\,\beta}(\Om_\del)} + |d+p-\beta|^{-1/p}\| u\|_{L^p(\Om_\del)}\right)  \,,
\]
where $C$ is independent of $\delta$ and $\beta$.
\qed

\section{Conclusion and Discussion}\label{sec:conclusion}


This work gives {suitable} characterizations of the trace spaces of a class of nonlocal function spaces denoted by $\cS^{\,\beta}_\del(\hat\Om)$, where the parameter $\del$ is the nonlocal interaction length and $\beta$ characterizes the singularity of the nonlocal interaction kernels. 
Such nonlocal function spaces have been {extensively} used recently as the energy spaces associated with nonlocal diffusion and nonlocal mechanics models \cite{du2012analysis,silling_2000,mengesha2014nonlocal,DDGG20,you2020data,FR2}.
{However, a clear understanding of the trace spaces of  $\cS^{\,\beta}_\del(\hat\Om)$ has been largely limited \cite{tian2017trace}. }
In the current work, we have introduced the function space $\cT^{\,\beta}_\del(\Om_\del)$ as the trace space of $\cS^{\,\beta}_\del(\hat\Om)$ and demonstrated that 
the trace map from $\cS^{\,\beta}_\del(\hat\Om)$  to $\cT^{\,\beta}_\del(\Om_\del)$  is continuous (given by the trace theorem) and conversely, there is a continuous extension operator from $\cT^{\,\beta}_\del(\Om_\del)$  to $\cS^{\,\beta}_\del(\hat\Om)$ (given by the inverse trace theorem). 
Moreover, the estimates on the trace and the inverse trace maps are uniform {with respect to the horizon parameter $\del$}, so that one can recover the classical trace and inverse trace theorems  {in the local limit} as the nonlocal interaction length $\del\to0$. This is also important {since there are many instances of nonlocal models} that recover the classical diffusion or elasticity equations as $\del\to 0$ \cite{du2012analysis,mengesha2014,mengesha2014nonlocal}. 

The investigation of trace spaces of Sobolev spaces has been a classical research area that has important implications in the mathematical and numerical studies of boundary value problems of local PDEs. {The results of this work therefore are expected to be helpful in the rigorous studies of nonlocal equations with possible nonlocal boundary constraints similar to their PDEs counterparts.  Studies in this direction are currently underway. } 
Moreover, nonlocal functions spaces on vector fields such as those appear in \cite{mengesha2014,mengesha2014nonlocal} can also be studied in the future. Another interesting direction for the future is to investigate, when $\delta\rightarrow \infty$, the consistency of suitably defined nonlocal spaces, similar to those discussed in this work, with their fractional limits.


\section*{Acknowledgments}

Y. Yu is supported by the National Science Foundation under award DMS 1753031.
{Q. Du is supported in part by the National Science Foundation under award DMS-2012562 and the ARO MURI Grant W911NF-15-1-0562.}
X. Tian is supported by the National Science Foundation under award DMS-2111608.

\bibliographystyle{apa}
\bibliography{yyu.bib}

\appendix

\section{Proof of Lemmas in Section \ref{sec:generaldomain}}\label{app:1}


\noindent\textit{Proof of Lemma \ref{lem:psiLip}.} 
First, note that the function $\psi$ is well-defined for each $\overline{\bx} \in \mathbb{R}^{d-1}$ since 
$$\dist((-\delta,\overline{\bx}),\Phi(\mathbb{R}^{d-1})) \ge |(-\delta,\overline{\bx})-(0,\overline{\bx})|=\delta,\quad\text{and }\quad \dist((\varphi(\overline{\bx}),\overline{\bx}),\Phi(\mathbb{R}^{d-1})) = 0,$$
there must exist at least one $\tilde{z} \in [-\delta,\varphi(\overline{\bx}))$ such that $\dist((\tilde{z},\overline{\bx}),\Phi(\mathbb{R}^{d-1}))=\delta$ by the continuity of $\tilde{x} \mapsto\dist((\tilde{x},\overline{\bx}),\Phi(\mathbb{R}^{d-1}))$.  The minimum of such $\tilde{z}$ is the value of $\psi(\overline{\bx})$ and clearly $\psi(\overline{\bx})< \varphi(\overline{\bx})$ by how the $\tilde{z}$ were chosen.

Let $\overline{\bx}, \overline{\by} \in \mathbb{R}^{d-1}, \ \overline{\bx} \neq \overline{\by}$ and without loss of generality assume $\psi(\overline{\bx}) \le \psi(\overline{\by})$.  By definition of $\psi$, there is a $\hat{\bx} \in \mathbb{R}^{d-1}$ such that $\dist((\psi(\overline{\bx}),\overline{\bx}),(\varphi(\hat{\bx}),\hat{\bx}))=\delta$.  Now define
\[
c_1:=\varphi(\hat{\bx})-\psi(\overline{\bx}), \text{ and } \overline{\bm c}:=\hat{\bx}-\overline{\bx}.
\]
Notice that we must have $\psi(\overline{\by})+c_1\le\varphi(\overline{\by}+\overline{\bm c})$.  If this was not the case, then consider function 
\[
\bm F(t)=(t\psi(\overline{\by})+(1-t)(\psi(\overline{\by})+c_1),t\overline{\by}+(1-t)(\overline{\by}+\overline{\bm c}))
\]for $t \in [0,1]$.  By the assumption, we have
\begin{align*}
\bm F(0)&=(\psi(\overline{\by})+c_1,\overline{\by}+\overline{\bm c}) \in \{\bx=(\tilde{x},\overline{\bx}): \tilde{x} > \varphi(\overline{\bx})\};\\ \bm F(1)&=(\psi(\overline{\by}),\overline{\by}) \in \{\bx=(\tilde{x},\overline{\bx}): \tilde{x} < \varphi(\overline{\bx})\}.
\end{align*}
Since $\bm F$ is continuous, there must be a $t_0 \in (0,1)$ such that $\bm F(t_0)=(\varphi(\overline{\bm z}),\overline{\bm z})$ for some $\overline{\bm z} \in \mathbb{R}^{d-1}$. Then,
\begin{align*}
|(\psi(\overline{\by}),\overline{\by})-(\varphi(\overline{\bz}),\overline{\bz})|&= |(\psi(\overline{\by}),\overline{\by})-(t_0\psi(\overline{\by})+(1-t_0)(\psi(\overline{\by})+c_1),t_0\overline{\by}+(1-t_0)(\overline{\by}+\overline{\bm c}))|\\
&=|((1-t_0)c_1,(1-t_0)\overline{\bm c})|< |(c_1,\overline{\bm{c}})|= \delta.
\end{align*}
which is a contradiction of the definition of $\psi$. With this inequality we finally have
\begin{align*}
    \frac{|\psi(\overline{\by})-\psi(\overline{\bx})|}{|\overline{\by}-\overline{\bx}|}= \frac{\psi(\overline{\by})-\psi(\overline{\bx})}{|\overline{\by}-\overline{\bx}|}= \frac{\psi(\overline{\by})+c_1-\varphi(\hat{\bx})}{|\overline{\by}+\overline{\bm c}-\hat{\bx}|} \le \frac{\varphi(\overline{\by}+\overline{\bm c})-\varphi(\hat{\bx})}{|\overline{\by}+\overline{\bm c}-\hat{\bx}|} < L.
\end{align*}
\qed

\begin{figure}
\centering
\begin{tikzpicture}
\draw (3,0) to [out=0, in=225] (6,3);
\draw (3,0) to [out=180, in=-45] (0,3);
\draw (3,1) to [out=0, in=225] (5,3);
\draw (3,1) to [out=180, in=-45] (1,3);
\draw [fill] (2,1.55) circle [radius=0.05];
\draw [fill] (2,0.3) circle [radius=0.05];
\draw [fill] (1.3,1.1) circle [radius=0.05];
\node [above] at (2,1.55) {$\bm{x}$};
\node [left] at (1.3,1.1) {$\bm{x}_n$};
\node [right] at (2,0.3) {$\bm{y}$};
\node at (3.5,.5) {$\Omega_\delta$};
\node at (3,2) {$\Omega$};
\node [right] at (2,.95) {$h$};
\node [above left] at (1.7,1.3) {$\delta$};
\draw [thick] (2,1.55)--(2,0.3);
\draw [thick] (1.3,1.1)--(2,0.3);
\draw [thick] (2,1.55)--(1.3,1.1);
\draw [dashed] (1.3,1.1)--(1.3,0.3);
\draw [densely dashed] (1.3,0.3)--(2,0.3);
\node [left] at (1.3,0.6) {\scriptsize $|\psi(\bm{z}_n)-\psi(\bm{z})|$};
\node [below] at (1.65,0.3) {\scriptsize $|\bm{z}_n-\bm{z}|$};
\draw (2,.6) arc (90:130:.3);
\draw (1.45,0.97) arc (-40:30:.2);
\node [left] at (2.05,0.7) {\scriptsize $\theta$};
\node [right] at (1.4,1.1) {\scriptsize $\gamma$};
\end{tikzpicture}
\caption{Settings and definitions for the proof of Lemma \ref{lem:phipsisupnorm}.}\label{fig:lemma}
\end{figure}
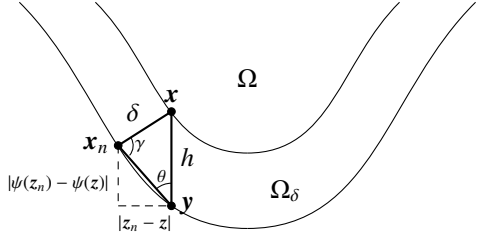

\noindent\textit{Proof of Lemma \ref{lem:phipsisupnorm}.} 
Let $\bm{z} \in \mathbb{R}^{d-1}$ and define $\bm{x}:=(\varphi(\bm{z}),\bm{z})$ and $\bm{y}:=(\psi(\bm{z}),\bm{z})$. By the definition of $\psi$, there is a $\bm{z}_n \in \mathbb{R}^{d-1}$ such that $|\bm{x}-\bm{x}_n|=\delta$ where $\bm{x}_n:=(\psi(\bm{z}_n),\bm{z}_n)$. Note that $\bm{x}=(h+\psi(\bm{z}),\bm{z})$ where $h=\varphi(\bm{z})-\psi(\bm{z})$. If $\bm{x}_n=\bm{y}$, then we are done since then 
\[
h=|(h,0)|=|\bm{x}-\bm{y}|=|\bm{x}-\bm{x}_n|=\delta \le \delta\sqrt{L^2+1}.
\] 
Otherwise, consider the triangle made by the points $\bm{x}_n$, $\bm{y}$, and $\bm{x}$ as shown in Figure \ref{fig:lemma}. \footnote{Here we are viewing the triangle by taking the natural isometry to $\mathbb{R}^2$ from an orthonormal basis for span$(\bm{y}-\bm{x},\bm{y}-\bm{x_n})$.  Note that $\bm{y}-\bm{x}$ and $\bm{y}-\bm{x}_n$ are linearly independent under the assumption $\bm{y} \neq \bm{x}_n$.  If this was not the case, then
\[
h\sqrt{|\varphi(\overline{\bm{y}})-\varphi(\bx_n)|^2+|\overline{\bm{y}}-\overline{x_n}|^2}=|\bm{y}-\bm{x}||\bm{y}-\bm{x}_n|=|(\bm{y}-\bm{x})\cdot(\bm{y}-\bm{x}_n)|=h|\varphi(\overline{\bm{y}})-\varphi(\overline{\bx_n})|
\]
requiring $\overline{\bm{y}}=\overline{\bm{x}_n}$. Since $\varphi$ is a single-valued function, $\varphi(\overline{\bm{y}})=\varphi(\overline{\bx_n})$, hence $\bm{y}=\bm{x}_n$, giving a contradiction.}  
Let $\theta$ be the angle associated with the vertex $\bm{y}$ and note that $\psi$ is a Lipschitz function with Lipschitz constant $L$. Then,
\[
0 <  \frac{\pi}{2}-\arctan(L) \le \frac{\pi}{2}-\arctan\left(\frac{|\psi(\bm{z})-\psi(\bm{z}_n)|}{|\bm{z}-\bm{z}_n|}\right)=\theta\le \frac{\pi}{2}.
\]
So if $\gamma$ is the angle associated with the vertex $\bm{x}_n$, then by the law of sines
\[
\frac{1}{\delta\sqrt{L^2+1}} = \frac{\cos(\arctan(L))}{\delta} = \frac{\sin\left(\frac{\pi}{2}-\arctan(L)\right)}{\delta}\le\frac{\sin(\theta)}{\delta}=\frac{\sin(\gamma)}{h} \le \frac{1}{h}.
\]
and so $h \le \delta\sqrt{L^2+1}$. By the definition of $\psi$ we clearly have $\delta \le h$.  Since $\bm z \in \mathbb{R}^{d-1}$ was arbitrary we finish the proof.

\qed



\noindent\textit{Proof of Lemma \ref{lem:projectiondist}.} 
By Lemma \ref{lem:psiLip} $\varphi-\psi$ is a Lipschitz function with Lipschitz constant $C_L$.  Then in combination with Lemma \ref{lem:phipsisupnorm},
\begin{align*}
    |\bx'-\by'|&\le|\overline{\bx}-\overline{\by}|+|\tilde{x}'-\tilde{y}'|\le |\overline{\bx}-\overline{\by}|+\left|\varphi(\overline{\bx})-\varphi(\overline{\by})\right|+\left|\frac{\tilde{x}}{\delta}(\varphi(\overline{\bx})-\psi(\overline{\bx}))-\frac{\tilde{y}}{\delta}(\varphi(\overline{\by})-\psi(\overline{\by}))\right|\\
    &\le |\overline{\bx}-\overline{\by}|+L|\overline{\bx}-\overline{\by}|+\left|\frac{\tilde{x}-\tilde{y}}{\delta}(\varphi(\overline{\bx})-\psi(\overline{\bx}))\right|+\left|\frac{\tilde{y}}{\delta}(\varphi(\overline{\by})-\psi(\overline{\by})-(\varphi(\overline{\bx})-\psi(\overline{\bx}))\right|\\
    &\le (L+1)|\overline{\bx}-\overline{\by}|+\sqrt{L^2+1}|\tilde{x}-\tilde{y}|+2\sqrt{L^2+1}|\overline{\bx}-\overline{\by}|\le K_L|\bx-\by|.
\end{align*}
For left inequality, first notice from \eqref{eq:changeofvariable_x} that 
\begin{align*}
    \tilde{x}&=\frac{\del(\tilde{x}'-\varphi(\overline{\bx}))}{\varphi(\overline{\bx})-\psi(\overline{\bx})},\quad \text{and}\quad
     \tilde{y}=\frac{\del(\tilde{y}'-\varphi(\overline{\by}))}{\varphi(\overline{\by})-\psi(\overline{\by})}. 
\end{align*}
Additionally, $-\delta<\tilde{y}<0$ implies $\psi(\overline{\by})<\tilde{y}'<\varphi(\overline{\by})$ so then,
\begin{align*}
    |\bx-\by| &\le |\tilde{x}-\tilde{y}| + |\overline{\bx}-\overline{\by}|\\
    &\leq\delta\frac{|\varphi(\overline{\bx})-\varphi(\overline{\by})-(\tilde{x}'-\tilde{y}')|}{|\varphi(\overline{\bx})-\psi(\overline{\bx})|}+\delta|\varphi(\overline{\by})-\tilde{y}'|\left|\frac{\varphi(\overline{\by})-\psi(\overline{\by})-(\varphi(\overline{\bx})-\psi(\overline{\bx}))}{(\varphi(\overline{\bx})-\psi(\overline{\bx}))(\varphi(\overline{\by})-\psi(\overline{\by}))}\right|+|\overline{\bx}-\overline{\by}|\\
    &\le \delta\frac{L|\overline{\bx}-\overline{\by}|+|\tilde{x}'-\tilde{y}'|}{\inf |\varphi-\psi|}+\delta \frac{|\tilde{y}(\varphi(\overline{\by})-\psi(\overline{\by}))|}{\delta}\left|\frac{\varphi(\overline{\by})-\varphi(\overline{\bx})+\psi(\overline{\bx}) -\psi(\overline{\by})}{(\varphi(\overline{\bx})-\psi(\overline{\bx}))(\varphi(\overline{\by})-\psi(\overline{\by}))}\right|+|\overline{\bx}-\overline{\by}|\\
    & \le \delta\frac{L|\overline{\bx}-\overline{\by}|+|\tilde{x}'-\tilde{y}'|}{\inf |\varphi-\psi|}+\delta \frac{2L|\overline{\bx}-\overline{\by}|}{\inf |\varphi-\psi|}+|\overline{\bx}-\overline{\by}|\le K_L'^{-1}|\bx'-\by'|.
    \end{align*}
\qed


\noindent\textit{Proof of Lemma \ref{lem:KernelEstimate}.} 
For (a) if $|\by-\bx|<\del/M$, then
\begin{align*}
|\bm w-\bz|&=|(\tilde{y}+\varphi(\overline{\bm{y}}),\overline{\bm{y}})-(\tilde{x}+\varphi(\overline{\bm{x}}),\overline{\bm{x}})|\leq|\tilde{y}-\tilde{x}|+|\varphi(\overline{\bm{y}})-\varphi(\overline{\bm{x}})|+|\overline{\bm{y}}-\overline{\bm{x}}|\\
&\le |\tilde{y}-\tilde{x}|+(L+1)|\overline{\bm{y}}-\overline{\bm{x}}|\le (L+2)|\bm{y}-\bm{x}|< M\frac{\delta}{M}=\delta.
\end{align*}
Therefore we have
\[
\gamma^{\,\beta}_{\delta/M}(|\bm{y}-\bm{x}|) \le \frac{1}{(\delta/M)^{d+p-\beta}}\frac{1}{|\bm{y}-\bm{x}|^{\,\beta}}\bm{1}_{|\bm{y}-\bm{x}|<\delta/M}\leq \frac{M^{d+p-\beta}}{\delta^{d+p-\beta}}\frac{M^\beta}{|\bm{w}-\bm{z}|^{\,\beta}}\bm{1}_{|\bm{w}-\bm{z}|<\delta} \leq M^{d+p}\gamma^{\,\beta}_{\delta}(|\bm w-\bz|).
\]
For (b),  with $|\by-\bx|<\del/M$ and noting $|\tilde{x}|<\delta$ we have
\begin{equation*}\label{eqn:mixeddomaineq}
\begin{split}
 &|\bm w - \bx'|\leq\left| \tilde{y}+\varphi(\overline{\bm{y}})-\left[\left(1+\frac{\tilde{x}}{\delta}\right)\varphi(\overline{\bm{x}})-\frac{\tilde{x}}{\delta}\psi(\overline{\bm{x}})\right] \right| + |\overline{\by}-\overline{\bx}| =\left| \varphi(\overline{\bm{y}})-\varphi(\overline{\bm{x}})+\tilde{y}-\frac{\tilde{x}}{\delta}(\varphi(\overline{\bm{x}})-\psi(\overline{\bm{x}})) \right| + |\overline{\by}-\overline{\bx}|  \\ 
 &\leq \left| \varphi(\overline{\bm{y}})-\varphi(\overline{\bm{x}})\right| +\left|\tilde{y}-\frac{\tilde{x}}{\delta}(\varphi(\overline{\bm{x}})-\psi(\overline{\bm{x}})) \right|+|\overline{\by}-\overline{\bx}|\leq L |\overline{\by}-\overline{\bx}| + \sqrt{L^2+1} |\tilde{y}-\tilde{x}| + |\overline{\by}-\overline{\bx}| \\
 &\leq (L+1+\sqrt{L^2+1}) |\bm{y}-\bm{x}| < \del.
\end{split}
\end{equation*}
By the same reason, we have 
\[
\gamma_{\delta/M}^\beta(|\bx-\by|) \le M^{d+p}\gamma_{\delta}^\beta(|\bm x' - \bm w|).
\]
We similarly have (c) by noticing that $|\by'-\bx'|\leq K_L |\by-\bx|\leq M |\by-\bx|$.
\qed

\noindent\textit{Proof of Lemma \ref{lem:InverseKernelEstimate}.} 
For (a), 
\[
|\bm w -\bm z|\leq | \tilde{y}'-\tilde{x}'|+ |\varphi(\overline{\by}) - \varphi(\overline{\bx})|+ |\overline{\by}-\overline{\bx}|\leq (L+2)|\by'-\bx'| \le M\delta
\]
and so 
\[
\gamma^{\,\beta}_{\delta}(|\bm{y}'-\bm{x}'|) \le \frac{1}{\delta^{d+p-\beta}}\frac{1}{|\bm{y}'-\bm{x}'|^{\,\beta}}\bm{1}_{|\bm{y}'-\bm{x}'|<\delta}\leq \frac{M^{d+p-\beta}}{(M\delta)^{d+p-\beta}}\frac{M^\beta}{|\bm{w}-\bm{z}|^{\,\beta}}\bm{1}_{|\bm{w}-\bm{z}|<M\delta} \leq M^{d+p}\gamma^{\,\beta}_{M\delta}(|\bm w-\bz|).
\]
For (b), 
\begin{align*}
    |\bm w-\bx|&\leq \left|\tilde{y}'-\varphi(\overline{\by})+\delta\left(\dfrac{\varphi(\overline{\bx})-\tilde{x}'}{\varphi(\overline{\bx})-\psi(\overline{\bx})}\right)\right| + |\overline{\by}-\overline{\bx}|\le |\tilde{y}'-\varphi(\overline{\by})+\varphi(\overline{\bx})-\tilde{x}'|+|\overline{\by}-\overline{\bx}|\le (L+2)|\by'-\bx'|
\end{align*}
and so similarly we have $\gamma_{\delta}^\beta(|\bx'-\by'|) \le M^{d+p}\gamma_{
M\delta}^\beta(|\bm x - \bm w|)$. \newline
For (c) note  $K_L'|\by-\bx|\leq  |\by'-\bx'|$ for $\bx'\in \Om_\del$ and $\by'\in \Om_\del$ where $K_L'$ is as in Lemma \ref{lem:projectiondist} and we can conclude the result in a similar manner.
\qed

\end{document}